\documentclass[twoside]{amsart}
\usepackage{amsmath}
\usepackage{amssymb}
\usepackage[cp850]{inputenc}
\usepackage[bookmarksnumbered,plainpages,hypertex]{hyperref}

\input{xy}
 \xyoption{all}
\newtheorem{theorem}{\sc Theorem}[section]
\newtheorem{proposition}[theorem]{\sc Proposition}
\newtheorem{notation}[theorem]{\sc Notation}
\newtheorem{lemma}[theorem]{\sc Lemma}
\newtheorem{corollary}[theorem]{\sc Corollary}
\theoremstyle{definition}
\newtheorem{definition}[theorem]{\sc Definition}
\newtheorem{definitions}[theorem]{\sc Definitions}

\theoremstyle{remark}
\newtheorem{remark}[theorem]{\sc Remark}

\newtheorem{claim}[theorem]{}

\def\ot{\otimes}
\def\cot{\square}

\def\w{\wedge }

\newcommand{\diagUnivTens}{\xymatrix@R=20pt@C=60pt{
  M \ar[d]_{f_M} \ar[r]^{i_1} & T_A(M) \ar@{.>}[dl]|f \\
  B  & A\ar[l]^{f_A} \ar[u]_{i_0 }}}
\newcommand{\diagUnivCot}{\xymatrix@C=2cm{
  T_{C}^{c}(M) \ar[d]_{p_0} \ar[r]^{p_1} & M   \\
  C& \widetilde{D}  \ar[u]_{f_M}\ar@{.>}[ul]|-{f}\ar[l]^{f_C} &D\ar[l]^{\xi_1}  \ar[ul]_0        }  }
\newcommand{\diagMultil}{\xymatrix@C=2cm{
  M \ar[d]_{\left( i_{0}\otimes i_{1}\right) \rho _{M}^{l}+\left(
i_{1}\otimes i_{0}\right) \rho _{M}^{r}} \ar[r]^{i_1} & T
\ar@{.>}[dl]|-{\Delta_T}  \\
  T\ot T& H  \ar[u]_{i_0} \ar[l]^{\left( i_{0}\otimes i_{0}\right) \Delta _{H}}         }                }
\newcommand{\diagMultir}{\xymatrix@C=2cm{
  M \ar[d]_{0} \ar[r]^{i_1} & T
\ar@{.>}[dl]|-{\varepsilon_T}  \\
  \mathbf{1}& H  \ar[u]_{i_0} \ar[l]^{\varepsilon _{H}}        }                 }
\newcommand{\diagUnivTensBialg}{\xymatrix@R=20pt@C=60pt{
  M \ar[d]_{f_M} \ar[r]^{i_1} & T_H(M) \ar@{.>}[dl]|f \\
  E  & H\ar[l]^{f_H} \ar[u]_{i_0 }}}
\newcommand{\diagComultil}{ \xymatrix@C=2cm{
  T^c \ar[d]_{p_0} \ar[r]^{p_1} & M   \\
  H& T^c\ot T^c  \ar[u]_{\mu _{M}^{l}\left( p_{0}\otimes p_{1}\right) +\mu
_{M}^{r}\left( p_{1}\otimes
p_{0}\right)}\ar@{.>}[ul]|-{m_{T^c}}\ar[l]^{m_{H}\left(
p_{0}\otimes p_{0}\right)}          }}
\newcommand{\diagComultir}{ \xymatrix@C=2cm{
  T^c \ar[d]_{p_0} \ar[r]^{p_1} & M   \\
  H& \mathbf{1}  \ar[u]_{0}\ar@{.>}[ul]|-{u_{T^c}}\ar[l]^{u_H}   }       }
\newcommand{\diagUnivCotensorBialg}{ \xymatrix@R=20pt@C=2cm{
  T_{H}^{c}(M) \ar[d]_{p_0} \ar[r]^{p_1} & M   \\
  H& \widetilde{D}  \ar[u]_{f_M}\ar@{.>}[ul]|-{f}\ar[l]^{f_H} &D\ar[l]^{\xi_1}  \ar[ul]_0    }      }
\newcommand{\diagFofTypeOne}{\xymatrix@R=20pt@C=60pt{
  M \ar[d]_{i^{T^c}_{1}} \ar[r]^{i^T_1} & T_H(M) \ar@{.>}[dl]|F \\
  T_{H}^c(M)  & H\ar[l]^{i^{T^c}_{0}} \ar[u]_{i^T_0 }}}
\newcommand{\diagTau}{\xymatrix@R=40pt@C=80pt{
                &         \textbf{1}^{\wedge _{E}^{n+1}}\ar@{.>}[dl]_{\tau_n} \ar[d]^{\alpha_{\textbf{1}^{\wedge _{E}^{n+1}}}}     \\
  {{\textbf{1}^{\wedge _{E}^{n}}}}\otimes {\textbf{1}^{\wedge _{E}^{n}}} \ar[r]_{{\xi_n^{n+1}}\otimes {\xi_n^{n+1}}} &
{\textbf{1}^{\wedge _{E}^{n+1}}}\otimes {\textbf{1}^{\wedge
_{E}^{n+1}}}}}
\newcommand{\diagTauA}{\xymatrix@C=0.5cm{
  0 \ar[r] & {\textbf{1}^{n}} \ar[rr]^{{\xi
_{n}^{n+1}}} && {\textbf{1}^{n+1}}
\ar[rr]^{p_{\mathbf{1}^{n}}^{\mathbf{1}^{n+1}}} &&
\frac{\textbf{1}^{n+1}}{\textbf{1}^{n}} \ar[r] & 0 }}
\newcommand{\diagTauB}{\xymatrix@C=1cm{
  0 \ar[r] & {\textbf{1}^{n+1}}\ot {\textbf{1}^{n}} \ar[rr]^{{\textbf{1}^{n+1}}\ot {\xi
_{n}^{n+1}}} &&{\textbf{1}^{n+1}}\ot  {\textbf{1}^{n+1}}
\ar[rr]^{{\textbf{1}^{n+1}}\ot
p_{\mathbf{1}^{n}}^{\mathbf{1}^{n+1}}} && {\textbf{1}^{n+1}}\ot
\frac{\textbf{1}^{n+1}}{\textbf{1}^{n}}
\\&&&\textbf{1}^{n+1}\ar@{.>}[ull]^{\beta_n}\ar[u]_{\alpha_{\textbf{1}^{n}}}}}
\newcommand{\diagTauC}{\xymatrix@C=1cm{
  0 \ar[r] & {\textbf{1}^{n}} \ot {\textbf{1}^{n}}\ar[rr]^{{\xi
_{n}^{n+1}} \ot {\textbf{1}^{n}}} && {\textbf{1}^{n+1}} \ot
{\textbf{1}^{n}} \ar[rr]^{p_{\mathbf{1}^{n}}^{\mathbf{1}^{n+1}}
\ot {\textbf{1}^{n}}} && \frac{\textbf{1}^{n+1}}{\textbf{1}^{n}}
\ot {\textbf{1}^{n}}
\\&&&\textbf{1}^{n+1}\ar@{.>}[ull]^{\tau_n}\ar[u]_{\beta_n}}}
\newcommand{\diagWedge}{\xymatrix@C=80pt{
  F\w_E G \ar[d]_{p_G^F} \ar[r]^{\Delta_{F\w_E G}} & (F\w_E G)\ot (F\w_E G) \ar[dd]^{(F\w_E G)\ot p_G^F} \\
  \frac{F\w_E G}{G}\ar[d]_{^F\rho_{\frac{F\w_E G}{G}}}  \\
  F\ot  \frac{F\w_E G}{G}\ar[r]^{i_F^{F\w_E G}\ot \frac{F\w_E G}{G} } & (F\w_E G)\ot\frac{F\w_E G}{G}   }}
\newcommand{\diagUnivTensBialgA}{\xymatrix@R=30pt@C=60pt{
  V\ot H \ar[d]_{i_1^{T_{\mathbf{1}}(V)}\ot H} \ar[r]^{i^{T_H(V\ot H)}_1} & T_H(V\ot H) \ar@{.>}[dl]|f \\
  T_{\mathbf{1}}(V)\rtimes H  &  H\ar[l]^{(i_0^{T_{\mathbf{1}}(V)}\ot H)l_H^{-1}} \ar[u]_{i^{T_H(V\ot H)}_0 }}}

\newcommand{\diagMultIdeal}{\xymatrix{ 0\ar[r]& IJ \ar[rr]^{i_{IJ}} &&  A \ar[rr]^{\pi_{I,J}} &&Q_{I,J}\ar[r]&0  \\ && I\ot J  \ar@{.>}[lu]^{\overline{m}_{I,J}}\ar[ur]_{m_{I,J}}}}
\newcommand{\diagComultWedge}{\xymatrix{ 0\ar[r]& C\w_E  D \ar[rr]^{i_{C\w_E  D}^E} &&  E\ar[dr]_{\Delta_{C,D}} \ar[rr]^{p^E_{C\w_E  D}} &&\frac{E}{C\w_E  D}\ar[r]\ar@{.>}[dl]^{\overline{\Delta}_{C,D}}&0  \\
 &&&& \frac{E}{C}\ot\frac{E}{D}}}

\begin{document}
\title{Braided Bialgebras of Type One}
\author{A. Ardizzoni}
\address{University of Ferrara, Department of Mathematics, Via Machiavelli
35, Ferrara, I-44100, Italy}
\email{rdzlsn@unife.it}
\urladdr{http://www.unife.it/utenti/alessandro.ardizzoni}
\author{C. Menini}
\address{University of Ferrara, Department of Mathematics, Via Machiavelli
35, Ferrara, I-44100, Italy}
\email{men@unife.it}
\urladdr{http://www.unife.it/utenti/claudia.menini}
\thanks{This paper was written while both authors were members of
G.N.S.A.G.A. with partial financial support from M.I.U.R..}

\begin{abstract}
Braided bialgebras of type one in abelian braided monoidal categories are
characterized as braided graded bialgebras which are strongly $\mathbb{N} $%
-graded both as an algebra and as a coalgebra.
\end{abstract}

\maketitle
\tableofcontents

\pagestyle{headings}

\section*{Introduction}

\markboth{\sc{A. ARDIZZONI AND C. MENINI}} {\sc{Bialgebras of Type
One in Monoidal Categories}}Bialgebras of type one were introduced by
Nichols in \cite{Ni}. They came out to play a relevant role in the theory of
Hopf Algebras. In particular, their "coinvariant" part, called Nichols
algebra, has been deeply investigated, see e.g. \cite{Ro}, \cite{AS} and the
references therein.

The main purpose of this paper is to study some relevant properties of
bialgebras of type one in the framework of abelian braided monoidal
categories.

Let $H$ be a braided bialgebra in a cocomplete and complete abelian braided
monoidal category $\left( \mathcal{M},c\right) $ satisfying $AB5$. Assume
that the tensor product commutes with direct sums and is two-sided exact.
Let $M$ be in $_{H}^{H}\mathcal{M}_{H}^{H}$. Let $T=T_{H}(M)$ be the
relative tensor algebra and let $T^{c}=T_{H}^{c}(M)$ be the relative
cotensor coalgebra as introduced in \cite{AMS:Cotensor}. We prove that both $%
T$ and $T^{c}$ have a natural structure of graded braided bialgebra and that
the natural algebra morphism $F:T\rightarrow T^{c}$, which coincide with the
canonical injections on $H$ and $M$, is a graded bialgebra homomorphism.
Thus its image is a graded braided bialgebra which we denote by $H[M]$ and
call, accordingly to \cite{Ni}, \emph{the braided bialgebra of type one
associated to }$H$\emph{\ and }$M$.

We would like to outline that a treatment of bialgebras of type one in the
case $H=\mathbf{1}$ can be found in \cite{Schauenburg1}. Following a
different approach, that traces back to some Nichols' original ideas, in
order to characterize braided bialgebras of type one, we develop some
relevant properties of graded algebras and coalgebras. We introduce the
definition of strongly $%
\mathbb{N}
$-graded algebras and coalgebras which is inspired in the second case to
\cite{NT}.

Given a graded coalgebra $\left( C=\oplus _{n\in
\mathbb{N}
}C_{n},\Delta ,\varepsilon \right) $ in $\mathcal{M}$, we can write $\Delta
_{\mid C_{n}}$ as the sum of unique components $\Delta
_{i,j}:C_{i+j}\rightarrow C_{i}\otimes C_{j}$ where $i+j=n.$ The coalgebra $%
C $ is defined to be a \emph{strongly }$%
\mathbb{N}
$\emph{-graded coalgebra} when $\Delta _{i,j}^{C}:C_{i+j}\rightarrow
C_{i}\otimes C_{j}$ is a monomorphism for every $i,j\in \mathbb{N}$. In
Theorem \ref{teo: psi mono}, we prove that the following assertions are
equivalent.

$\left( a\right) $ $C$ is a strongly $%
\mathbb{N}
$-graded coalgebra.

$\left( b\right) $ the canonical morphism $\psi :C\rightarrow
T_{C_{0}}^{c}(C_{1})$ is a monomorphism.

$\left( c\right) $ $C_{0}\oplus C_{1}=C_{0}\wedge _{C}C_{0}.$

Also the "dual" notion of\emph{\ strongly }$%
\mathbb{N}
$\emph{-graded algebra} is introduced and a similar result is achieved.

Braided bialgebras of type one are then characterized as braided graded
bialgebras which are strongly $%
\mathbb{N}
$-graded both as an algebra and as a coalgebra.

Assume that $H$ is a braided Hopf algebra in $\left( \mathcal{M},c\right) $
with bijective antipode. Using general results by Bespalov and Drabant (see
\cite{Bes} and \cite{BD}), we prove that $H[M]$ is the Radford-Majid
bosonization
\begin{equation*}
H[M]=\mathbf{1}\left[ M^{co\left( H\right) }\right] \rtimes H
\end{equation*}%
where $\mathbf{1}\left[ M^{co\left( H\right) }\right] $ is the braided
bialgebra of type one associated to $\mathbf{1}$ and $M^{co\left( H\right)
}\ $in $_{H}^{H}\mathcal{YD}\left( \mathcal{M}\right) $.

Finally we would like to outline that many results are firstly stated and
proved in the coalgebra case. Our choice is motivated by the special feature
of this case, due to the lack of the existence of the coradical in the
monoidal context. As far as the algebra case is concerned, we would like to
point out that here we can drop the assumptions of completeness and $AB5$%
-condition essentially because the analogue of Lemma \ref{lem: complete}
doesn't require these extra conditions. Also proofs in the algebra case are
not given whenever they would have been an easy adaptation of the coalgebra
ones.

Elsewhere we will apply the results of the present paper to study strictly
graded bialgebras and graded bialgebras associated to a (co)algebra
homomorphism.

\section{Preliminaries and Notations\label{sec: Preliminaries and Notations}}

Let $[(X,i_{X})]$ be a subobject of an object $E$ in an abelian category $%
\mathcal{C},$ where $i_{X}=i_{X}^{E}:X\hookrightarrow E$ is a monomorphism
and $[(X,i_{X})]$ is the associated equivalence class. By abuse of language,
we will say that $(X,i_{X})$ is a subobject of $E$ and we will write $%
(X,i_{X})=(Y,i_{Y})$ to mean that $(Y,i_{Y}) \in [(X,i_{X})]$. The same
convention applies to cokernels. If $(X,i_{X})$ is a subobject of $E$ then
we will write $(E/X,p_X)=\text{Coker}(i_X)$, where $p_{X}=p_{X}^{E}:E%
\rightarrow E/X$.

\begin{claim}
\textbf{Monoidal Categories.} Recall that (see \cite[Chap. XI]{Kassel}) a
\emph{monoidal category}\textbf{\ }is a category $\mathcal{M}$ endowed with
an object $\mathbf{1}\in \mathcal{M}$\textbf{\ } (called \emph{unit}), a
functor $\otimes :\mathcal{M}\times \mathcal{M}\rightarrow \mathcal{M}$
(called \emph{tensor product}), and functorial isomorphisms $%
a_{X,Y,Z}:(X\otimes Y)\otimes Z\rightarrow $ $X\otimes (Y\otimes Z)$, $l_{X}:%
\mathbf{1}\otimes X\rightarrow X,$ $r_{X}:X\otimes \mathbf{1}\rightarrow X,$
for every $X,Y,Z$ in $\mathcal{M}$. The functorial morphism $a$ is called
the \emph{associativity constraint }and\emph{\ }satisfies the \emph{Pentagon
Axiom, }that is the following relation
\begin{equation*}
(U\otimes a_{V,W,X})\circ a_{U,V\otimes W,X}\circ (a_{U,V,W}\otimes
X)=a_{U,V,W\otimes X}\circ a_{U\otimes V,W,X}
\end{equation*}%
holds true, for every $U,V,W,X$ in $\mathcal{M}.$ The morphisms $l$ and $r$
are called the \emph{unit constraints} and they obey the \emph{Triangle
Axiom, }that is $(V\otimes l_{W})\circ a_{V,\mathbf{1},W}=r_{V}\otimes W$,
for every $V,W$ in $\mathcal{M}$.
\end{claim}

It is well known that the Pentagon Axiom completely solves the consistency
problem arising out of the possibility of going from $((U\otimes V)\otimes
W)\otimes X$ to $U\otimes (V\otimes (W\otimes X))$ in two different ways
(see \cite[page 420]{Majid}). This allows the notation $X_{1}\otimes \cdots
\otimes X_{n}$ forgetting the brackets for any object obtained from $%
X_{1},\cdots X_{n}$ using $\otimes $. Also, as a consequence of the
coherence theorem, the constraints take care of themselves and can then be
omitted in any computation involving morphisms in $\mathcal{M}$.\newline
Thus, for sake of simplicity, from now on, we will omit the associativity
constraints.\medskip \newline
The notions of algebra, module over an algebra, coalgebra and comodule over
a coalgebra can be introduced in the general setting of monoidal categories.
Given an algebra $A$ in $\mathcal{M}$ on can define the categories $_{A}%
\mathcal{M}$, $\mathcal{M}_{A}$ and $_{A}\mathcal{M}_{A}$ of left, right and
two-sided modules over $A$ respectively. Similarly, given a coalgebra $C$ in
$\mathcal{M}$, one can define the categories of $C$-comodules $^{C}\mathcal{M%
},\mathcal{M}^{C},{^{C}\mathcal{M}^{C}}$. For more details, the reader is
refereed to \cite{AMS}.

\begin{definitions}
\label{abelian assumptions}Let $\mathcal{M}$ be a monoidal category.\newline
We say that $\mathcal{M}$ is an \textbf{abelian monoidal category }if $%
\mathcal{M}$ is abelian and both the functors $X\otimes (-):\mathcal{M}%
\rightarrow \mathcal{M}$ and $(-)\otimes X:\mathcal{M}\rightarrow \mathcal{M}
$ are additive and right exact, for any $X\in \mathcal{M}.$\newline
We say that $\mathcal{M}$ is an \textbf{coabelian monoidal category }if $%
\mathcal{M}^o$ is an abelian monoidal category, where $\mathcal{M}^{o}$
denotes the dual monoidal category of $\mathcal{M}$. Recall that $\mathcal{M}%
^{o}$ and $\mathcal{M}$ have the same objects but $\mathcal{M}^{o}(X,Y)=%
\mathcal{M}(Y,X) $ for any $X,Y$ in $\mathcal{M}$.
\end{definitions}

Given an algebra $A$ in an abelian monoidal category $\mathcal{M}$, there
exist a suitable functor $\otimes _{A}:{_{A}\mathcal{M}_{A}}\times {_{A}%
\mathcal{M}_{A}}\rightarrow {_{A}\mathcal{M}}_{A}$ and constraints that make
the category $({_{A}\mathcal{M}}_{A},\otimes _{A},A)$ abelian monoidal, see
\cite[1.11]{AMS}. The tensor product over $A$ in $\mathcal{M}$ of a right $A$%
-module $V$ and a left $A$-module $W$ is defined to be the coequalizer:
\begin{equation*}
\xymatrix{ (V\otimes A)\otimes W \ar@<.5ex>[rr] \ar@<-.5ex>[rr]&& V\otimes W
\ar[rr]^{_{A}\chi _{V,W}} && V\otimes _{A}W \ar[r] & 0 }
\end{equation*}%
Note that, since $\otimes $ preserves coequalizers, then $V\otimes _{A}W$ is
also an $A$-bimodule, whenever $V$ and $W$ are $A$-bimodules.\medskip
\newline
Dually, given a coalgebra $(C,\Delta ,\varepsilon )$ in a coabelian monoidal
category $\mathcal{M}$, there exist a suitable functor $\square _{C}:{^{C}%
\mathcal{M}^{C}}\times {^{C}\mathcal{M}^{C}}\rightarrow {^{C}\mathcal{M}^{C}}
$ and constraints that make the category $({^{C}\mathcal{M}^{C}},\square
_{C},C)$ coabelian monoidal. The cotensor product over $C$ in $\mathcal{M}$
of a right $C$-bicomodule $V$ and a left $C$-comodule $W$ is defined to be
the equalizer:
\begin{equation*}
\xymatrix{ 0 \ar[r] & V\cot_{C}W \ar[rr]^{_C\varsigma_{V,W}} && V\otimes W
\ar@<.5ex>[rr] \ar@<-.5ex>[rr]&&V\ot(C\ot W) }
\end{equation*}%
Note that, since $\otimes $ preserves equalizers, then $V\square _{C}W$ is
also a $C$-bicomodule, whenever $V$ and $W$ are $C$-bicomodules.

\begin{claim}
\textbf{Graded Objects.} \label{claim 4.2}Let $\left( X_{n}\right) _{n\in
\mathbb{N}
}$ be a sequence of objects a cocomplete abelian monoidal category $\mathcal{%
M}$ and let
\begin{equation*}
X=\bigoplus_{n\in \mathbb{N}}X_{n}
\end{equation*}%
be their coproduct in $\mathcal{M}$. In this case we also say that $X$ is a
\emph{graded object of }$\mathcal{M}$ and that the sequence $\left(
X_{n}\right) _{n\in
\mathbb{N}
}$ defines a \emph{graduation} on $X.$ A morphism
\begin{equation*}
f:X=\bigoplus_{n\in \mathbb{N}}X_{n}\rightarrow Y=\bigoplus_{n\in \mathbb{N}%
}Y_{n}
\end{equation*}%
is called a \emph{graded homomorphism} whenever there exists a family of
morphisms $\left( f_{n}:X_{n}\rightarrow Y_{n}\right) _{n\in
\mathbb{N}
}$ such that $f=\oplus _{n\in \mathbb{N}}f_{n}$ i.e. such that%
\begin{equation*}
f\circ i_{X_{n}}^{X}=i_{Y_{n}}^{Y}\circ f_{n},\text{ for every }n\in \mathbb{%
N}\text{.}
\end{equation*}%
We fix the following notations:%
\begin{equation*}
X\left( n\right) =\left\{
\begin{tabular}{ll}
$0,$ & for $n=0$ \\
$\underset{0\leq i\leq n-1}{\oplus }X_{i},$ & for $n\geq 1$%
\end{tabular}%
\right. \qquad \text{and}\qquad X\left[ n\right] =\underset{i\geq n}{\oplus }%
X_{i},\text{ for }n\in
\mathbb{N}
.
\end{equation*}%
Let $\sigma _{i}^{i+1}:X\left( i\right) \rightarrow X\left( i+1\right) $ be
the canonical inclusion and for any $j>i,$ define:
\begin{equation*}
\sigma _{i}^{j}=\sigma _{j-1}^{j}\sigma _{j-2}^{j-1}\cdots \sigma
_{i+1}^{i+2}\sigma _{i}^{i+1}:X\left( i\right) \rightarrow X\left( j\right) .
\end{equation*}%
Throughout let
\begin{eqnarray*}
\pi _{n}^{m} &:&X\left( n\right) \rightarrow X\left( m\right) \text{ }(m\leq
n),\text{\qquad }\pi _{n}:X\rightarrow X\left( n\right) , \\
p_{n}^{m} &:&X\left( n\right) \rightarrow X_{m}\text{ }(m<n),\text{\qquad }%
p_{n}:X\rightarrow X_{n},
\end{eqnarray*}%
be the canonical projections and let
\begin{eqnarray*}
\sigma _{m}^{n} &:&X\left( m\right) \rightarrow X\left( n\right) \text{ }%
(m\leq n),\text{\qquad }\sigma _{n}:X\left( n\right) \rightarrow X, \\
i_{m}^{n} &:&X_{m}\rightarrow X\left( n\right) \text{ }(m<n),\text{\qquad }%
i_{n}:X_{n}\rightarrow X,
\end{eqnarray*}%
be the canonical injections for any $m,n\in \mathbb{N}$. For technical
reasons we set $\pi _{n}^{m}=0,$ $\sigma _{m}^{n}=0$ for any $n<m$ and $%
p_{n}^{m}=0$, $i_{m}^{n}=0$ for any $n\leq m.$ Then, we have the following
relations:%
\begin{equation*}
p_{n}\sigma _{k}=p_{k}^{n},\qquad p_{n}i_{k}=\delta _{n,k}\mathrm{Id}%
_{X_{k}},\qquad \pi _{n}i_{k}=i_{k}^{n}.
\end{equation*}%
Moreover, we have:
\begin{equation*}
\begin{tabular}{lll}
$\pi _{n}^{m}\sigma _{k}^{n}=\sigma _{k}^{m},\text{ if }k\leq m\leq n,$ & $%
\qquad $and$\qquad $ & $\pi _{n}^{m}\sigma _{k}^{n}=\pi _{k}^{m},\text{ if }%
m\leq k\leq n,$ \\
$p_{n}^{m}\pi _{k}^{n}=p_{k}^{m},\text{ if }m<n\leq k,$ & $\qquad $and$%
\qquad $ & $\sigma _{n}^{m}i_{k}^{n}=i_{k}^{m},\text{ if }k<n\leq m,$ \\
$p_{n}^{m}\sigma _{k}^{n}=p_{k}^{m},\text{ if }m<k\leq n,$ & $\qquad $and$%
\qquad $ & $\pi _{n}^{m}i_{k}^{n}=i_{k}^{m},\text{ if }k<m\leq n,$ \\
$p_{n}^{m}\pi _{n}=p_{m},\text{ if }m<n,$ & $\qquad $and$\qquad $ & $\sigma
_{n}i_{m}^{n}=i_{m},\text{ if }m<n\text{,}$ \\
$\pi _{n}\sigma _{k}=\sigma _{k}^{n},\text{ if }k\leq n,$ & $\qquad $and$%
\qquad $ & $\pi _{n}\sigma _{k}=\pi _{k}^{n},\text{ if }n\leq k,$ \\
$p_{n}^{m}i_{m}^{n}=\mathrm{Id}_{X_{m}},\text{ if }m<n.$ &  &
\end{tabular}%
\end{equation*}%
In the other cases, these compositions are zero.\newline
Similarly let
\begin{eqnarray*}
\tau _{n}^{m} &:&X\left[ n\right] \rightarrow X\left[ m\right] \text{ }%
(n\leq m),\text{\qquad }\tau _{n}:X\rightarrow X\left[ n\right] , \\
\zeta _{n}^{m} &:&X\left[ n\right] \rightarrow X_{m}\text{ }(n\leq m),
\end{eqnarray*}%
be the canonical projections and let
\begin{eqnarray*}
\nu _{m}^{n} &:&X\left[ m\right] \rightarrow X\left[ n\right] \text{ }(n\leq
m),\text{\qquad }\nu _{n}:X\left[ n\right] \rightarrow X, \\
\lambda _{m}^{n} &:&X_{m}\rightarrow X\left[ n\right] \text{ }(n\leq m),
\end{eqnarray*}%
be the canonical injection for any $m,n\in \mathbb{N}$. For technical
reasons we set $\tau _{n}^{m}=0,$ $\nu _{m}^{n}=0$ for any $m<n$ and $\zeta
_{n}^{m}=0$, $\lambda _{m}^{n}=0$ for any $m<n.$ Then, we have the following
relations:%
\begin{equation*}
p_{n}\nu _{k}=\zeta _{k}^{n},\qquad \tau _{n}i_{k}=\lambda _{k}^{n}.
\end{equation*}%
Moreover, we have:
\begin{equation*}
\begin{tabular}{lll}
$\tau _{n}^{m}\nu _{k}^{n}=\nu _{k}^{m},\text{ if }n\leq m\leq k,$ & $\qquad
$and$\qquad $ & $\tau _{n}^{m}\nu _{k}^{n}=\tau _{k}^{m},\text{ if }n\leq
k\leq m,$ \\
$\zeta _{n}^{m}\tau _{k}^{n}=\zeta _{k}^{m},\text{ if }k\leq n\leq m,$ & $%
\qquad $and$\qquad $ & $\nu _{n}^{m}\lambda _{k}^{n}=\lambda _{k}^{m},\text{
if }m\leq n\leq k,$ \\
$\zeta _{n}^{m}\nu _{k}^{n}=\zeta _{k}^{m},\text{ if }n\leq k\leq m,$ & $%
\qquad $and$\qquad $ & $\tau _{n}^{m}\lambda _{k}^{n}=\lambda _{k}^{m},\text{
if }n\leq m\leq k,$ \\
$\zeta _{n}^{m}\tau _{n}=p_{m},\text{ if }n\leq m,$ & $\qquad $and$\qquad $
& $\nu _{n}\lambda _{m}^{n}=i_{m},\text{ if }n\leq m\text{,}$ \\
$\tau _{n}\nu _{k}=\nu _{k}^{n},\text{ if }n\leq k,$ & $\qquad $and$\qquad $
& $\tau _{n}\nu _{k}=\tau _{k}^{n},\text{ if }k\leq n,$ \\
$\zeta _{n}^{m}\lambda _{m}^{n}=\mathrm{Id}_{X_{m}},\text{ if }n\leq m.$ &
&
\end{tabular}%
\end{equation*}%
In the other cases, these compositions are zero. \newline
Given graded objects $X,Y$ in $\mathcal{M}$ we set%
\begin{equation*}
\left( X\otimes Y\right) _{n}=\oplus _{a+b=n}\left( X_{a}\otimes
Y_{b}\right) .
\end{equation*}%
Then this defines a graduation on $X\otimes Y$ whenever{\Huge \ }the tensor
product commutes with direct sums. We denote by
\begin{equation*}
X_{a}\otimes Y_{b}\overset{\gamma _{a,b}^{X,Y}}{\rightarrow }\left( X\otimes
Y\right) _{a+b}\text{\qquad and\qquad }\left( X\otimes Y\right) _{a+b}%
\overset{\omega _{a,b}^{X,Y}}{\rightarrow }X_{a}\otimes Y_{b}
\end{equation*}%
the canonical injection and projection respectively. We have%
\begin{eqnarray}
\sum\limits_{a+b=n}\left( i_{a}^{X}\otimes i_{b}^{Y}\right) \omega
_{a,b}^{X,Y} &=&\nabla \left[ \left( i_{a}^{X}\otimes i_{b}^{Y}\right)
_{a+b=n}\right]  \label{form nablasum1} \\
\sum\limits_{a+b=n}\gamma _{a,b}^{X,Y}\left( p_{a}^{X}\otimes
p_{b}^{Y}\right) &=&\Delta \left[ \left( p_{a}^{X}\otimes p_{b}^{Y}\right)
_{a+b=n}\right]  \label{form deltasum1}
\end{eqnarray}
where $\nabla \left[ \left( i_{a}^{X}\otimes i_{b}^{Y}\right) _{a+b=n}\right]
$ denotes the codiagonal morphism associated to the family $\left(
i_{a}^{X}\otimes i_{b}^{Y}\right) _{a+b=n}$ and $\Delta \left[ \left(
p_{a}^{X}\otimes p_{b}^{Y}\right) _{a+b=n}\right] $ denotes the diagonal
morphism associated to the family $\left( p_{a}^{X}\otimes p_{b}^{Y}\right)
_{a+b=n}.$
\end{claim}

\begin{lemma}
\label{lem: exactSeq}Let $X=\bigoplus_{n\in \mathbb{N}}X_{n}$ be a graded
object in a cocomplete abelian category $\mathcal{M}$. Then, for every $n\in
\mathbb{N}
$%
\begin{equation*}
0\rightarrow X\left[ n\right] \overset{\nu _{n}}{\longrightarrow }X\overset{%
\pi _{n}}{\longrightarrow }X\left( n\right) \rightarrow 0\qquad \text{and}%
\qquad 0\rightarrow X\left( n\right) \overset{\sigma _{n}}{\longrightarrow }X%
\overset{\tau _{n}}{\longrightarrow }X\left[ n\right] \rightarrow 0
\end{equation*}%
are splitting exact sequences.
\end{lemma}

\section{Graded and Strongly $\mathbb{N}$-graded Coalgebras\label{sec:
Graded Coalg}}

In this section $\left( \mathcal{M},\otimes ,\mathbf{1}\right) $ will denote
a cocomplete coabelian monoidal category such that the tensor product
commutes with direct sums.

\begin{claim}
\label{def: grCoalg} Recall that a \emph{graded coalgebra} in $\mathcal{M}$
is a coalgebra $\left( C,\Delta ,\varepsilon \right) $ where
\begin{equation*}
C=\oplus _{n\in \mathbb{N}}C_{n}
\end{equation*}%
is a graded object of $\mathcal{M}$ such that $\Delta :C\rightarrow C\otimes
C$ is a graded homomorphism i.e. there exists a family $\left( \Delta
_{n}\right) _{_{n\in \mathbb{N}}}$ of morphisms
\begin{equation*}
\Delta _{n}^{C}=\Delta _{n}:C_{n}\rightarrow \left( C\otimes C\right)
_{n}=\oplus _{a+b=n}\left( C_{a}\otimes C_{b}\right) \text{ such that }%
\Delta =\oplus _{n\in \mathbb{N}}\Delta _{n}.
\end{equation*}%
We set
\begin{equation*}
\Delta _{a,b}^{C}=\Delta _{a,b}:=\left( C_{a+b}\overset{\Delta _{a+b}}{%
\rightarrow }\left( C\otimes C\right) _{a+b}\overset{\omega _{a,b}^{C,C}}{%
\rightarrow }C_{a}\otimes C_{b}\right) .
\end{equation*}%
A homomorphism $f:\left( C,\Delta _{C},\varepsilon _{C}\right) \rightarrow
\left( D,\Delta _{D},\varepsilon _{D}\right) $ of coalgebras is a graded
coalgebra homomorphism if it is a graded homomorphism too.
\end{claim}

\begin{lemma}
\label{lem: complete}Let $I$ be a set and let $(X_{i})_{i\in I}$ be a family
of objects in a cocomplete and complete abelian category $\mathcal{C}$
satisfying $AB5$. Let $Y$ be an object in $\mathcal{C}$ and let $%
f:Y\rightarrow \oplus _{i\in I}X_{i}$ be a morphism such that
\begin{equation*}
p_{k}f=0\text{ for any }k\in I,
\end{equation*}%
where $p_{k}:\oplus _{i\in I}X_{i}\rightarrow X_{k}$ denotes the canonical
projection. Then $f=0.$
\end{lemma}

\begin{proof}
See \cite[page 54 and 61]{Po}.
\end{proof}

\begin{proposition}
\label{pro: grCoalg} Let $\left( C,\Delta ,\varepsilon \right) $ be a
coalgebra in $\mathcal{M}$ which is a graded object with graduation defined
by $\left( C_{k}\right) _{k\in \mathbb{N}}.$ Fix $n\in
\mathbb{N}
.$ Consider the following assertions.

\begin{enumerate}
\item[(i)] There exists a morphism
\begin{equation*}
\Delta _{n}:C_{n}\rightarrow \left( C\otimes C\right) _{n}=\oplus
_{a+b=n}\left( C_{a}\otimes C_{b}\right) \text{ such that }\Delta
i_{n}=\nabla \left[ \left( i_{a}\otimes i_{b}\right) _{a+b=n}\right] \Delta
_{n}.
\end{equation*}

\item[(ii)] There exists a family $\left( \Delta _{a,b}\right) _{_{a+b=n}}$
of morphisms $\Delta _{a,b}:C_{n}\rightarrow C_{a}\otimes C_{b}$ such that%
\begin{equation}
\Delta i_{n}=\sum\limits_{a+b=n}\left( i_{a}\otimes i_{b}\right) \Delta
_{a,b}.  \label{form: grCoalg}
\end{equation}

\item[(iii)] $(p_{a}\otimes p_{b})\Delta i_{n}=0,$ for every $a,b\in
\mathbb{N}
,a+b\neq n.$
\end{enumerate}

Then $\left( i\right) \Leftrightarrow \left( ii\right) \Rightarrow \left(
iii\right) .$ Furthermore, if $\mathcal{M}$ is also complete and satisfies $%
AB5,$ then $\left( ii\right) \Leftrightarrow \left( iii\right) $.
\end{proposition}

\begin{proof}
$\left( i\right) \Rightarrow \left( ii\right) $ Set%
\begin{equation*}
\Delta _{a,b}=\left( C_{n}\overset{\Delta _{n}}{\rightarrow }\left( C\otimes
C\right) _{n}\overset{\omega _{a,b}^{C,C}}{\rightarrow }C_{a}\otimes
C_{b}\right)
\end{equation*}%
for every $a,b\in
\mathbb{N}
$ such that $a+b=n$. We have
\begin{equation*}
\sum\limits_{a+b=n}\left( i_{a}\otimes i_{b}\right) \Delta
_{a,b}=\sum\limits_{a+b=n}\left( i_{a}\otimes i_{b}\right) \omega
_{a,b}^{C,C}\Delta _{n}\overset{(\ref{form nablasum1})}{=}\nabla \left[
\left( i_{a}\otimes i_{b}\right) _{a+b=n}\right] \Delta _{n}=\Delta i_{n}.
\end{equation*}%
$\left( ii\right) \Rightarrow \left( iii\right) $ We have%
\begin{equation}
\left( p_{a}\otimes p_{b}\right) \sum\limits_{u+v=n}\left( i_{u}\otimes
i_{v}\right) \Delta _{u,v}=\delta _{a+b,n}\Delta _{a,b}=\Delta
_{a,b}p_{a+b}i_{n},  \label{formula: grcoalg1}
\end{equation}%
so that, in view of (\ref{form: grCoalg}), we obtain%
\begin{equation*}
\left( p_{a}\otimes p_{b}\right) \Delta i_{n}=\delta _{a+b,n}\Delta _{a,b}%
\text{, for every }a,b\in
\mathbb{N}
\text{.}
\end{equation*}%
$\left( ii\right) \Rightarrow \left( i\right) $ Set%
\begin{equation*}
\Delta _{n}=\sum\limits_{a+b=n}\left( C_{n}\overset{\Delta _{a,b}}{%
\rightarrow }C_{a}\otimes C_{b}\overset{\gamma _{a,b}^{C,C}}{\rightarrow }%
\left( C\otimes C\right) _{n}\right) .
\end{equation*}%
Then%
\begin{equation*}
\nabla \left[ \left( i_{u}\otimes i_{v}\right) _{u+v=n}\right] \Delta
_{n}=\sum\limits_{a+b=n}\nabla \left[ \left( i_{u}\otimes i_{v}\right)
_{u+v=n}\right] \gamma _{a,b}^{C,C}\Delta _{a,b}=\sum\limits_{a+b=n}\left(
i_{a}\otimes i_{b}\right) \Delta _{a,b}=\Delta i_{n}.
\end{equation*}%
Assume $\mathcal{M}$ is also complete and satisfies $AB5$.

$\left( iii\right) \Rightarrow \left( ii\right) $ Set $\Delta
_{a,b}=(p_{a}\otimes p_{b})\Delta i_{a+b}$ for every $a,b\in
\mathbb{N}
$. For every $u,v\in
\mathbb{N}
,$ we have%
\begin{equation*}
(p_{u}\otimes p_{v})\sum\limits_{a+b=n}\left( i_{a}\otimes i_{b}\right)
\Delta _{a,b}\overset{\text{(\ref{formula: grcoalg1})}}{=}\delta
_{u+v,n}\Delta _{u,v}=(p_{u}\otimes p_{v})\Delta i_{n}.
\end{equation*}%
since the codomain of $\Delta i_{n}$ is $C\otimes C=\oplus _{u,v\in
\mathbb{N}
}C_{u}\otimes C_{v}$ and $\mathcal{M}$ satisfies $AB5$, by Lemma \ref{lem:
complete}, we conclude.
\end{proof}

\begin{proposition}
\label{coro: grCoalg}1) If $\left( C=\oplus _{n\in
\mathbb{N}
}C_{n},\Delta ,\varepsilon \right) $ is a graded coalgebra in $\mathcal{M}$,
then
\begin{equation}
(p_{a}\otimes p_{b})\Delta =\Delta _{a,b}p_{n},\text{ for every }a,b\in
\mathbb{N}
.  \label{form: coro grCoalg1}
\end{equation}%
2) Assume that $\mathcal{M}$ is also complete and satisfies $AB5.$ Let $%
\left( C,\Delta ,\varepsilon \right) $ be a coalgebra in $\mathcal{M}$ which
is a graded object with graduation defined by $\left( C_{n}\right) _{n\in
\mathbb{N}
}.$ If there exists a family $\left( \Delta _{a,b}\right) _{_{a,b\in
\mathbb{N}
}}$ of morphisms $\Delta _{a,b}:C_{a+b}\rightarrow C_{a}\otimes C_{b}$ such
that (\ref{form: coro grCoalg1}) holds, then $\left( C=\oplus _{n\in
\mathbb{N}
}C_{n},\Delta ,\varepsilon \right) $ is a graded coalgebra in $\mathcal{M}$.
\end{proposition}

\begin{proof}
$1)$ For every $a,b,t\in
\mathbb{N}
$ we have%
\begin{equation*}
(p_{a}\otimes p_{b})\Delta i_{t}\overset{\text{(\ref{form: grCoalg})}}{=}%
(p_{a}\otimes p_{b})\sum\limits_{u+v=n}\left( i_{u}\otimes i_{v}\right)
\Delta _{u,v}\overset{\text{(\ref{formula: grcoalg1})}}{=}\Delta
_{a,b}p_{a+b}i_{n}.
\end{equation*}
Since this equality holds for an arbitrary $n\in
\mathbb{N}
$, we get (\ref{form: coro grCoalg1}).

$2)$ By (\ref{form: coro grCoalg1}) we have $(p_{a}\otimes p_{b})\Delta
i_{n}=0,$ for every $a,b\in
\mathbb{N}
,a+b\neq n.$ Thus, by Proposition \ref{pro: grCoalg}, $C$ is a graded
coalgebra.
\end{proof}

\begin{proposition}
\label{lem: graded Deltaij}1) Let $C=\oplus _{n\in \mathbb{N}}C_{n}$ be a
graded object of $\mathcal{M}$ such that there exists a family $\left(
\Delta _{a,b}\right) _{_{a,b\in
\mathbb{N}
}}$
\begin{equation*}
\Delta _{a,b}:C_{a+b}\rightarrow C_{a}\otimes C_{b},
\end{equation*}%
of morphisms and a morphism $\varepsilon _{0}:C_{0}\rightarrow \mathbf{1}$
which satisfy
\begin{gather}
\left( \Delta _{a,b}\otimes C_{c}\right) \Delta _{a+b,c}=\left( C_{a}\otimes
\Delta _{b,c}\right) \Delta _{a,b+c}\text{,}  \label{form: locDelta} \\
\left( C_{d}\otimes \varepsilon _{0}\right) \Delta
_{d,0}=r_{C_{d}}^{-1},\qquad \left( \varepsilon _{0}\otimes C_{d}\right)
\Delta _{0,d}=l_{C_{d}}^{-1}\text{,}  \label{form: locEps}
\end{gather}%
for every $a,b,c\in
\mathbb{N}
$. Then there exists a unique morphism $\Delta :C\rightarrow C\otimes C$
such that (\ref{form: grCoalg}) holds.\newline
Moreover $\left( C=\oplus _{n\in \mathbb{N}}C_{n},\Delta ,\varepsilon
=\varepsilon _{0}p_{0}\right) $ is a graded coalgebra.

2) If $C$ is a graded coalgebra, then $\varepsilon =\varepsilon i_{0}p_{0}$
so that $\varepsilon $ is a graded homomorphism, and we have that (\ref%
{form: locDelta}) and (\ref{form: locEps}) hold for every $a,b,c\in
\mathbb{N}
$, where $\varepsilon _{0}=\varepsilon i_{0}$.\newline
Moreover $\left( C_{0},\Delta _{0}=\Delta _{0,0},\varepsilon
_{0}=\varepsilon i_{0}\right) $ is a coalgebra in $\mathcal{M}$, $i_{0}$ is
a coalgebra homomorphism and, for every $n\in
\mathbb{N}
$, $\left( C_{n},\Delta _{0,n},\Delta _{n,0}\right) $ is a $C_{0}$%
-bicomodule such that $p_{n}:C\rightarrow C_{n}$ is a morphism of $C_{0}$%
-bicomodules ($C$ is a $C_{0}$-bicomodule through $p_{0}$).
\end{proposition}

\begin{proof}
1) Since the tensor product commutes with direct sums, we have that $%
C\otimes C=\oplus _{a,b\in
\mathbb{N}
}\left( C_{a}\otimes C_{b}\right) .$\newline
Thus, by the universal property of coproduct there exists a unique morphism $%
\Delta :C\rightarrow C\otimes C$ such that (\ref{form: grCoalg}) holds. By
using (\ref{form: grCoalg}) and (\ref{form: locDelta}) we get $\left( \Delta
\otimes C\right) \Delta i_{n}=\left( C\otimes \Delta \right) \Delta i_{n},$
for every $n\in
\mathbb{N}
$ and hence $\Delta $ is coassociative. Similarly, by using (\ref{form:
grCoalg}) and (\ref{form: locEps}), we get $\left( C\otimes \varepsilon
\right) \Delta =r_{C}^{-1}$ and $\left( \varepsilon \otimes C\right) \Delta
=l_{C}^{-1}.$ In view of (\ref{form: grCoalg}), by Proposition \ref{pro:
grCoalg}, $C$ is a graded coalgebra.

2) By using (\ref{form: coro grCoalg1}), we obtain%
\begin{eqnarray*}
\left( p_{a}\otimes p_{b}\otimes p_{c}\right) \left( \Delta \otimes C\right)
\Delta i_{a+b+c} &=&\left( \Delta _{a,b}\otimes C_{c}\right) \Delta _{a+b,c}
\\
\left( p_{a}\otimes p_{b}\otimes p_{c}\right) \left( C\otimes \Delta \right)
\Delta i_{a+b+c} &=&\left( C_{a}\otimes \Delta _{b,c}\right) \Delta _{a,b+c}
\end{eqnarray*}%
so that (\ref{form: locDelta}) holds true.

From $\left( C\otimes \varepsilon _{C}\right) \Delta =r_{C}^{-1}$, we deduce%
\begin{equation*}
\left( p_{b}\otimes \mathbf{1}\right) \left( C\otimes \varepsilon
_{C}\right) \Delta i_{a}=\left( p_{b}\otimes \mathbf{1}\right)
r_{C}^{-1}i_{a}.
\end{equation*}%
We compute the first term:%
\begin{equation*}
\left( p_{b}\otimes \mathbf{1}\right) \left( C\otimes \varepsilon \right)
\Delta i_{a}\overset{(\ref{form: grCoalg})}{=}\left( p_{b}\otimes \mathbf{1}%
\right) \left( C\otimes \varepsilon \right) \sum\limits_{u+v=a}\left(
i_{u}\otimes i_{v}\right) \Delta _{u,v}=\left( C_{b}\otimes \varepsilon
i_{a-b}\right) \Delta _{b,a-b}
\end{equation*}%
and the second one:%
\begin{equation*}
\left( p_{b}\otimes \mathbf{1}\right)
r_{C}^{-1}i_{a}=r_{C_{b}}^{-1}p_{b}i_{a}=\delta _{b,a}r_{C_{b}}^{-1}\text{.}
\end{equation*}%
Thus we obtain
\begin{equation*}
\left( C_{b}\otimes \varepsilon i_{a-b}\right) \Delta _{b,a-b}=\delta
_{b,a}r_{C_{b}}^{-1}.
\end{equation*}%
Similarly one gets $\left( \varepsilon _{C}i_{a-b}\otimes C_{b}\right)
\Delta _{a-b,b}=\delta _{a,b}l_{C_{b}}^{-1}.$ In particular one has (\ref%
{form: locEps}).

Let us prove that $\varepsilon $ is a graded homomorphism. Let us check that%
\begin{equation*}
\varepsilon =\varepsilon i_{0}p_{0}.
\end{equation*}%
Set $\varepsilon ^{\prime }:=\varepsilon i_{0}p_{0}.$ Then%
\begin{equation*}
r_{C}\left( A\otimes \varepsilon ^{\prime }\right) \Delta i_{n}\overset{%
\text{(\ref{form: grCoalg})}}{=}r_{C}\left( A\otimes \varepsilon
i_{0}p_{0}\right) \sum\limits_{a+b=n}\left( i_{a}\otimes i_{b}\right) \Delta
_{a,b}=r_{C}\left( i_{n}\otimes \varepsilon i_{0}\right) \Delta _{n,0}%
\overset{\text{(\ref{form: locEps})}}{=}r_{C}\left( i_{n}\otimes \mathbf{1}%
\right) r_{C_{n}}^{-1}=i_{n}\text{.}
\end{equation*}%
Since this holds for every $n\in
\mathbb{N}
$, we deduce that $r_{C}\left( A\otimes \varepsilon ^{\prime }\right) \Delta
=\mathrm{Id}_{C}$ so that $\varepsilon ^{\prime }$ is a right counit for $C.$
Since $\varepsilon $ is a counit for $C$, one gets $\varepsilon =\varepsilon
^{\prime }.$

Hence $\varepsilon $ is a graded homomorphism.

By applying (\ref{form: locDelta}) and (\ref{form: locEps}) to the cases
\begin{equation*}
\left( a,b,c\right) =\left( 0,0,0\right) ,\text{ }d=0,
\end{equation*}%
we get that $\left( C_{0},\Delta _{0}=\Delta _{0,0},\varepsilon
_{0}=\varepsilon i_{0}\right) $ is a coalgebra in $\mathcal{M}$. \newline
By applying (\ref{form: locDelta}) and (\ref{form: locEps}) to the cases
\begin{equation*}
\left( a,b,c\right) =\left( n,0,0\right) ,\left( 0,n,0\right) ,\left(
0,0,n\right) \text{, }d=n,
\end{equation*}%
we get that $\left( C_{n},\Delta _{0,n},\Delta _{n,0}\right) $ is a $C_{0}$%
-bicomodule for every $n\in
\mathbb{N}
$. By applying (\ref{form: coro grCoalg1}) to the cases $\left( a,b\right)
=\left( n,0\right) $ and $\left( a,b\right) =\left( 0,n\right) ,$ we get the
last assertion.
\end{proof}

\begin{lemma}
\label{lem: grcoalHomo}Let $f:\left( C,\Delta _{C},\varepsilon _{C}\right)
\rightarrow \left( D,\Delta _{D},\varepsilon _{D}\right) $ be a homomorphism
of graded coalgebras in $\mathcal{M}$. Then
\begin{equation}
\Delta _{a,b}^{D}\circ f_{a+b}=\left( f_{a}\otimes f_{b}\right) \circ \Delta
_{a,b}^{C}\text{, for every }n,a,b\in
\mathbb{N}
.  \label{form: grSubcoalg1}
\end{equation}
\end{lemma}

\begin{proof}
By definition of graded homomorphism, $f\circ i_{t}^{C}=i_{t}^{D}\circ f$
and hence $p_{n}^{D}\circ f\circ i_{t}^{C}=f\circ p_{n}^{C}\circ i_{t}^{C}$
for every $n,t\in
\mathbb{N}
$. Therefore $p_{n}^{D}\circ f=f_{n}\circ p_{n}^{C}$ for every $n\in
\mathbb{N}
$. Using this relation and (\ref{form: coro grCoalg1}), for every $a,b\in
\mathbb{N}
$, the equality%
\begin{equation*}
(p_{a}^{D}\otimes p_{b}^{D})\circ \Delta _{D}\circ f\circ
i_{a+b}^{C}=(p_{a}^{D}\otimes p_{b}^{D})\circ \left( f\otimes f\right) \circ
\Delta _{C}\circ i_{a+b}^{C}
\end{equation*}%
rewrites as $\Delta _{a,b}^{D}\circ f_{a+b}=\left( f_{a}\otimes f_{b}\right)
\circ \Delta _{a,b}^{C}.$
\end{proof}

\begin{claim}
\label{cl: wedge} Let $C$ be a coalgebra in $\mathcal{M}$. As in the case of
vector spaces, we can introduce the wedge product of two subobjects $X,Y$ of
$C$ in $\mathcal{M}:$%
\begin{equation*}
(X\wedge _{C}Y,i_{X\wedge Y}^{C}):=Ker[(p_{X}\otimes p_{Y})\circ \triangle
_{C}],
\end{equation*}%
where $p_{X}:C\rightarrow C/X$ and $p_{Y}:C\rightarrow C/Y$ are the
canonical quotient maps. In particular we have the following exact sequence:
\begin{equation*}
\xymatrix{ 0 \ar[r] & X\wedge_C Y \ar[rr]^{i_{X\w Y}^C} && C \ar[rr]^(.4){(p
_{X}\otimes p _{Y}) \circ \triangle _{C}} && C/X \ot C/Y.}
\end{equation*}%
Let $\delta :D\rightarrow C$ be a monomorphism which is a homomorphism of
coalgebras in $\mathcal{M}$. Denote by $(L,p)$ the cokernel of $\delta $ in $%
\mathcal{M}$. Regard $D$ as a $C$-bicomodule via $\delta $ and observe that $%
L$ is a $C$-bicomodule and $p$ is a morphism of bicomodules. Let
\begin{equation*}
(D^{\wedge _{C}^{n}},\delta _{n}):=\ker (p^{\otimes {n}}\Delta _{C}^{n-1})
\end{equation*}%
for any $n\in \mathbb{N}\setminus \{0\}.$ Note that $(D^{\wedge
_{C}^{1}},\delta _{1})=(D,\delta )$ and $(D^{\wedge _{C}^{2}},\delta
_{2})=D\wedge _{C}D.$ \newline
In order to simplify the notations we set $(D^{\wedge _{C}^{0}},\delta
_{0})=(0,0).$\newline
Now, since $\mathcal{M}$ has left exact tensor functors and since $%
p^{\otimes {n}}\Delta _{C}^{n-1}$ is a morphism of $C$-bicomodules (as a
composition of morphisms of $C$-bicomodules), we get that $D^{\wedge
_{C}^{n}}$ is a coalgebra and $\delta _{n}:D^{\wedge _{C}^{n}}\rightarrow C$
is a coalgebra homomorphism for any $n>0$ and hence for any $n\in \mathbb{N}$%
.
\end{claim}

\begin{proposition}
\cite[Proposition 2.17]{AMS:Cotensor} \label{pro: D^2}Let $\delta
:D\rightarrow C$ be a monomorphism which is a coalgebra homomorphism in $%
\mathcal{M}$. Then we have
\begin{equation}
D^{\wedge _{C}^{m}}\wedge _{C}D^{\wedge _{C}^{n}}=D^{\wedge _{C}^{m+n}}.
\label{formula 2 pro: D^2}
\end{equation}
\end{proposition}

\begin{definition}
\label{def: strongly grCoalg}Let $(C=\oplus _{n\in
\mathbb{N}
}C_{n},\Delta ,\varepsilon )$ be a graded coalgebra in $\mathcal{M}$. In
analogy with the group graded case (see \cite{NT}), we say that $C$ is a
\emph{strongly }$%
\mathbb{N}
$\emph{-graded coalgebra} whenever

$\Delta _{i,j}^{C}:C_{i+j}\rightarrow C_{i}\otimes C_{j}$ is a monomorphism
for every $i,j\in \mathbb{N},$

where $\Delta _{i,j}^{C}$ is the morphism defined in Definition \ref{def:
grCoalg}.
\end{definition}

\begin{theorem}
\label{teo: DeltaInj} Let $(C=\oplus _{n\in
\mathbb{N}
}C_{n},\Delta ,\varepsilon )$ be a graded coalgebra in a cocomplete and
complete coabelian monoidal category $\mathcal{M}$ satisfying $AB5.$ \newline
Then:

1) Let $n\geq 2$ be such that
\begin{equation}
\left\{
\begin{tabular}{c}
there exist $0\leq a,b\leq n-1$ such that $a+b=n$ \\
and $\Delta _{u,v}:C_{u+v}\rightarrow C_{u}\otimes C_{v}$ is a monomorphism
for every $a\leq u$ and $b\leq v$%
\end{tabular}%
\right\}  \label{condition DeltaInj}
\end{equation}%
\begin{equation*}
0\rightarrow C\left( n\right) \overset{\sigma _{n}}{\longrightarrow }C%
\overset{\left( \tau _{a}\otimes \tau _{b}\right) \Delta }{\longrightarrow }C%
\left[ a\right] \otimes C\left[ b\right]
\end{equation*}%
is an exact sequence and%
\begin{equation*}
C\left( n\right) =C\left( a\right) \wedge _{C}C\left( b\right) \text{.}
\end{equation*}%
2) Assume that $C$ is a strongly $%
\mathbb{N}
$-graded coalgebra. Then
\begin{equation*}
(C\left( n\right) ,\sigma _{n})=C_{0}^{\wedge _{C}^{n}},\text{ for every }%
n\in \mathbb{N}.
\end{equation*}
\end{theorem}

\begin{proof}
1) Let $n\geq 2$ be such that condition (\ref{condition DeltaInj}) holds
true. Let us prove that the following sequence%
\begin{equation*}
0\rightarrow C\left( n\right) \overset{\sigma _{n}}{\longrightarrow }C%
\overset{\left( \tau _{a}\otimes \tau _{b}\right) \Delta }{\longrightarrow }C%
\left[ a\right] \otimes C\left[ b\right]
\end{equation*}%
is exact. For every $0\leq t\leq n-1,$ we compute%
\begin{equation*}
\left( \tau _{a}\otimes \tau _{b}\right) \Delta \sigma _{n}i_{t}^{n}=\left(
\tau _{a}\otimes \tau _{b}\right) \Delta i_{t}\overset{\text{(\ref{form:
grCoalg})}}{=}\sum\limits_{u+v=t}\left( \tau _{a}i_{u}\otimes \tau
_{b}i_{v}\right) \Delta _{u,v}=\sum\limits_{u+v=t}\left( \lambda
_{u}^{a}\otimes \lambda _{v}^{b}\right) \Delta _{u,v}=0.
\end{equation*}%
Let $g:Y\rightarrow C$ be a morphism such that $\left( \tau _{a}\otimes \tau
_{b}\right) \Delta g=0.$ Now, for every $a\leq u$ and $b\leq v,$ we have%
\begin{equation*}
0=\left( \zeta _{a}^{u}\tau _{a}\otimes \zeta _{b}^{v}\tau _{b}\right)
\Delta g=\left( p_{u}\otimes p_{v}\right) \Delta g\overset{\text{(\ref{form:
coro grCoalg1})}}{=}\Delta _{u,v}p_{u+v}g.
\end{equation*}%
By hypothesis $\Delta _{u,v}$ is a monomorphism so that $p_{u+v}g=0.$ We
conclude that%
\begin{equation*}
p_{c}g=0,\text{ for every }c\geq n.
\end{equation*}%
Set $\overline{g}=\pi _{n}g$ and let us prove that $g=\sigma _{n}\overline{g}%
.$ By Lemma \ref{lem: complete} this is the case if and only if%
\begin{equation*}
p_{c}g=p_{c}\sigma _{n}\overline{g},\text{ for every }c\in \mathbb{N}.
\end{equation*}%
We calculate%
\begin{equation*}
p_{c}\sigma _{n}\overline{g}=p_{n}^{c}\overline{g}=p_{n}^{c}\pi
_{n}g=\left\{
\begin{tabular}{ll}
$p_{c}g$ & for every $c<n$ \\
$0=p_{c}g$ & otherwise.%
\end{tabular}%
\right.
\end{equation*}%
Thus the required sequence is exact.

By Lemma \ref{lem: exactSeq}, the following sequence%
\begin{equation*}
0\rightarrow C\left( n\right) \overset{\sigma _{n}}{\longrightarrow }C%
\overset{\tau _{n}}{\longrightarrow }C\left[ n\right] \rightarrow 0
\end{equation*}%
is split exact. Then we have%
\begin{equation*}
\left( C\left( n\right) ,\sigma _{n}\right) =\ker \left[ \left( \tau
_{a}\otimes \tau _{b}\right) \Delta \right] =\ker \left( \tau _{a}\right)
\wedge _{C}\ker \left( \tau _{b}\right) =C\left( a\right) \wedge _{C}C\left(
b\right) .
\end{equation*}%
2) Assume that $C$ is a strongly $%
\mathbb{N}
$-graded coalgebra. Let us prove that $C\left( n\right) =C_{0}^{\wedge
_{C}^{n}},$ for every $n\in \mathbb{N}.$ \newline
The case $n=0$ is trivial. Let us prove the equality above for every $n\geq
1 $ by induction on $n.$\newline
If $n=1,$ by definition, we have $C\left( 1\right) =C_{0}=C_{0}^{\wedge
_{C}^{1}}.$\newline
Let $n\geq 2$ and assume that $C\left( t\right) =C_{0}^{\wedge _{C}^{t}}$
for every $0\leq t\leq n-1.$ By hypothesis we get that condition (\ref%
{condition DeltaInj}) holds true for $n$ and hence by 2) we have $C\left(
n\right) =C\left( a\right) \wedge _{C}C\left( b\right) $ so that, by
Proposition \ref{pro: D^2}, we have%
\begin{equation*}
C_{0}^{\wedge _{C}^{n}}=C_{0}^{\wedge _{C}^{a}}\wedge _{C}C_{0}^{\wedge
_{C}^{b}}=C\left( a\right) \wedge _{C}C\left( b\right) =C\left( n\right) .
\end{equation*}
\end{proof}

\begin{proposition}
\label{pro: grSubcoalg}Assume that $\mathcal{M}$ is also complete and
satisfies $AB5.$ Let $(C=\oplus _{n\in
\mathbb{N}
}C_{n},\Delta _{C},\varepsilon _{C})$ and $(D=\oplus _{n\in
\mathbb{N}
}D_{n},\Delta _{D},\varepsilon _{D})$ be graded coalgebras in $\mathcal{M}$.
Let $f:C\rightarrow D$ be a graded coalgebra homomorphism. Assume that $D$
is a strongly $%
\mathbb{N}
$-graded coalgebra. Then the following assertions are equivalent

$\left( a\right) $ $f_{0}$ and $f_{1}$ are monomorphisms and $C$ is a
strongly $\mathbb{N}$-graded coalgebra.

$\left( a^{\prime }\right) $ $f_{0}$ and $f_{1}$ are monomorphisms and $%
\Delta _{a,1}^{C}:C_{a+1}\rightarrow C_{a}\otimes C_{1}$ is a monomorphism
for every $a\in
\mathbb{N}
$.

$\left( b\right) $ $f_{n}$ is a monomorphism for every $n\in
\mathbb{N}
$.

$\left( c\right) $ $f$ is a monomorphism.

$\left( d\right) $ $f_{0}$ and $f_{1}$ are monomorphisms and $(C\left(
n\right) ,\sigma _{n})=C_{0}^{\wedge _{C}^{n}},$ for every $n\in \mathbb{N}.$

$\left( e\right) $ $f_{0}$ and $f_{1}$ are monomorphisms and $(C\left(
2\right) ,\sigma _{2})=C_{0}^{\wedge _{C}^{2}}$.
\end{proposition}

\begin{proof}
Recall that, by Lemma \ref{lem: grcoalHomo}, we have (\ref{form: grSubcoalg1}%
).

$\left( a\right) \Rightarrow \left( a^{\prime }\right) $ It is trivial.

$\left( a^{\prime }\right) \Rightarrow \left( b\right) $ We proceed by
induction on $a\geq 1$. For $a=1$ there is nothing to prove. Let $a\geq 2$
and assume that $f_{a-1}$ is a monomorphism. Since $\Delta
_{a-1,1}^{C}:C_{a+1}\rightarrow C_{a}\otimes C_{1}$ is a monomorphism, in
view of (\ref{form: grSubcoalg1}), we get that $f_{a}=f_{\left( a-1\right)
+1}$ is a monomorphism.

$\left( b\right) \Rightarrow \left( c\right) $ Since $f=\oplus _{n\in
\mathbb{N}
}f_{n}$ and $\mathcal{M}$ satisfies $AB5$, it is trivial (see e.g. \cite[%
page 53]{Po}).

$\left( c\right) \Rightarrow \left( b\right) $ Since $f$ is a monomorphism
and $f\circ i_{n}^{C}=i_{n}^{D}\circ f_{n}$ for every $n\in
\mathbb{N}
$, it is clear that each $f_{n}$ is a monomorphism too.

$\left( b\right) \Rightarrow \left( a\right) $ Let $a,b\geq 0$. In view of (%
\ref{form: grSubcoalg1}), $\Delta _{a,b}^{C}$ is a monomorphism as $\Delta
_{a,b}^{D}$ is a monomorphism.

$\left( a\right) \Rightarrow \left( d\right) $ It follows by Theorem \ref%
{teo: DeltaInj}.\newline
$\left( d\right) \Rightarrow \left( e\right) $ It is trivial.\newline
$\left( e\right) \Rightarrow \left( a\right) $ By hypothesis,%
\begin{equation*}
\ker \left[ \left( \tau _{1}^{C}\otimes \tau _{1}^{C}\right) \Delta _{C}%
\right] =\left( C_{0}^{\wedge _{C}^{2}},\delta _{2}^{C}\right) =(C\left(
2\right) ,\sigma _{2}^{C})=\ker \left( \tau _{2}^{C}\right)
\end{equation*}%
so that
\begin{equation*}
\mathrm{Im}\left[ \left( \tau _{1}^{C}\otimes \tau _{1}^{C}\right) \Delta
_{C}\right] =\mathrm{\mathrm{\mathrm{\mathrm{\mathrm{\mathrm{\mathrm{Coker}}}%
}}}}\left[ \ker \left( \tau _{2}^{C}\right) \right] =\mathrm{Im}\left( \tau
_{2}^{C}\right) =\left( C\left[ 2\right] ,\tau _{2}^{C}\right) .
\end{equation*}%
Hence there is a monomorphism $\beta :C\left[ 2\right] \rightarrow C\left[ 1%
\right] \otimes C\left[ 1\right] $ such that%
\begin{equation*}
\beta \circ \tau _{2}^{C}=\left( \tau _{1}^{C}\otimes \tau _{1}^{C}\right)
\circ \Delta _{C}.
\end{equation*}%
For every $n\in
\mathbb{N}
$, we have%
\begin{equation*}
\beta \circ \lambda _{n}^{2}=\beta \circ \tau _{2}^{C}\circ i_{n}^{C}=\left(
\tau _{1}^{C}\otimes \tau _{1}^{C}\right) \circ \Delta _{C}\circ i_{n}^{C}%
\overset{\text{(\ref{form: grCoalg})}}{=}\left( \tau _{1}^{C}\otimes \tau
_{1}^{C}\right) \sum\limits_{a+b=n}\left( i_{a}^{C}\otimes i_{b}^{C}\right)
\circ \Delta _{a,b}^{C}=\sum\limits_{a+b=n}\left( \lambda _{a}^{1}\otimes
\lambda _{b}^{1}\right) \circ \Delta _{a,b}^{C}
\end{equation*}%
so that%
\begin{equation}
\beta \circ \lambda _{n}^{2}=\sum\limits_{a+b=n}\left( \lambda
_{a}^{1}\otimes \lambda _{b}^{1}\right) \circ \Delta _{a,b}^{C}.
\label{form: grSubcoalg0}
\end{equation}

Let us prove that $\Delta _{a,b}^{C}:C_{a+b}\rightarrow C_{a}\otimes C_{b}$
is a monomorphism for every $0\leq a,b$.

For $a=0$ or $b=0$ there is nothing to prove as $\left( C_{n},\Delta
_{0,n}^{C},\Delta _{n,0}^{C}\right) $ is a $C_{0}$-bicomodule in view of
Proposition \ref{lem: graded Deltaij}. Let us check that $\Delta
_{a,b}^{C}:C_{a+b}\rightarrow C_{a}\otimes C_{b}$ is a monomorphism for
every $a,b$ by induction on $n=a+b$.

For $n=1$ there is nothing to prove.

For $n=2$, since $\lambda _{0}^{1}=0,$ we have%
\begin{equation*}
\beta \circ \lambda _{2}^{2}\overset{\text{(\ref{form: grSubcoalg0})}}{=}%
\sum\limits_{a+b=2}\left( \lambda _{a}^{1}\otimes \lambda _{b}^{1}\right)
\circ \Delta _{a,b}^{C}=\left( \lambda _{1}^{1}\otimes \lambda
_{1}^{1}\right) \circ \Delta _{1,1}^{C}.
\end{equation*}%
So that, since $\beta $ and $\lambda _{i}^{j}$ are monomorphisms for every $%
j\leq i,$ we have%
\begin{equation*}
\mathrm{\ker }\left( \Delta _{1,1}^{C}\right) =\mathrm{\ker }\left[ \left(
\lambda _{1}^{1}\otimes \lambda _{1}^{1}\right) \circ \Delta _{1,1}^{C}%
\right] =\mathrm{\ker }\left( \beta \circ \lambda _{2}^{2}\right) =0.
\end{equation*}

Let $n\geq 3$ and assume that $\Delta _{u,v}^{C}:C_{u+v}\rightarrow
C_{u}\otimes C_{v}$ is a monomorphism for every $u,v\in
\mathbb{N}
$ with $u+v<n$. Fix $a,b\in
\mathbb{N}
$ such that $n=a+b$ and let $h$ be a morphism in $\mathcal{M}$ such that $%
\Delta _{a,b}^{C}\circ h=0.$ Let us prove that $h=0.$

By (\ref{form: locDelta}), for every $i,j>0$ with $n=i+j,$ we have%
\begin{eqnarray*}
\left( \Delta _{i,a-i}^{C}\otimes C_{b}\right) \Delta _{a,b}^{C} &=&\left(
C_{i}\otimes \Delta _{j-b,b}^{C}\right) \Delta _{i,j}^{C},\text{ if }a\geq i%
\text{ and hence }j\geq b \\
\left( \Delta _{a,i-a}^{C}\otimes C_{j}\right) \Delta _{i,j}^{C} &=&\left(
C_{a}\otimes \Delta _{b-j,j}^{C}\right) \Delta _{a,b}^{C},\text{ if }i\geq a%
\text{ and hence }b\geq j.
\end{eqnarray*}%
In the first case (the other one is analogous) we have%
\begin{equation*}
0=\left( \Delta _{i,a-i}^{C}\otimes C_{b}\right) \Delta _{a,b}^{C}h=\left(
C_{i}\otimes \Delta _{j-b,b}^{C}\right) \Delta _{i,j}^{C}h.
\end{equation*}%
Since $\left( j-b\right) +b=j<n,$ we have that $\Delta _{j-b,b}^{C}$ is a
monomorphism and hence $C_{i}\otimes \Delta _{j-b,b}^{C}$ is a monomorphism
too. Thus $\Delta _{i,j}^{C}h=0$ for every $i,j>0$ with $n=i+j.$

Then, since $\lambda _{0}^{1}=0,$ we have%
\begin{equation*}
\beta \circ \lambda _{n}^{2}\circ h\overset{\text{(\ref{form: grSubcoalg0})}}%
{=}\sum\limits_{i+j=n}\left( \lambda _{i}^{1}\otimes \lambda _{j}^{1}\right)
\circ \Delta _{i,j}^{C}\circ h=\sum\limits_{\substack{ i+j=n  \\ i\neq
0,j\neq 0}}\left( \lambda _{i}^{1}\otimes \lambda _{j}^{1}\right) \circ
\Delta _{i,j}^{C}\circ h=0
\end{equation*}

Since $\beta $ and $\lambda _{n}^{2}$ are monomorphisms, we get $h=0.$
\end{proof}

\begin{claim}
\textbf{The Cotensor Coalgebra.} \label{claim: cotensor coalg} For any
object $X\in \mathcal{M},$ we can define the cotensor coalgebra $T^{c}\left(
X\right) $ as follows. As an object in $\mathcal{M}$ it is defined by
setting $\ $
\begin{equation*}
T^{c}=T^{c}\left( X\right) :=\oplus _{p\in \mathbb{N}}X^{\otimes p}.
\end{equation*}

Let $p_{a}:T(X)\rightarrow X^{\otimes a}$ be the canonical projection and
let $\varepsilon _{0}:=\mathrm{Id}_{\mathbf{1}}.$

Let%
\begin{equation*}
\Delta _{p,q}:X^{\otimes p+q}\rightarrow X^{\otimes p}\otimes X^{\otimes q},%
\text{ for every }p,q\in \mathbb{N}
\end{equation*}%
be the canonical isomorphism, which is unique by Coherence Theorem. \newline
Still by Coherence Theorem and the definitions of $\Delta _{d,0}$ and $%
\Delta _{0,d}$, one has%
\begin{gather*}
\lbrack \Delta _{p,q}\otimes X^{\otimes r}]\Delta _{p+q,r}=[X^{\otimes
p}\otimes \Delta _{q,r}]\Delta _{p,q+r}. \\
\left( X^{\otimes d}\otimes \varepsilon _{0}\right) \Delta
_{d,0}=r_{X^{\otimes d}}^{-1},\qquad \left( u_{0}\otimes X^{\otimes
d}\right) \Delta _{0,d}=l_{X^{\otimes d}}^{-1}
\end{gather*}

so that (\ref{form: locDelta}) and (\ref{form: locEps}) hold.

By Proposition \ref{lem: graded Deltaij}, there there exists a unique
morphism $\Delta _{T^{c}}:T^{c}\rightarrow T^{c}\otimes T^{c}$ such that (%
\ref{form: grCoalg}) holds i.e.%
\begin{equation}
\Delta \circ i_{n}=\sum\limits_{a+b=n}\left( i_{a}\otimes i_{b}\right) \circ
\Delta _{a,b}.  \label{eq: deltaCotens}
\end{equation}

Moreover $\left( T^{c},\Delta _{T^{c}},\varepsilon _{T^{c}}=\varepsilon
_{0}p_{0}=p_{0}\right) $ is a graded coalgebra with graduation defined by $%
\left( X^{\otimes p}\right) _{p\in \mathbb{N}}$. Let $(C,\Delta
_{C},\varepsilon _{C})$ be a coalgebra in $\mathcal{M}$.\newline
As explained in Section \ref{sec: Preliminaries and Notations}, we have that
$({^{C}\mathcal{M}^{C}},\square _{C},C)$ is a coabelian monoidal category.
Furthermore one can see that ${^{C}\mathcal{M}^{C}}$ has also arbitrary
direct sums which commute with $\square _{C}$. \newline
Therefore we can consider, in the monoidal category $({^{C}\mathcal{M}^{C}}%
,\square _{C},C),$ the cotensor coalgebra of an arbitrary $C$-bimodule $M.$
We will denote it by
\begin{equation*}
\left( T^{c}=T_{C}^{c}(M),\overline{\Delta }_{T^{c}},\overline{\varepsilon
_{T^{c}}}\right) .
\end{equation*}%
Note that%
\begin{equation*}
T^{c}=\oplus _{p\in \mathbb{N}}M^{\square _{C}p},\qquad \overline{\Delta }%
_{T^{c}}:T^{c}\rightarrow T^{c}\square _{C}T^{c}\qquad \text{and}\qquad
\overline{\varepsilon _{T^{c}}}:T^{c}\rightarrow C.
\end{equation*}%
Set%
\begin{equation*}
\Delta _{T}={_{C}\varsigma _{V,W}}\circ \overline{\Delta }%
_{T}:T^{c}\rightarrow T^{c}\otimes T^{c}\qquad \text{and}\qquad \varepsilon
_{T^{c}}=\varepsilon _{C}\circ \overline{\varepsilon _{T^{c}}}%
:T^{c}\rightarrow \mathbf{1},
\end{equation*}%
where ${_{C}\varsigma _{V,W}}:T^{c}\square _{C}T^{c}\longrightarrow
T^{c}\otimes T^{c}$ is the canonical morphism introduced in Section \ref%
{sec: Preliminaries and Notations}.
\end{claim}

\begin{proposition}
\label{pro: cotGrad}Let $\left( C,\Delta _{C},\varepsilon _{C}\right) $ be a
coalgebra in $\mathcal{M}$ and let $(M,\rho _{M}^{r},\rho _{M}^{l})$ be a $C$%
-bicomodule. Let $T^{c}=T_{C}^{c}(M)$ be the cotensor coalgebra. Then $%
T^{c}=\oplus _{n\in \mathbb{N}}T_{n}^{c}$ is a graded coalgebra where $%
T_{n}^{c}=M^{\square _{C}n}.$ Moreover we have%
\begin{equation*}
\Delta _{m,n}=\left\{
\begin{tabular}{ll}
$(M^{\square m-1}\square \varsigma _{M}\square M^{\square n-1}),$ & $\text{%
for any }m,n\geq 1;$ \\
$(M^{\square m-1}\square \rho _{M}^{r}),$ & $\text{for any }m\geq 1,n=0;$ \\
$(\rho _{M}^{l}\square M^{\square n-1}),$ & $\text{for any }m=0,n\geq 1;$ \\
$\Delta _{C},$ & $\text{for any }m=0,n=0.$%
\end{tabular}%
\right.
\end{equation*}%
In particular $T^{c}$ is a strongly $%
\mathbb{N}
$-graded coalgebra.
\end{proposition}

\begin{proof}
It follows by (\ref{claim: cotensor coalg}) and by Proposition \ref{pro:
grCoalg} as (\ref{form: grCoalg}) is clearly satisfied.
\end{proof}

\begin{remark}
Note that the coalgebra structure $T_{C}^{c}(M)$ in Proposition \ref{pro:
cotGrad} is the same introduced in \cite[Theorem 2.9]{AMS:Cotensor}.
\end{remark}

We have the following result.

\begin{notation}
\label{notation tilde}Let $\delta :D\rightarrow C$ be a homomorphism of
coalgebras in $\mathcal{M}$. By Proposition \cite[Proposition 1.10]%
{AMS:Cotensor}, $((D^{\wedge _{C}^{i}})_{i\in \mathbb{N}},(\xi
_{i}^{j})_{i,j\in \mathbb{N}})$ is a direct system in $\mathcal{M}$ whose
direct limit carries a natural coalgebra structure that makes it the direct
limit of $((D^{\wedge _{C}^{i}})_{i\in \mathbb{N}},(\xi _{i}^{j})_{i,j\in
\mathbb{N}})$ as a direct system of coalgebras.\newline
From now on we set: $(\widetilde{D}_{C},(\xi _{i})_{i\in \mathbb{N}})=%
\underrightarrow{\lim }(D^{\wedge _{C}^{i}})_{i\in \mathbb{N}}$, where $\xi
_{i}:D^{\wedge _{C}^{i}}\rightarrow \widetilde{D}_{C}$ denotes the
structural morphism of the direct limit. We simply write $\widetilde{D}$ if
there is no danger of confusion. We note that, since $\widetilde{D}$ is a
direct limit of coalgebras, the canonical (coalgebra) homomorphisms $(\delta
_{i}:D^{\wedge _{C}^{i}}\rightarrow C)_{i\in \mathbb{N}}$, which are
compatible, factorize to a unique coalgebra homomorphism $\widetilde{\delta }%
:\widetilde{D}\rightarrow C$ such that $\widetilde{\delta }\xi _{i}=\delta
_{i}$ for any $i\in \mathbb{N}$.
\end{notation}

\begin{theorem}
\cite[Theorems 2.13 and 2.15]{AMS:Cotensor}\label{coro: univ property of
cotensor coalgebra}\label{teo: pre univ property of cotensor coalgebra} Let $%
(C,\Delta ,\varepsilon )$ be a coalgebra in a cocomplete and complete
coabelian monoidal category $\mathcal{M}$ satisfying AB5. Let $(M,\rho
_{M}^{r},\rho _{M}^{l})$ be a $C$-bicomodule. Let $\delta :D\rightarrow E$
be a monomorphism which is a homomorphism of coalgebras. Let $f_{C}:%
\widetilde{D}\rightarrow C$ be a coalgebra homomorphism and let $f_{M}:%
\widetilde{D}\rightarrow M$ be a morphism of $C$-bicomodules such that $%
f_{M}\xi _{1}=0$, where $\widetilde{D}$ is a bicomodule via $f_{C}.$ Then
there is a unique coalgebra homomorphism $f:\widetilde{D}\rightarrow
T_{C}^{c}(M)$ such that $p_{0}f=f_{C}$ and $p_{1}f=f_{M}$, where $%
p_{n}:T_{C}^{c}(M)\rightarrow M^{\square n}$ denotes the canonical
projection.
\begin{equation*}
\diagUnivCot%
\end{equation*}%
Moreover%
\begin{equation}
p_{n}\circ f=f_{M}^{\square _{C}n}\circ \overline{\Delta }_{\widetilde{D}%
}^{n-1}\text{ for any }n\in \mathbb{N}.  \label{relation f}
\end{equation}
\end{theorem}

\begin{proposition}
\label{pro: ciocco}\cite[Proposition 2.11]{AMS:Cotensor} Let $(C,\Delta
,\varepsilon )$ be a coalgebra in $\mathcal{M}$ and let $(M,\rho
_{M}^{r},\rho _{M}^{l})$ be a $C$-bicomodule. Let $T^{c}:=T_{C}^{c}(M).$ Let
$E$ be a coalgebra and let $\alpha :E\rightarrow T^{c}$ and $\beta
:E\rightarrow T^{c}$ be coalgebra homomorphisms. If $p_{1}\alpha =p_{1}\beta
,$ then $p_{n}\alpha =p_{n}\beta $ for any $n\geq 1.$
\end{proposition}

\begin{corollary}
\label{coro: uniqueness} Assume that $\mathcal{M}$ is also complete and
satisfies $AB5$. Let $(C,\Delta ,\varepsilon )$ be a coalgebra in $\mathcal{M%
}$ and let $(M,\rho _{M}^{r},\rho _{M}^{l})$ be a $C$-bicomodule. Let $%
T^{c}:=T_{C}^{c}(M).$ Let $E$ be a coalgebra and let $\alpha :E\rightarrow
T^{c}$ and $\beta :E\rightarrow T^{c}$ be coalgebra homomorphisms.\newline
Then $\alpha =\beta $ whenever $p_{n}\alpha =p_{n}\beta ,$ for $n=0,1$.
\end{corollary}

\begin{proof}
It follows by Proposition \ref{pro: ciocco} and Lemma \ref{lem: complete}.
\end{proof}

\begin{proposition}
\cite[Proposition 3.3]{AMS:Cotensor} \label{pro: lim for graded coalg}Let $%
C=\oplus _{t\in \mathbb{N}}C_{t}$ be a graded coalgebra in $\mathcal{M}$.
Then%
\begin{equation}
\tau _{1}^{\otimes n+1}\Delta _{C}^{n}i_{b}=0,\text{ for every }0\leq b\leq
n.  \label{formula lem: c^n in wedge}
\end{equation}%
Moreover%
\begin{equation*}
C=\underrightarrow{\lim }(C_{0}^{\wedge _{C}^{t}})_{t\in \mathbb{N}}.
\end{equation*}
\end{proposition}

\begin{theorem}
\label{coro: GRuniv property of cotensor coalgebra}Let $(\mathcal{M},\otimes
,\mathbf{1})$ be a cocomplete and complete coabelian monoidal category
satisfying $AB5$ and such that the tensor product commutes with direct sums.%
\newline
Let $(C,\Delta ,\varepsilon )$ be a coalgebra in $\mathcal{M}$ and let $%
(M,\rho _{M}^{r},\rho _{M}^{l})$ be a $C$-bicomodule. \newline
Let $B$ be a graded coalgebra, let $g_{C}:B_{0}\rightarrow C$ be a coalgebra
homomorphism and let $g_{M}:B_{1}\rightarrow M$ be a morphism in ${^{C}%
\mathcal{M}^{C}},$ where $B_{1}\in {^{C}\mathcal{M}^{C}}$ via $g_{C}.$%
\newline
Then there is a unique coalgebra homomorphism $f:B\rightarrow T_{C}^{c}(M)$
such that
\begin{equation*}
p_{0}^{T^{c}}\circ f=g_{C}\circ p_{0}^{B}\text{\qquad and\qquad }%
p_{1}^{T^{c}}\circ f=g_{M}\circ p_{1}^{B}.
\end{equation*}%
Moreover $f$ is a morphism of graded coalgebras where%
\begin{equation*}
f_{t}=\left( g_{M}\circ p_{1}^{B}\right) ^{\square _{C}t}\circ \overline{%
\Delta }_{B}^{t-1}\circ i_{t}^{B}\text{ for every }t\in
\mathbb{N}%
\end{equation*}%
where $\overline{\Delta }_{B}^{t-1}:B\rightarrow B^{\square _{C}t}$ is the $%
t $-th iterated comultiplication of $B$.
\end{theorem}

\begin{proof}
In view of Proposition \ref{lem: graded Deltaij}, $\left( B_{0},\Delta
_{0}=\Delta _{0,0},\varepsilon _{0}=\varepsilon i_{0}\right) $ is a
coalgebra in $\mathcal{M}$, $i_{0}^{B}$ is a coalgebra homomorphism and $%
\left( B_{1},\Delta _{0,1},\Delta _{1,0}\right) $ is a $B_{0}$-bicomodule.
Moreover, $p_{0}^{B}:B\rightarrow B_{0}$ is a coalgebra homomorphism (see %
\ref{form: coro grCoalg1}) and $p_{1}^{B}:B\rightarrow B_{1}$ is a morphism
of $B_{0}$-bicomodules ($B$ is a $B_{0}$-bicomodule through $p_{0}$).\newline
By Theorem \ref{coro: univ property of cotensor coalgebra} applied in the
case when $\delta $ is the morphism $i_{0}^{B}:B_{0}\rightarrow B$ and by
Proposition \ref{pro: lim for graded coalg}, there is a unique coalgebra
homomorphism $f:B\rightarrow T_{C}^{c}(M)$ such that
\begin{equation*}
p_{0}^{T^{c}}\circ f=g_{C}\circ p_{0}^{B}\text{\qquad and\qquad }%
p_{1}^{T^{c}}\circ f=g_{M}\circ p_{1}^{B}.
\end{equation*}
By (\ref{relation f}), we have%
\begin{equation*}
p_{t}^{T^{c}}\circ f=\left( g_{M}\circ p_{1}^{B}\right) ^{\square
_{C}t}\circ \overline{\Delta }_{B}^{t-1}\text{ for any }t\in \mathbb{N}.
\end{equation*}

Let us prove there is a family $\left( f_{t}\right) _{t\in
\mathbb{N}
}$ of morphisms such that%
\begin{equation}
p_{t}^{T^{c}}\circ f=f_{t}\circ p_{t}^{B}.  \label{form: univgrCo}
\end{equation}%
By (\ref{form: coro grCoalg1}) we have%
\begin{equation*}
(p_{a}^{B}\otimes p_{b}^{B})\Delta _{B}=\Delta _{a,b}^{B}p_{a+b}^{B}
\end{equation*}%
By Proposition \ref{lem: graded Deltaij}, for every $n\in
\mathbb{N}
,$ $p_{n}^{B}:B\rightarrow B_{n}$ is a morphism of $B_{0}$-bicomodules and
hence of $C$-bicomodules via $g_{C}.$ Thus%
\begin{equation*}
(p_{a}^{B}\square _{C}p_{b}^{B})\overline{\Delta }_{B}=\overline{\Delta }%
_{a,b}^{B}p_{a+b}^{B}
\end{equation*}%
Hence there is a morphism $\alpha _{t}:\left( B_{1}\right) ^{\otimes
_{A}t}\rightarrow B_{t}$ such that%
\begin{equation*}
\left( p_{1}^{B}\right) ^{\square _{C}t}\circ \overline{\Delta }%
_{B}^{t-1}=\alpha _{t}\circ p_{t}^{B}.
\end{equation*}%
Then $\alpha _{t}=\left( p_{1}^{B}\right) ^{\square _{C}t}\circ \overline{%
\Delta }_{B}^{t-1}\circ i_{t}^{B}.$ We have%
\begin{equation*}
p_{t}^{T^{c}}\circ f=\left( g_{M}\circ p_{1}^{B}\right) ^{\square
_{C}t}\circ \overline{\Delta }_{B}^{t-1}=\left( g_{M}\right) ^{\square
_{C}t}\circ \left( p_{1}^{B}\right) ^{\square _{C}t}\circ \overline{\Delta }%
_{B}^{t-1}=\left( g_{M}\right) ^{\square _{C}t}\circ \alpha _{t}\circ
p_{t}^{B}.
\end{equation*}%
so that (\ref{form: univgrCo}) holds true for $f_{t}:=\left( g_{M}\right)
^{\square _{C}t}\circ \alpha _{t}=\left( g_{M}\right) ^{\square _{C}t}\circ
\left( p_{1}^{B}\right) ^{\square _{C}t}\circ \overline{\Delta }%
_{B}^{t-1}\circ i_{t}^{B}.$

Now, $f$ is a morphism of graded coalgebras if and only if $ f\circ
i_{t}^{B}=i_{t}^{T^{c}}\circ f_{t}\text{, for every }t\in
\mathbb{N}
\text{.} $ In view of Lemma \ref{lem: complete}, this is equivalent to prove
that $p_{s}^{T^{c}}\circ f\circ i_{t}^{B}=p_{s}^{T^{c}}\circ
i_{t}^{T^{c}}\circ f_{t}\text{, for every }s,t\in
\mathbb{N}%
$ which is trivially true in view of (\ref{form: univgrCo}).
\end{proof}

\begin{remark}
In the case $C=\mathbf{1}$, the universal property of the cotensor coalgebra
stated in Theorem \ref{coro: GRuniv property of cotensor coalgebra}, differs
from \cite[Remark 2.3]{Schauenburg1} where it is proved that there exists a
unique \textbf{graded }coalgebra homomorphism $f:B\rightarrow T_{C}^{c}(M)$
such that $f_{1}=g_{M}$. In fact we need extra assumptions as $AB5$
condition in order to establish the uniqueness of $f$ as a coalgebra
homomorphism. This uniqueness allows us to prove that $f$ is graded whenever
$B$ is graded.
\end{remark}

\begin{theorem}
\label{teo: psi mono}Let $(C=\oplus _{n\in
\mathbb{N}
}C_{n},\Delta _{C},\varepsilon _{C})$ be a graded coalgebra in a cocomplete
and complete coabelian monoidal category $\mathcal{M}$ satisfying $AB5.$
Assume that that the tensor product commutes with direct sums.\newline
Let $T^{c}:=T_{C_{0}}^{c}(C_{1})$ be the cotensor coalgebra. Then there is a
unique coalgebra homomorphism
\begin{equation*}
\psi :C\rightarrow T_{C_{0}}^{c}(C_{1}),
\end{equation*}%
such that $p_{0}^{T^{c}}\circ \psi =p_{0}^{C}$ and $p_{1}^{T^{c}}\circ \psi
=p_{1}^{C}$. \newline
Moreover $\psi $ is a graded coalgebra homomorphism with%
\begin{equation*}
\psi _{m}=\left( p_{1}^{C}\right) ^{\square m}\circ \overline{\Delta }%
_{C}^{m-1}\circ i_{m}^{C}\text{ for every }m\in
\mathbb{N}%
\end{equation*}%
and the following assertions are equivalent.

$\left( a\right) $ $C$ is a strongly $%
\mathbb{N}
$-graded coalgebra.

$\left( a^{\prime }\right) $ $\Delta _{a,1}^{C}:C_{a+1}\rightarrow
C_{a}\otimes C_{1}$ is a monomorphism for every $a\in
\mathbb{N}
$.

$\left( b\right) $ $\psi _{n}$ is a monomorphism for every $n\in
\mathbb{N}
$.

$\left( c\right) $ $\psi $ is a monomorphism.

$\left( d\right) $ $(C\left( n\right) ,\sigma _{n}^{C})=C_{0}^{\wedge
_{C}^{n}},$ for every $n\in \mathbb{N}.$

$\left( e\right) $ $(C\left( 2\right) ,\sigma _{2}^{C})=C_{0}^{\wedge
_{C}^{2}}$.
\end{theorem}

\begin{proof}
In view of Proposition \ref{lem: graded Deltaij}, $C_{1}$ is a bicomodule
over the coalgebra $C_{0}$ so that we can consider $T_{C_{0}}^{c}(C_{1})$.
By Theorem \ref{coro: GRuniv property of cotensor coalgebra}, there is a
unique coalgebra homomorphism $\psi :C\rightarrow $ $T_{C_{0}}^{c}(C_{1})\ $%
such that $p_{0}^{T^{c}}\circ \psi =p_{0}^{C}$ and $p_{1}^{T^{c}}\circ \psi
=p_{1}^{C}$. Moreover $\psi $ is a morphism of graded coalgebras where%
\begin{equation*}
\psi _{t}=\left( p_{1}^{C}\right) ^{\square t}\circ \overline{\Delta }%
_{C}^{t-1}\circ i_{t}^{C}\text{ for every }t\in
\mathbb{N}%
\end{equation*}%
where $\overline{\Delta }_{C}^{t-1}:C\rightarrow C^{\square _{C_{0}}t}$ is
the $t$-th iterated comultiplication of $C$.\newline
Since, by Proposition \ref{pro: cotGrad}, $T^{c}$ is a strongly $%
\mathbb{N}
$-graded coalgebra, we can apply Proposition \ref{pro: grSubcoalg}.

We conclude by observing that%
\begin{equation*}
\psi _{0}=\left( p_{1}^{C}\right) ^{\square 0}\circ \overline{\Delta }%
_{C}^{0-1}\circ i_{0}^{C}=p_{0}^{C}\circ i_{0}^{C}=\mathrm{Id}_{C_{0}}\quad
\text{and}\quad \psi _{1}=\left( p_{1}^{C}\right) ^{\square 1}\circ
\overline{\Delta }_{C}^{1-1}\circ i_{1}^{C}=p_{1}^{C}\circ i_{1}^{C}=\mathrm{%
Id}_{C_{1}}\text{.}
\end{equation*}
\end{proof}

\section{Graded and Strongly $\mathbb{N}$-graded Algebras}

In this section $\left( \mathcal{M},\otimes ,\mathbf{1}\right) $ will denote
a cocomplete abelian monoidal category such that the tensor product commutes
with direct sums. Analogous results to those included in Section \ref{sec:
Graded Coalg} can be obtained for graded algebras by a careful
"dualization". Still we would like to point out that here we don't need to
assume that $\mathcal{M}$ is also complete and satisfies $AB5$ essentially
because the analogue of Lemma \ref{lem: complete} with the injections
doesn't require these extra conditions.

\begin{definition}
\label{def: grAlg} Recall that a \emph{graded algebra} in $\mathcal{M}$ is
an algebra $\left( A,m,u\right) $ where
\begin{equation*}
A=\oplus _{n\in \mathbb{N}}A_{n}
\end{equation*}%
is a graded object of $\mathcal{M}$ such that $m:A\otimes A\rightarrow A$ is
a graded homomorphism i.e. there exists a family $\left( m_{n}\right)
_{_{n\in \mathbb{N}}}$ of morphisms
\begin{equation*}
m_{n}^{A}=m_{n}:\oplus _{a+b=n}\left( A_{a}\otimes A_{b}\right) =\left(
A\otimes A\right) _{n}\rightarrow A_{n}\text{ such that }m=\oplus _{n\in
\mathbb{N}}m_{n}.
\end{equation*}%
We set
\begin{equation*}
m_{a,b}^{A}:=\left( A_{a}\otimes A_{b}\overset{\gamma _{a,b}^{A,A}}{%
\rightarrow }\left( A\otimes A\right) _{a+b}\overset{m_{a+b}}{\rightarrow }%
A_{a+b}\right) .
\end{equation*}%
A homomorphism $f:\left( A,m_{A},u_{A}\right) \rightarrow \left(
B,m_{B},u_{B}\right) $ of algebras is a graded algebra homomorphism if it is
a graded homomorphism too.
\end{definition}

\begin{proposition}
\label{pro: grAlg}Let $\left( A,m,u\right) $ be an algebra in $\mathcal{M}$
which is a graded object with graduation defined by $\left( A_{k}\right)
_{k\in \mathbb{N}}.$ Fix $n\in
\mathbb{N}
.$ The following assertions are equivalent.

\begin{enumerate}
\item[(i)] There exists a morphism
\begin{equation*}
m_{n}:\oplus _{a+b=n}\left( A_{a}\otimes A_{b}\right) =\left( A\otimes
A\right) _{n}\rightarrow A_{n}\text{ such that }p_{n}m=m_{n}\Delta \left[
\left( p_{a}\otimes p_{b}\right) _{a+b=n}\right] \text{.}
\end{equation*}

\item[(ii)] There exists a family $\left( m_{a,b}\right) _{_{a+b=n}}$ of
morphisms $m_{a,b}:A_{a}\otimes A_{b}\rightarrow A_{n}$ such that%
\begin{equation}
p_{n}m=\sum\limits_{a+b=n}m_{a,b}\left( p_{a}\otimes p_{b}\right) .
\label{form: grAlg}
\end{equation}

\item[(iii)] $p_{n}m(i_{a}\otimes i_{b})=0,$ for every $a,b\in
\mathbb{N}
,a+b\neq n.$
\end{enumerate}
\end{proposition}

\begin{proof}
It is analogous to that of Proposition \ref{pro: grCoalg}.
\end{proof}

\begin{proposition}
\label{coro: grAlg}$1)$ If $\left( A=\oplus _{n\in \mathbb{N}%
}A_{n},m,u\right) $ is a graded algebra in $\mathcal{M}$, then
\begin{equation}
m(i_{a}\otimes i_{b})=i_{a+b}m_{a,b},\text{ for every }a,b\in
\mathbb{N}
.  \label{form: coro grAlg1}
\end{equation}%
$2)$ Let $\left( A,m,u\right) $ be an algebra in $\mathcal{M}$ which is a
graded object with graduation defined by $\left( A_{k}\right) _{k\in \mathbb{%
N}}.$ If there exists a family $\left( m_{a,b}\right) _{_{a,b\in
\mathbb{N}
}}$ of morphisms $m_{a,b}:A_{a}\otimes A_{b}\rightarrow A_{n}$ such that (%
\ref{form: coro grAlg1}) holds, then $\left( A=\oplus _{n\in \mathbb{N}%
}A_{n},m,u\right) $ is a graded algebra in $\mathcal{M}$.
\end{proposition}

\begin{proof}
$1)$ By definition, we have%
\begin{equation*}
i_{a+b}m_{a,b}=i_{a+b}m_{a+b}\gamma _{a,b}^{A,A}=m\nabla \left[ \left(
i_{u}\otimes i_{v}\right) _{u+v=a+b}\right] \gamma
_{a,b}^{A,A}=m(i_{a}\otimes i_{b}).
\end{equation*}

$2)$ By the universal property of coproducts there is a unique
\begin{equation*}
m_{n}:\oplus _{a+b=n}\left( A_{a}\otimes A_{b}\right) =\left( A\otimes
A\right) _{n}\rightarrow A_{n}
\end{equation*}%
such that $m_{n}\gamma _{a,b}^{A,A}=m_{a,b}$ for every $a,b\in
\mathbb{N}
$ such that $a+b=n.$ For every $n,a,b\in
\mathbb{N}
$ such that $a+b=n$, we have%
\begin{eqnarray*}
&&m\nabla \left[ \left( i_{u}\otimes i_{v}\right) _{u+v=n}\right] \gamma
_{a,b}^{A,A} \\
&=&m\nabla \left[ \left( i_{u}\otimes i_{v}\right) _{u+v=a+b}\right] \gamma
_{a,b}^{A,A}=m(i_{a}\otimes i_{b})\overset{(\ref{form: coro grAlg1})}{=}%
i_{a+b}m_{a,b}=i_{n}m_{n}\gamma _{a,b}^{A,A}
\end{eqnarray*}%
so that $m\nabla \left[ \left( i_{u}\otimes i_{v}\right) _{u+v=n}\right]
=i_{n}m_{n}.$
\end{proof}

\begin{proposition}
\label{lem: graded mij}1) Let $A=\oplus _{n\in \mathbb{N}}A_{n}$ be a graded
object of $\mathcal{M}$ such that there exists a family $\left(
m_{a,b}\right) _{_{a,b\in
\mathbb{N}
}}$
\begin{equation*}
m_{a,b}:A_{a}\otimes A_{b}\rightarrow A_{a+b},
\end{equation*}%
of morphisms and a morphism $u_{0}:\mathbf{1}\rightarrow A_{0}$ which
satisfy
\begin{gather}
m_{a+b,c}\left( m_{a,b}\otimes A_{c}\right) =m_{a,b+c}\left( A_{a}\otimes
m_{b,c}\right) \text{,}  \label{form: locMulti} \\
m_{d,0}\left( A_{d}\otimes u_{0}\right) =r_{A_{d}},\qquad m_{0,d}\left(
u_{0}\otimes A_{d}\right) =l_{A_{d}}\text{,}  \label{form: locUnit}
\end{gather}
for every $a,b,c\in
\mathbb{N}
$. Then there exists a unique morphism $m:A\otimes A\rightarrow A$ such that
(\ref{form: coro grAlg1}) holds.\newline
Moreover $\left( A\oplus _{n\in \mathbb{N}}A_{n},m,u=i_{0}u_{0}\right) $ is
a graded algebra.

$2)$ If $A$ is a graded algebra then (\ref{form: grAlg}) holds, $%
u=i_{0}p_{0}u\ $so that $u$ is a graded homomorphism, and we have that (\ref%
{form: locMulti}) and (\ref{form: locUnit}) hold for every $a,b,c\in
\mathbb{N}
$, where $u_{0}=p_{0}u$.

Moreover $\left( A_{0},m_{0}=m_{0,0},u_{0}=p_{0}u\right) $ is an algebra in $%
\mathcal{M}$, $p_{0}$ is an algebra homomorphism and, for every $n\in
\mathbb{N}
$, $\left( A_{n},m_{0,n},m_{n,0}\right) $ is an $A_{0}$-bimodule such that $%
i_{n}:A_{n}\rightarrow A$ is a morphism of $A_{0}$-bimodules ($A$ is an $%
A_{0}$-bimodule through $i_{0}$).
\end{proposition}

\begin{proof}
$1)$ Since the tensor product commutes with direct sums, we have that $%
A\otimes A=\oplus _{a,b\in
\mathbb{N}
}\left( A_{a}\otimes A_{b}\right) .$ Thus, by the universal property of
coproduct there exists a unique morphism $m:A\otimes A\rightarrow A$ such
that (\ref{form: coro grAlg1}) holds.\newline
The remaining part of the proof is similar to that of Proposition \ref{lem:
graded Deltaij}.

$2)$ It is analogous to Proposition \ref{lem: graded Deltaij}.
\end{proof}

Given an ideal $\left( I,i_{I}\right) $ of an algebra $A$ in $\mathcal{M}$,
one can define as for ordinary algebras, the iterated $n$-th power $I^{n}$
of $I$ in $A$ (see e.g. \cite[Example 3.2]{AMS}).

\begin{definition}
\label{def: strongly grAlg}Let $(A=\oplus _{n\in
\mathbb{N}
}A_{n},m,u)$ be a graded algebra in $\mathcal{M}$. In analogy with the group
graded case, we say that $A$ is a \emph{strongly }$%
\mathbb{N}
$\emph{-graded algebra} whenever

$m_{i,j}^{A}:A_{i}\otimes A_{j}\rightarrow A_{i+j}$ is an epimorphism for
every $i,j\in \mathbb{N},$

where $m_{i,j}^{A}$ is the morphism of Definition \ref{def: grAlg}.
\end{definition}

\begin{proposition}
\label{pro: grSubalg}Let $(A=\oplus _{n\in
\mathbb{N}
}A_{n},m_{A},u_{A})$ and $(B=\oplus _{n\in
\mathbb{N}
}B_{n},m_{B},u_{B})$ be graded algebras in $\mathcal{M}$. Let $%
f:A\rightarrow B$ be a graded algebra homomorphism. Assume that $A$ is a
strongly $%
\mathbb{N}
$-graded algebra.

Then the following assertions are equivalent

$\left( a\right) $ $f_{0}$ and $f_{1}$ are epimorphisms and $B$ is is a
strongly $%
\mathbb{N}
$-graded algebra.

$\left( a^{\prime }\right) $ $f_{0}$ and $f_{1}$ are epimorphisms and $%
m_{a,1}^{B}:B_{a}\otimes B_{1}\rightarrow B_{a+1}$ is an epimorphism for
every $a\in
\mathbb{N}
$.

$\left( b\right) $ $f_{n}$ is an epimorphism for every $n\in
\mathbb{N}
$.

$\left( c\right) $ $f$ is an epimorphism.

$\left( d\right) $ $f_{0}$ and $f_{1}$ are epimorphisms and $(B\left[ n%
\right] ,\nu _{n}^{B})=B\left[ 1\right] ^{n},$ for every $n\in \mathbb{N}.$

$\left( e\right) $ $f_{0}$ and $f_{1}$ are epimorphisms and $(B\left[ 2%
\right] ,\nu _{2}^{B})=B\left[ 1\right] ^{2}$.
\end{proposition}

\begin{proof}
It is analogous to that of Proposition \ref{pro: grSubcoalg}.
\end{proof}

\begin{claim}
\textbf{The tensor algebra.} \label{claim: tensor alg} For any object $X\in
\mathcal{M},$ we can define the tensor algebra $T\left( X\right) $ as
follows. As an object in $\mathcal{M}$ it is defined by setting $\ $
\begin{equation*}
T=T\left( X\right) :=\oplus _{p\in \mathbb{N}}X^{\otimes p}.
\end{equation*}
Let $i_{p}:X^{\otimes p}\rightarrow T(X)$ be the canonical injection and let
$u_{0}:=\mathrm{Id}_{\mathbf{1}}.$

Let%
\begin{equation*}
m_{p,q}:X^{\otimes p}\otimes X^{\otimes q}\rightarrow X^{\otimes p+q},\text{
for every }p,q\in \mathbb{N}
\end{equation*}%
be the canonical isomorphism, which is unique by Coherence Theorem. \newline
Still by Coherence Theorem and the definitions of $m_{d,0}$ and $m_{0,d}$,
one has%
\begin{gather}
m_{p+q,r}[m_{p,q}\otimes X^{\otimes r}]=m_{p,q+r}[X^{\otimes p}\otimes
m_{q,r}].  \label{ec:tensor algebra} \\
m_{d,0}\left( X^{\otimes d}\otimes u_{0}\right) =r_{X^{\otimes d}},\qquad
m_{0,d}\left( u_{0}\otimes X^{\otimes d}\right) =l_{X^{\otimes d}}  \notag
\end{gather}

so that (\ref{form: locMulti}) and (\ref{form: locUnit}) hold.

By Proposition \ref{lem: graded mij}, there exists a unique morphism $%
m_{T}:T\otimes T\rightarrow T$ such that (\ref{form: coro grAlg1}) holds,
i.e.%
\begin{equation}
m_{T}\circ (i_{p}\otimes i_{q})=i_{p+q}\circ m_{p,q},\text{ for every }%
p,q\in \mathbb{N}.  \label{eq: multTensor}
\end{equation}

Moreover $\left( T,m_{T},u_{T}=i_{0}u_{0}=i_{0}\right) $ is a graded algebra
with graduation defined by $\left( X^{\otimes p}\right) _{p\in \mathbb{N}}$.
Note that
\begin{equation*}
m_{T}:=\oplus _{p\in \mathbb{N}}\left( \nabla \lbrack
(m_{i,j})_{i+j=p}]\right) :T(X)\otimes T(X)\rightarrow T(X),
\end{equation*}%
where $\nabla \lbrack (m_{i,j})_{i+j=p}]:\oplus _{i+j=p}\left( X^{\otimes
i}\otimes X^{\otimes j}\right) \longrightarrow X^{\otimes p}$ denotes the
codiagonal morphism associated to the family $(m_{i,j})_{i+j=p}.$ Let $%
(A,m_{A},u_{A})$ be an algebra in $\mathcal{M}.$\newline
As explained in Section \ref{sec: Preliminaries and Notations}, we have that
$({_{A}\mathcal{M}_{A}},\otimes _{A},A)$ is an abelian monoidal category.
Furthermore one can see that $_{A}\mathcal{M}_{A}$ has also arbitrary direct
sums which commute with $\otimes _{A}$. \newline
Therefore we can consider, in the monoidal category $(_{A}\mathcal{M}%
_{A},\otimes _{A},A),$ the tensor algebra of an arbitrary $A$-bimodule $M.$
We will denote it by
\begin{equation*}
\left( T=T_{A}(M),\overline{m}_{T},\overline{u_{T}}\right) .
\end{equation*}%
Note that%
\begin{equation*}
T=\oplus _{p\in \mathbb{N}}M^{\otimes _{A}p},\qquad \overline{m}%
_{T}:T\otimes _{A}T\rightarrow T\qquad \text{and}\qquad \overline{u_{T}}%
:A\rightarrow T.
\end{equation*}%
Set%
\begin{equation*}
m_{T}=\overline{m}_{T}\circ {_{A}\chi _{T,T}}:T\otimes T\rightarrow T\qquad
\text{and}\qquad u_{T}=\overline{u_{T}}\circ u_{A}:\mathbf{1}\rightarrow T,
\end{equation*}%
where ${_{A}\chi _{T,T}}:T\otimes T\longrightarrow T\otimes _{A}T$ is the
canonical morphism introduced in Section \ref{sec: Preliminaries and
Notations}.
\end{claim}

We have the following result.

\begin{proposition}
\label{pro: tensGrad}Let $\left( A,m_{A},u_{A}\right) $ be an algebra in $%
\mathcal{M}$ and let $(M,\mu _{M}^{r},\mu _{M}^{l})$ be an $A$-bimodule. Let
$T=T_{A}(M)$ be the tensor algebra. Then $\left( T=\oplus _{n\in \mathbb{N}%
}T_{n},m_{T},u_{T}\right) $ is a graded algebra in $(\mathcal{M},\otimes ,%
\mathbf{1})$, where $T_{n}=M^{\otimes _{A}n}.$ Moreover we have%
\begin{equation*}
m_{u,v}=\left\{
\begin{tabular}{ll}
$(M^{\otimes _{A}u-1}\otimes _{A}\chi _{M}\otimes _{A}M^{\otimes _{A}v-1}),$
& $\text{for any }u,v\geq 1;$ \\
$(M^{\otimes _{A}u-1}\otimes _{A}\mu _{M}^{r}),$ & $\text{for any }u\geq
1,v=0;$ \\
$(\mu _{M}^{l}\otimes _{A}M^{\otimes _{A}v-1}),$ & $\text{for any }u=0,v\geq
1;$ \\
$m_{A},$ & $\text{for any }u=0,v=0.$%
\end{tabular}%
\right.
\end{equation*}%
In particular $T$ is a strongly $%
\mathbb{N}
$-graded algebra.
\end{proposition}

\begin{proof}
Follows by (\ref{claim: tensor alg}) and by Proposition \ref{coro: grAlg} as
(\ref{form: coro grAlg1}) is clearly satisfied.
\end{proof}

We are now able to state the Universal property of the relative tensor
algebra.

\begin{theorem}[Universal property of the relative tensor algebra]
\label{teo: univ property of tensor algebra} Let $(\mathcal{M},\otimes ,%
\mathbf{1})$ be a cocomplete abelian monoidal category such that the tensor
product commutes with direct sums.\newline
Let $A,B$ be algebras in $\mathcal{M}$ and let $f_{A}:A\rightarrow B$ be an
algebra homomorphism. \newline
Let $M\in {_{A}\mathcal{M}_{A}}$, and let $f_{M}:M\rightarrow B$ be a
morphism in ${_{A}\mathcal{M}_{A}},$ where $B\in {_{A}\mathcal{M}_{A}}$ via $%
f_{0}.$\newline
Then there is a unique algebra homomorphism $f:T_{A}(M)\rightarrow B$ such
that
\begin{equation*}
\diagUnivTens%
\end{equation*}%
where $i_{0}:A\rightarrow T_{A}(M)$ and $i_{1}:M\rightarrow T_{A}(M)$ are
the canonical injections. \newline
Moreover%
\begin{equation*}
f\circ i_{t}^{T}=\overline{m}_{B}^{t-1}\circ \left( f_{M}\right) ^{\otimes
_{A}t}\text{ for every }t\in
\mathbb{N}
,
\end{equation*}%
where $\overline{m}_{B}^{t-1}:B^{\otimes _{A}t}\rightarrow B$ is the $t$-th
iterated multiplication of $B$.
\end{theorem}

When $B$ is a graded algebra the universal property of $T_{A}(M)$ has a
graded version which is the following.

\begin{theorem}
\label{coro: GRuniv property of tensor algebra}Let $(\mathcal{M},\otimes ,%
\mathbf{1})$ be a cocomplete abelian monoidal category such that the tensor
product commutes with direct sums. Let $A$ be an algebra in $\mathcal{M}$
and let $M\in {_{A}\mathcal{M}_{A}}$.\newline
Let $B$ be a graded algebra, let $g_{A}:A\rightarrow B_{0}$ be an algebra
homomorphism and let $g_{M}:M\rightarrow B_{1}$ be a morphism in ${_{A}%
\mathcal{M}_{A}},$ where $B_{1}\in {_{A}\mathcal{M}_{A}}$ via $g_{A}.$%
\newline
Then there is a unique algebra homomorphism $f:T_{A}(M)\rightarrow B$ such
that
\begin{equation*}
f\circ i_{0}^{T}=i_{0}^{B}\circ g_{A}\text{\qquad and\qquad }f\circ
i_{1}^{T}=i_{1}^{B}\circ g_{M}.
\end{equation*}%
Moreover $f$ is a morphism of graded algebras where%
\begin{equation*}
f_{t}=p_{t}^{B}\circ \overline{m}_{B}^{t-1}\circ \left( i_{1}^{B}\circ
g_{M}\right) ^{\otimes _{A}t}\text{ for every }t\in
\mathbb{N}%
\end{equation*}%
and $\overline{m}_{B}^{t-1}:B^{\otimes _{A}t}\rightarrow B$ is the $t$-th
iterated multiplication of $B$.
\end{theorem}

\begin{proof}
It is analogous to that of Theorem \ref{coro: GRuniv property of cotensor
coalgebra}.
\end{proof}

\begin{theorem}
\label{teo: phi epi}Let $(A=\oplus _{n\in
\mathbb{N}
}A_{n},m_{A},u_{A})$ be a graded algebra in a cocomplete abelian monoidal
category $\mathcal{M}.$ Assume that the tensor product commutes with direct
sums. \newline
Let $T:=T_{A_{0}}(A_{1})$ be the tensor algebra. Then there is a unique
algebra homomorphism
\begin{equation*}
\varphi :T_{A_{0}}(A_{1})\rightarrow A,
\end{equation*}%
such that $\varphi \circ i_{0}^{T}=i_{0}^{A}$ and $\varphi \circ
i_{1}^{T}=i_{1}^{A}$. \newline
Moreover $\varphi $ is a graded algebra homomorphism with%
\begin{equation*}
\varphi _{t}=p_{t}^{A}\circ \overline{m}_{A}^{t-1}\circ \left(
i_{1}^{A}\right) ^{\otimes _{A_{0}}t}\text{ for every }t\in
\mathbb{N}%
\end{equation*}%
and the following assertions are equivalent.

$\left( a\right) $ $A$ is a strongly $%
\mathbb{N}
$-graded algebra.

$\left( a^{\prime }\right) $ $m_{a,1}^{A}:A_{a}\otimes A_{1}\rightarrow
A_{a+1}$ is an epimorphism for every $a\in
\mathbb{N}
$.

$\left( b\right) \ \varphi _{n}$ is an epimorphism for every $n\in
\mathbb{N}
$.

$\left( c\right) $ $\varphi $ is an epimorphism.

$\left( d\right) $ $(A\left[ n\right] ,\nu _{n}^{A})=A\left[ 1\right] ^{n},$
for every $n\in \mathbb{N}.$

$\left( e\right) $ $(A\left[ 2\right] ,\nu _{2}^{A})=A\left[ 1\right] ^{2}$.
\end{theorem}

\begin{proof}
It is analogous to that of Theorem \ref{teo: psi mono}.
\end{proof}

\section{Braided Bialgebra Structure of the Cotensor Coalgebra}

The main aim of this section is to provide a braided bialgebra structure for
the cotensor coalgebra inside a braided monoidal category. This structure is
used to extend the notion of bialgebra of type one, introduced in the
classical case by Nichols in \cite{Ni}, to the wider context of a braided
monoidal category (see Definition \ref{def: type1}). A universal property
for the cotensor bialgebra is also proven (see Theorem \ref{teo: univ
property of cotensor bialgebra}).

\begin{claim}
\textbf{Braided bialgebras.} \label{cl: brdBialg} A \emph{braided bialgebra}
in $(\mathcal{M},c)$, is a sextuple $(H,m,u,\Delta ,\varepsilon )$ such that
$(H,m,u)$ is an algebra in $\mathcal{M}$, $(H,\Delta ,\varepsilon )$ is a
coalgebra in $\mathcal{M}$ and this two structure are compatible in the
sense that
\begin{equation}
\Delta \circ m=\left( m\otimes m\right) \circ \left( H\otimes c\otimes
H\right) \circ \left( \Delta \otimes \Delta \right) \qquad \text{and}\qquad
m_{\mathbf{1}}\circ \left( \varepsilon \otimes \varepsilon \right)
=\varepsilon \circ m.  \label{form: def Braided}
\end{equation}
\end{claim}

\begin{definitions}
Let $(B,m_{B},u_{B},\Delta _{B},\varepsilon _{B})$ be a braided bialgebra in
a cocomplete braided monoidal category $\left( \mathcal{M},c\right) $. Then $%
B$ is called a \emph{graded braided bialgebra} if $B$ is a graded object
with graduation defined by $\left( B_{k}\right) _{k\in \mathbb{N}}$ and if,
with respect to this graduation, $\left( B,m_{B},u_{B}\right) $ is a graded
algebra and $(B,\Delta _{B},\varepsilon _{B})$ is a graded coalgebra. A
morphism $f:A\rightarrow B$ between graded braided bialgebras is called a
\emph{morphism of graded (braided) bialgebras} whenever it is both a
morphism of graded algebras and a morphism of graded coalgebras i.e. if it
is a bialgebra homomorphism which is also a graded homomorphism.
\end{definitions}

\begin{definition}
\label{def: HopfBimod} Let $H$ be a braided bialgebra in $\left( \mathcal{M}%
,c\right) .$

An object in $_{H}^{H}\mathcal{M}_{H}^{H}$ is a $5$-tuple $(M,\mu
_{M}^{r},\mu _{M}^{l},\rho _{M}^{r},\rho _{M}^{l})$ such that

\begin{itemize}
\item $\left( M,\mu _{M}^{r},\mu _{M}^{l}\right) $ is an $H$-bimodule;

\item $\left( M,\rho _{M}^{r},\rho _{M}^{l}\right) $ is an $H$-bicomodule;

\item the following compatibility conditions are fulfilled:%
\begin{eqnarray}
\rho _{M}^{l}\mu _{M}^{l} &=&\left( m_{H}\otimes \mu _{M}^{l}\right) \left(
H\otimes c_{H,H}\otimes M\right) \left( \Delta _{H}\otimes \rho
_{M}^{l}\right) ,  \label{form: HopfBimod 1} \\
\rho _{M}^{l}\mu _{M}^{r} &=&\left( m_{H}\otimes \mu _{M}^{r}\right) \left(
H\otimes c_{M,H}\otimes H\right) \left( \rho _{M}^{l}\otimes \Delta
_{H}\right) ,  \label{form: HopfBimod 2} \\
\rho _{M}^{r}\mu _{M}^{l} &=&\left( \mu _{M}^{l}\otimes m_{H}\right) \left(
H\otimes c_{H,M}\otimes H\right) \left( \Delta _{H}\otimes \rho
_{M}^{r}\right) ,  \label{form: HopfBimod 3} \\
\rho _{M}^{r}\mu _{M}^{r} &=&\left( \mu _{M}^{r}\otimes m_{H}\right) \left(
M\otimes c_{H,H}\otimes H\right) \left( \rho _{M}^{r}\otimes \Delta
_{H}\right) .  \label{form: HopfBimod 4}
\end{eqnarray}
\end{itemize}
\end{definition}

\begin{proposition}
\label{pro: BotB is graded}Let $\left( \mathcal{M},c\right) $ be a
cocomplete coabelian braided monoidal category and let $C=\oplus _{n\in
\mathbb{%
\mathbb{N}
}}C_{n}$ be a graded coalgebra in $\mathcal{M}$. Assume that the tensor
product commutes with direct sums. Then $\left( C\otimes C,\Delta _{C\otimes
C},\varepsilon _{C\otimes C}\right) $ is a graded coalgebra where%
\begin{equation*}
\Delta _{C\otimes C}:C\otimes C\overset{\Delta _{C}\otimes \Delta _{C}}{%
\longrightarrow }C\otimes C\otimes C\otimes C\overset{C\otimes
c_{C,C}\otimes C}{\longrightarrow }C\otimes C\otimes C\otimes C,\quad
\varepsilon _{C\otimes C}:C\otimes C\overset{\varepsilon _{C}\otimes
\varepsilon _{C}}{\longrightarrow }\mathbf{1}\otimes \mathbf{1}\overset{m_{%
\mathbf{1}}}{\longrightarrow }\mathbf{1}\text{.}
\end{equation*}%
and with graduation given by $\left( C\otimes C\right)
_{n}=\bigoplus_{a+b=n}C_{a}\otimes C_{b}.$

Moreover, for every $s,t\in
\mathbb{N}
$ we have%
\begin{equation}
\Delta _{s,t}^{C\otimes C}:=\nabla \left[ \sum_{\substack{ i+j=a,u+v=b  \\ %
i+u=s,j+v=t}}\left( \gamma _{i,u}^{C,C}\otimes \gamma _{j,v}^{C,C}\right)
\left( C_{i}\otimes c_{C_{j},C_{u}}\otimes C_{v}\right) \left( \Delta
_{i,j}\otimes \Delta _{u,v}\right) \right] _{\substack{ a,b,  \\ a+b=s+t}}.
\label{form: DeltaCotC}
\end{equation}
\end{proposition}

\begin{proof}
It is well known that $\left( C\otimes C,\Delta _{C\otimes C},\varepsilon
_{C\otimes C}\right) $ is a coalgebra (it is dual to \cite[Lemma 9.2.12,
page 438]{Majid}). Let us check the part of the statement concerning the
graduation. In view of Proposition \ref{pro: grCoalg}, since $C$ is a graded
coalgebra, for every $n\in
\mathbb{N}
$ there exists a family $\left( \Delta _{a,b}\right) _{_{a+b=n}}$ of
morphisms $\Delta _{a,b}:C_{n}\rightarrow C_{a}\otimes C_{b},$ such that (%
\ref{form: grCoalg}) holds true. For every $s,t\in
\mathbb{N}
$ define $\Delta _{s,t}^{C\otimes C}:\left( C\otimes C\right)
_{s+t}\rightarrow \left( C\otimes C\right) _{s}\otimes \left( C\otimes
C\right) _{t}$ accordingly to (\ref{form: DeltaCotC}). It is straightforward
to prove that%
\begin{equation*}
\Delta _{C\otimes C}\circ i_{n}^{C\otimes C}=\sum\limits_{s+t=n}\left(
i_{s}^{C\otimes C}\otimes i_{t}^{C\otimes C}\right) \Delta _{s,t}^{C\otimes
C}.
\end{equation*}
\end{proof}

\begin{theorem}
\label{teo: cotensor bial} Let $H$ be a braided bialgebra in a cocomplete
and complete coabelian braided monoidal category $\left( \mathcal{M}%
,c\right) $ satisfying $AB5.$ Assume that the tensor product commutes with
direct sums. \newline
Let $(M,\mu _{M}^{r},\mu _{M}^{l},\rho _{M}^{r},\rho _{M}^{l})$ be in $%
_{H}^{H}\mathcal{M}_{H}^{H}$. Let $T^{c}=T_{H}^{c}(M)$ be the cotensor
coalgebra. Then there are unique coalgebra homomorphisms%
\begin{equation*}
m_{T^{c}}:T^{c}\otimes T^{c}\rightarrow T^{c}\qquad \text{and}\qquad
u_{T^{c}}:\mathbf{1}\rightarrow T^{c}
\end{equation*}%
such that the diagrams%
\begin{equation}
\diagComultil%
\text{ \qquad }%
\diagComultir
\label{form: cotmulti}
\end{equation}%
are commutative, where $p_{n}:T^{c}\rightarrow M^{\square _{H}n}$ denotes
the canonical projection. Moreover $\left( T^{c},m_{T^{c}},u_{T^{c}},\Delta
_{T^{c}},\varepsilon _{T^{c}}\right) $ is a graded braided bialgebra in $%
\mathcal{M}$.
\end{theorem}

\begin{proof}
First of all, by Proposition \ref{pro: BotB is graded}, we have that $\left(
E=T^{c}\otimes T^{c},\Delta _{E},\varepsilon _{E}\right) $ is a graded
coalgebra where%
\begin{eqnarray*}
\Delta _{E} &:&T^{c}\otimes T^{c}\overset{\Delta _{T^{c}}\otimes \Delta
_{T^{c}}}{\longrightarrow }T^{c}\otimes T^{c}\otimes T^{c}\otimes T^{c}%
\overset{T^{c}\otimes c_{T^{c},T^{c}}\otimes T^{c}}{\longrightarrow }%
T^{c}\otimes T^{c}\otimes T^{c}\otimes T^{c}, \\
\varepsilon _{E} &:&T^{c}\otimes T^{c}\overset{\varepsilon _{T^{c}}\otimes
\varepsilon _{T^{c}}}{\longrightarrow }\mathbf{1}\otimes \mathbf{1}\overset{%
m_{\mathbf{1}}}{\longrightarrow }\mathbf{1}\text{.}
\end{eqnarray*}%
Since $T^{c}=\oplus _{n\in \mathbb{N}}T_{n}^{c}$ is a graded coalgebra, we
have that (\ref{form: coro grCoalg1}) holds so that, by Proposition \ref%
{pro: cotGrad}, we have%
\begin{equation*}
(p_{0}\otimes p_{0})\Delta _{T^{c}}=\Delta _{H}p_{0}\qquad \text{and}\qquad
(p_{0}\otimes p_{1})\Delta _{T^{c}}=\rho _{M}^{l}p_{1}.
\end{equation*}%
Set%
\begin{equation*}
g_{H}:=m_{H}:E_{0}=H\otimes H\rightarrow H\qquad \text{and}\qquad
g_{M}:=\nabla \left[ \left( \mu _{M}^{l},\mu _{M}^{r}\right) \right]
:E_{1}=\left( H\otimes M\right) \oplus \left( M\otimes H\right) \rightarrow
M.
\end{equation*}%
Then $g_{H}$ is a coalgebra homomorphism:%
\begin{equation*}
\Delta _{H}g_{H}=\Delta _{H}m_{H}=\left( m_{H}\otimes m_{H}\right) \left(
H\otimes c_{H,H}\otimes H\right) \left( \Delta _{H}\otimes \Delta
_{H}\right) =\left( g_{H}\otimes g_{H}\right) \Delta _{H\otimes H}.
\end{equation*}%
Moreover $g_{M}$ is a morphism of left $H$-comodules%
\begin{eqnarray*}
&&\left( g_{H}\otimes g_{M}\right) \Delta _{0,1}^{T^{c}\otimes T^{c}}\overset%
{(\ref{form: DeltaCotC})}{=}\left( g_{H}\otimes g_{M}\right) \nabla \left[
\begin{array}{c}
\left( \gamma _{0,0}^{T^{c},T^{c}}\otimes \gamma _{0,1}^{T^{c},T^{c}}\right)
\left( T_{0}^{c}\otimes c_{T_{0}^{c},T_{0}^{c}}\otimes T_{1}^{c}\right)
\left( \Delta _{0,0}\otimes \Delta _{0,1}\right) , \\
\left( \gamma _{0,0}^{T^{c},T^{c}}\otimes \gamma _{1,0}^{T^{c},T^{c}}\right)
\left( T_{0}^{c}\otimes c_{T_{1}^{c},T_{0}^{c}}\otimes T_{0}^{c}\right)
\left( \Delta _{0,1}\otimes \Delta _{0,0}\right)%
\end{array}%
\right] \\
&=&\nabla \left[
\begin{array}{c}
\left( m_{H}\otimes \mu _{M}^{l}\right) \left( H\otimes c_{H,H}\otimes
M\right) \left( \Delta _{H}\otimes \rho _{M}^{l}\right) , \\
\left( m_{H}\otimes \mu _{M}^{r}\right) \left( H\otimes c_{M,H}\otimes
H\right) \left( \rho _{M}^{r}\otimes \Delta _{H}\right)%
\end{array}%
\right] \overset{\text{(\ref{form: HopfBimod 1}),(\ref{form: HopfBimod 2})}}{%
=}\nabla \left[ \rho _{M}^{l}\mu _{M}^{l},\rho _{M}^{l}\mu _{M}^{r}\right]
=\rho _{M}^{l}g_{M}
\end{eqnarray*}%
Analogously $\rho _{M}^{r}g_{M}=\left( g_{M}\otimes g_{H}\right) \Delta
_{1,0}^{T^{c}\otimes T^{c}},$ i.e. $g_{M}$ is a morphism of right $H$%
-comodules and hence a morphism of $H$-bicomodules. By applying Theorem \ref%
{coro: GRuniv property of cotensor coalgebra} to the case $"C"=H$ and $%
"B"=H\otimes H$, using the maps $g_{H}$ and $g_{M}$ above, we get a unique
coalgebra homomorphism $m_{T^{c}}=f:E\rightarrow T^{c}$ such that the left
side of (\ref{form: cotmulti}) is commutative so that
\begin{equation*}
p_{0}m_{T^{c}}\left( m_{T^{c}}\otimes T^{c}\right) =p_{0}m_{T^{c}}\left(
T^{c}\otimes m_{T^{c}}\right) \qquad \text{and}\qquad p_{1}m_{T^{c}}\left(
m_{T^{c}}\otimes T^{c}\right) =p_{1}m_{T^{c}}\left( T^{c}\otimes
m_{T^{c}}\right) .
\end{equation*}%
Since $T^{c}\otimes T^{c}\otimes T^{c}$ is a coalgebra and $m_{T^{c}}\left(
m_{T^{c}}\otimes T^{c}\right) $ and $m_{T^{c}}\left( T^{c}\otimes
m_{T^{c}}\right) $ are coalgebra homomorphisms, then by Corollary \ref{coro:
uniqueness}, we obtain
\begin{equation*}
m_{T^{c}}\left( m_{T^{c}}\otimes T^{c}\right) =m_{T^{c}}\left( T^{c}\otimes
m_{T^{c}}\right) .
\end{equation*}%
Set%
\begin{equation*}
g_{H}^{\prime }:=u_{H}:\mathbf{1}\rightarrow H\qquad \text{and}\qquad
g_{M}^{\prime }:=0:0\rightarrow M.
\end{equation*}%
Then $g_{H}^{\prime }$ is a coalgebra homomorphism and $g_{M}^{\prime }$ is
a morphism of $H$-bicomodules. By applying Theorem \ref{coro: GRuniv
property of cotensor coalgebra} to the case $B=\mathbf{1}$, we get a unique
coalgebra homomorphism $u_{T^{c}}=f^{\prime }:\mathbf{1}\rightarrow T^{c}$
such that the right side of \ref{form: cotmulti} is commutative so that%
\begin{eqnarray*}
p_{0}m_{T^{c}}\left( u_{T^{c}}\otimes T^{c}\right) &=&m_{H}\left(
u_{H}\otimes H\right) \left( \mathbf{1}\otimes p_{0}\right) =l_{H}\left(
\mathbf{1}\otimes p_{0}\right) =p_{0}l_{T^{c}} \\
p_{1}m_{T^{c}}\left( u_{T^{c}}\otimes T^{c}\right) &=&\mu _{M}^{l}\left(
u_{H}\otimes M\right) \left( \mathbf{1}\otimes p_{1}\right) =l_{M}\left(
\mathbf{1}\otimes p_{1}\right) =p_{1}l_{T^{c}}.
\end{eqnarray*}%
Since $\mathbf{1}\otimes T^{c}$ is a coalgebra and $m_{T^{c}}\left(
u_{T^{c}}\otimes T^{c}\right) $ and $r_{T^{c}}$ are coalgebra homomorphisms,
then by Corollary \ref{coro: uniqueness}, we have $ m_{T^{c}}\left(
u_{T^{c}}\otimes T^{c}\right) =l_{T^{c}}. $ Analogously $m_{T^{c}}\left(
T^{c}\otimes u_{T^{c}}\right) =r_{T^{c}}.$ Now, in view of Theorem \ref%
{coro: GRuniv property of cotensor coalgebra}, $m_{T^{c}}$ and $u_{T^{c}}$
are graded homomorphisms so that $\left( T^{c},m_{T^{c}},u_{T^{c}},\Delta
_{T^{c}},\varepsilon _{T^{c}}\right) $ is a graded braided bialgebra in $%
\mathcal{M}$.
\end{proof}

\begin{remark}
Let $M$ be an object in a cocomplete and complete coabelian braided monoidal
category $\left( \mathcal{M},c\right) $ satisfying $AB5.$ Assume that the
tensor product commutes with direct sums. \newline
By applying Theorem \ref{teo: cotensor bial} to the case $H=\mathbf{1}$ we
endow the cotensor coalgebra $T^c=T_{H}^{c}(M)$ with an algebra structure
such that $T^c$ becomes a braided bialgebra. This algebra structure is the
braided analogue of the so called "Shuffle Algebra" in the category of
vector spaces.
\end{remark}

\begin{theorem}
\label{teo: univ property of cotensor bialgebra} Let $H$ be a braided
bialgebra in a cocomplete and complete coabelian braided monoidal category $%
\left( \mathcal{M},c\right) $ satisfying $AB5.$ Assume that the tensor
product commutes with direct sums. \newline
Let $(M,\mu _{M}^{r},\mu _{M}^{l},\rho _{M}^{r},\rho _{M}^{l})$ be in $%
_{H}^{H}\mathcal{M}_{H}^{H}$. Let $T^{c}=T_{H}^{c}(M)$ be the cotensor
coalgebra.\newline
Let $\delta :D\rightarrow E$ be a monomorphism which is a homomorphism of
coalgebras. Assume that there exist morphisms
\begin{equation*}
m_{\widetilde{D}}:\widetilde{D}\otimes \widetilde{D}\rightarrow \widetilde{D}%
\qquad \text{and}\qquad u_{\widetilde{D}}:\mathbf{1\rightarrow }\widetilde{D}
\end{equation*}%
such that $\left( \widetilde{D},m_{\widetilde{D}},u_{\widetilde{D}},\Delta _{%
\widetilde{D}},\varepsilon _{\widetilde{D}},\right) $ is a braided bialgebra
in $\mathcal{M}$.\newline
Let $f_{H}:\widetilde{D}\rightarrow H$ be a bialgebra homomorphism and let $%
f_{M}:\widetilde{D}\rightarrow M$ be a morphism of $H$-bicomodules such that
$f_{M}\xi _{1}=0$, where $\widetilde{D}$ is a bicomodule via $f_{H}.$ Assume
that%
\begin{equation}
f_{M}m_{\widetilde{D}}=\mu _{M}^{l}\left( f_{H}\otimes f_{M}\right) +\mu
_{M}^{r}\left( f_{M}\otimes f_{H}\right)  \label{form: UnivCotBialg}
\end{equation}%
(i.e. $f_{M}$ is a derivation of $\widetilde{D}$ with values in the $%
\widetilde{D}$-bimodule $M$, where $M$ is regarded as a bimodule via $f_{H}$%
). Then there is a unique coalgebra homomorphism $f:\widetilde{D}\rightarrow
T_{H}^{c}(M)$ such that $p_{0}f=f_{H}$ and $p_{1}f=f_{M}$, where $%
p_{n}:T_{C}^{c}(M)\rightarrow M^{\square n}$ denotes the canonical
projection. 
\begin{equation*}
\diagUnivCotensorBialg%
\end{equation*}%
Moreover $f$ is a bialgebra homomorphism.
\end{theorem}

\begin{proof}
By Theorem \ref{coro: univ property of cotensor coalgebra}, there is a
unique coalgebra homomorphism $f:\widetilde{D}\rightarrow T_{H}^{c}(M)$ such
that $p_{0}f=f_{H}$ and $p_{1}f=f_{M}$. We have
\begin{eqnarray*}
p_{0}m_{T^{c}}\left( f\otimes f\right) \overset{\text{(\ref{form: cotmulti})}%
}{=}m_{H}\left( f_{H}\otimes f_{H}\right) &=&f_{H}m_{\widetilde{D}}=p_{0}fm_{%
\widetilde{D}}, \\
p_{1}m_{T^{c}}\left( f\otimes f\right) \overset{\text{(\ref{form: cotmulti})}%
}{=}\mu _{M}^{l}\left( f_{H}\otimes f_{M}\right) +\mu _{M}^{r}\left(
f_{M}\otimes f_{H}\right) &=&f_{M}m_{\widetilde{D}}=p_{1}fm_{\widetilde{D}}.
\end{eqnarray*}%
We also obtain%
\begin{eqnarray*}
f_{M}u_{\widetilde{D}} &=&f_{M}m_{\widetilde{D}}\left( u_{\widetilde{D}%
}\otimes u_{\widetilde{D}}\right) m_{\mathbf{1}}^{-1}\overset{\text{(\ref%
{form: UnivCotBialg})}}{=}\mu _{M}^{l}\left( f_{H}u_{\widetilde{D}}\otimes
f_{M}u_{\widetilde{D}}\right) m_{\mathbf{1}}^{-1}+\mu _{M}^{r}\left( f_{M}u_{%
\widetilde{D}}\otimes f_{H}u_{\widetilde{D}}\right) m_{\mathbf{1}}^{-1} \\
&=&\mu _{M}^{l}\left( u_{H}\otimes f_{M}u_{\widetilde{D}}\right) m_{\mathbf{1%
}}^{-1}+\mu _{M}^{r}\left( f_{M}u_{\widetilde{D}}\otimes u_{H}\right) m_{%
\mathbf{1}}^{-1}=f_{M}u_{\widetilde{D}}+f_{M}u_{\widetilde{D}}
\end{eqnarray*}%
so that $f_{M}u_{\widetilde{D}}=0\ $and hence, by Theorem \ref{coro: univ
property of cotensor coalgebra},%
\begin{equation*}
p_{0}fu_{\widetilde{D}}=f_{H}u_{\widetilde{D}}=u_{H}=p_{0}u_{T^{c}}\qquad
\text{and}\qquad p_{1}fu_{\widetilde{D}}=f_{M}u_{\widetilde{D}%
}=0=p_{1}u_{T^{c}}.
\end{equation*}%
Since $m_{T^{c}}\left( f\otimes f\right) ,fm_{\widetilde{D}}:\widetilde{D}%
\otimes \widetilde{D}\rightarrow T^{c}$ and $fu_{\widetilde{D}},u_{T^{c}}:%
\mathbf{1}\rightarrow T^{c}$ are coalgebra homomorphisms and since, for $%
i=0,1,$
\begin{equation*}
p_{i}m_{T^{c}}\left( f\otimes f\right) =p_{i}fm_{\widetilde{D}}\qquad \text{%
and}\qquad p_{i}fu_{\widetilde{D}}=p_{i}u_{T^{c}},
\end{equation*}
then, by Corollary \ref{coro: uniqueness}, we get that $m_{T^{c}}\left(
f\otimes f\right) =fm_{\widetilde{D}}$ and $fu_{\widetilde{D}}=u_{T^{c}}$
i.e. that $f$ is an algebra homomorphism.
\end{proof}

\section{Graded Braided Bialgebra Structure of the Tensor Algebra \label%
{sec: braidtensor}}

In this section $\left( \left( \mathcal{M},\otimes ,\mathbf{1}\right)
,c\right) $ will denote a cocomplete abelian braided monoidal category such
that the tensor product commutes with direct sums.

Next aim is to provide a braided bialgebra structure (see \ref{cl: brdBialg}%
) for the tensor algebra inside a braided monoidal category.

\begin{proposition}
\label{pro: AotA is graded}Let $A=\oplus _{n\in \mathbb{%
\mathbb{N}
}}A_{n}$ be a graded algebra in $\mathcal{M}$. Then $\left( A\otimes
A,m_{A\otimes A},u_{A\otimes A}\right) $ is a graded algebra where%
\begin{equation*}
m_{A\otimes A}:A\otimes A\otimes A\otimes A\overset{A\otimes c_{A,A}\otimes A%
}{\longrightarrow }A\otimes A\otimes A\otimes A\overset{m\otimes m}{%
\longrightarrow }A\otimes A,\quad u_{A\otimes A}:\mathbf{1}\overset{m_{%
\mathbf{1}}^{-1}}{\longrightarrow }\mathbf{1}\otimes \mathbf{1}\overset{%
u\otimes u}{\longrightarrow }A\otimes A\text{.}
\end{equation*}%
and with graduation given by $\left( A\otimes A\right)
_{n}=\bigoplus_{a+b=n}A_{a}\otimes A_{b}.$
\end{proposition}

\begin{proof}
It is similar to the proof of Proposition \ref{pro: BotB is graded}.
\end{proof}

\begin{theorem}
\label{teo: tensor bial} Let $H$ be a braided bialgebra in a cocomplete
abelian braided monoidal category $\left( \mathcal{M},c\right) $. Assume
that the tensor product commutes with direct sums. \newline
Let $(M,\mu _{M}^{r},\mu _{M}^{l},\rho _{M}^{r},\rho _{M}^{l})$ be in $%
_{H}^{H}\mathcal{M}_{H}^{H}$. Let $T=T_{H}(M)$ be the tensor algebra. Then
there are unique coalgebra homomorphisms%
\begin{equation*}
\Delta _{T}:T\rightarrow T\otimes T\qquad \text{and}\qquad \varepsilon
_{T}:T\rightarrow \mathbf{1}
\end{equation*}%
such that the diagrams%
\begin{equation}
\diagMultil%
\text{ \qquad }%
\diagMultir
\label{form: multi}
\end{equation}%
are commutative, where $i_{n}:M^{\otimes _{H}n}\rightarrow T$ denotes the
canonical injection. Moreover $\left( T,m_{T},u_{T},\Delta _{T},\varepsilon
_{T}\right) $ is a graded braided bialgebra in $\mathcal{M}$.
\end{theorem}

\begin{proof}
It is analogous to the proof of Theorem \ref{teo: cotensor bial}.
\end{proof}

We are now able to state the universal property of the tensor bialgebra.

\begin{theorem}
\label{teo: univ property of tensor bialgebra} Let $H$ be a braided
bialgebra in a cocomplete abelian braided monoidal category $\left( \mathcal{%
M},c\right) .$ Assume that the tensor product commutes with direct sums.
\newline
Let $(M,\mu _{M}^{r},\mu _{M}^{l},\rho _{M}^{r},\rho _{M}^{l})$ be in $%
_{H}^{H}\mathcal{M}_{H}^{H}$. Let $T=T_{H}(M)$ be the tensor algebra.\newline
Let $E$ be a braided bialgebra in $\mathcal{M}$.\newline
Let $f_{H}:H\rightarrow E$ be a bialgebra homomorphism and let $%
f_{M}:M\rightarrow E$ be a morphism of $H$-bimodules, where $E$ is a
bimodule via $f_{H}.$ Assume that%
\begin{equation*}
\Delta _{E}f_{M}=\left( f_{H}\otimes f_{M}\right) \rho _{M}^{l}+\left(
f_{M}\otimes f_{H}\right) \rho _{M}^{r}
\end{equation*}%
(i.e. $f_{M}$ is a coderivation of $E$ with domain the $E$-bicomodule $M$,
where $M$ is regarded as a bicomodule via $f_{H}$). Then there is a unique
algebra homomorphism $f:T_{H}(M)\rightarrow E$ such that $fi_{0}=f_{H}$ and $%
fi_{1}=f_{M}$, where $i_{n}:M^{\otimes _{H}n}\rightarrow T_{H}(M)$ denotes
the canonical injection.
\begin{equation*}
\diagUnivTensBialg%
\end{equation*}%
Moreover $f$ is a bialgebra homomorphism.
\end{theorem}

\begin{proof}
By Theorem \ref{teo: univ property of tensor algebra}, there is a unique
algebra homomorphism $f:T\rightarrow E$ such that $fi_{0}=f_{H}$ and $%
fi_{1}=f_{M}$. The rest of the proof is similar to that of Theorem\ref{teo:
univ property of cotensor bialgebra}.
\end{proof}

\section{Braided Bialgebras of Type One}

Even if not all the assumptions are always needed, in this section $\left(
\left( \mathcal{M},\otimes ,\mathbf{1}\right) ,c\right) $ will denote a
cocomplete and complete abelian coabelian braided monoidal category
satisfying $AB5$ and such that the tensor product commutes with direct sums.

\begin{lemma}
\label{lem: forType1}Let $(B,m_{B},u_{B},\Delta _{B},\varepsilon _{B})$ be a
graded braided bialgebra in $\left( \mathcal{M},c\right) $. Then

1) $(B_{0},m_{0,0},u_{0},\Delta _{0,0},\varepsilon _{0})$ is a braided
bialgebra. Moreover both $p_{0}^{B}$ and $i_{0}^{B}$ are bialgebra
homomorphisms.

2) $(B_{1},\mu _{B_{1}}^{r}=m_{1,0},\mu _{B_{1}}^{l}=m_{0,1},\rho
_{B_{1}}^{r}=\Delta _{1,0},\rho _{B_{1}}^{l}=\Delta _{0,1})$ is an object in
$_{B_{0}}^{B_{0}}\mathcal{M}_{B_{0}}^{B_{0}}$.

3) $i_{1}^{B}$ is a morphism of $B_{0}$-bimodules, where $B$ is a $B_{0}$%
-bimodule via $i_{0}^{B}$, and we have
\begin{equation*}
\Delta _{B}i_{1}^{B}=\left( i_{0}^{B}\otimes i_{1}^{B}\right) \Delta
_{0,1}+\left( i_{1}^{B}\otimes i_{0}^{B}\right) \Delta _{1,0}
\end{equation*}%
(i.e. $i_{1}^{B}$ is a coderivation of $B$ with domain the $B$-bicomodule $%
B_{1}$, where $B_{1}$ is regarded as a bicomodule via $i_{0}^{B}$).

4) $p_{1}^{B}$ is a morphism of $B_{0}$-bicomodules, where $B$ is a $B_{0}$%
-bicomodule via $p_{0}^{B}$, and we have
\begin{equation*}
p_{1}^{B}m_{B}=m_{0,1}\left( p_{0}^{B}\otimes p_{1}^{B}\right)
+m_{1,0}\left( p_{1}^{B}\otimes p_{0}^{B}\right)
\end{equation*}%
(i.e. $p_{1}^{B}$ is a derivation of $B$ with values in the $B$-bimodule $%
B_{1}$, where $B_{1}$ is regarded as a bimodule via $p_{0}^{B}$).
\end{lemma}

\begin{proof}
We prove the statements by means of (\ref{form: def Braided}).

1) By Proposition \ref{lem: graded Deltaij}, $\left( B_{0},\Delta
_{0}=\Delta _{0,0},\varepsilon _{0}=\varepsilon _{B}i_{0}^{B}\right) $ is a
coalgebra in $\mathcal{M}$ and $i_{0}^{B}$ is a coalgebra homomorphism. By
Proposition \ref{lem: graded mij} $\left(
B_{0},m_{0}=m_{0,0},u_{0}=p_{0}^{B}u_{B}\right) $ is an algebra in $\mathcal{%
M}$ and $p_{0}^{B}$ is an algebra homomorphism.

By (\ref{form: grAlg}), (\ref{form: coro grCoalg1}) and by naturality of the
braiding we obtain
\begin{eqnarray*}
&&\left( m_{0,0}\otimes m_{0,0}\right) \left( B_{0}\otimes
c_{B_{0},B_{0}}\otimes B_{0}\right) \left( \Delta _{0,0}\otimes \Delta
_{0,0}\right) \\
&=&\left( p_{0}^{B}\otimes p_{0}^{B}\right) \left( m_{B}\otimes m_{B}\right)
\left( B\otimes c_{B,B}\otimes B\right) \left( \Delta _{B}\otimes \Delta
_{B}\right) \left( i_{0}^{B}\otimes i_{0}^{B}\right) \\
&\overset{\text{(\ref{form: def Braided})}}{=}&\left( p_{0}^{B}\otimes
p_{0}^{B}\right) \Delta _{B}m_{B}\left( i_{0}^{B}\otimes i_{0}^{B}\right)
=\Delta _{0,0}m_{0,0}.
\end{eqnarray*}%
Moreover
\begin{equation*}
m_{\mathbf{1}}\left( \varepsilon _{0}\otimes \varepsilon _{0}\right) =m_{%
\mathbf{1}}\left( \varepsilon _{B}\otimes \varepsilon _{B}\right) \left(
i_{0}^{B}\otimes i_{0}^{B}\right) \overset{\text{(\ref{form: def Braided})}}{%
=}\varepsilon _{B}m_{B}\left( i_{0}^{B}\otimes i_{0}^{B}\right) \overset{%
\text{(\ref{form: coro grAlg1})}}{=}\varepsilon
_{B}i_{0}^{B}m_{0,0}=\varepsilon _{0}m_{0,0}
\end{equation*}%
so that $(B_{0},m_{0,0},u_{0},\Delta _{0,0},\varepsilon _{0})$ is a braided
bialgebra.\newline
Let us check that $i_{0}^{B}$ is an algebra homomorphism. We have%
\begin{equation*}
m_{B}\left( i_{0}^{B}\otimes i_{0}^{B}\right) \overset{\text{(\ref{form:
coro grAlg1})}}{=}i_{0}^{B}m_{0,0}.
\end{equation*}%
Moreover, by Proposition \ref{lem: graded mij}, $%
u=i_{0}^{B}p_{0}u=i_{0}^{B}u_{0} $. Thus $i_{0}^{B}$ is an algebra
homomorphism. \newline
Similarly, by means of (\ref{form: coro grCoalg1}) and of Proposition \ref%
{lem: graded Deltaij}, we get that $p_{0}^{B}$ is also a coalgebra
homomorphism.

2) By (\ref{form: grAlg}), (\ref{form: coro grCoalg1}) and by naturality of
the braiding we obtain%
\begin{eqnarray}
&&\left( p_{0}^{B}\otimes p_{1}^{B}\right) \left( m_{B}\otimes m_{B}\right)
\left( B\otimes c_{B,B}\otimes B\right) \left( \Delta _{B}\otimes \Delta
_{B}\right)  \label{form: Hopfmod 1} \\
&=&\left[
\begin{array}{c}
\left( m_{0,0}\otimes m_{1,0}\right) \left( B_{0}\otimes
c_{B_{1},B_{0}}\otimes B_{0}\right) \left[ \Delta _{0,1}p_{1}^{B}\otimes
\Delta _{0,0}p_{0}^{B}\right] + \\
+\left( m_{0,0}\otimes m_{0,1}\right) \left( B_{0}\otimes
c_{B_{0},B_{0}}\otimes B_{1}\right) \left[ \Delta _{0,0}p_{0}^{B}\otimes
\Delta _{0,1}p_{1}^{B}\right]%
\end{array}%
\right]  \notag
\end{eqnarray}%
Hence, by (\ref{form: coro grCoalg1}), (\ref{form: coro grAlg1}) and by
naturality of the braiding we obtain
\begin{eqnarray*}
\Delta _{0,1}m_{0,1} &=&\left( p_{0}^{B}\otimes p_{1}^{B}\right) \Delta
_{B}m_{B}\left( i_{0}^{B}\otimes i_{1}^{B}\right) \\
&\overset{\text{(\ref{form: def Braided})}}{=}&\left( p_{0}^{B}\otimes
p_{1}^{B}\right) \left( m_{B}\otimes m_{B}\right) \left( B\otimes
c_{B,B}\otimes B\right) \left( \Delta _{B}\otimes \Delta _{B}\right) \left(
i_{0}^{B}\otimes i_{1}^{B}\right) \\
&\overset{\text{(\ref{form: Hopfmod 1})}}{=}&\left( m_{0,0}\otimes
m_{0,1}\right) \left( B_{0}\otimes c_{B_{0},B_{0}}\otimes B_{1}\right) \left[
\Delta _{0,0}\otimes \Delta _{0,1}\right]
\end{eqnarray*}%
Similarly we get (\ref{form: HopfBimod 2}), (\ref{form: HopfBimod 3}) and (%
\ref{form: HopfBimod 4}).

3)-4) follow in view of Proposition \ref{pro: grCoalg} and \ref{pro: grAlg}.
\end{proof}

\begin{theorem}
\label{coro: GRuniv property of cotensor bialgebra}Let $H$ be a braided
bialgebra in a cocomplete and complete abelian and coabelian braided
monoidal category $\left( \mathcal{M},c\right) $ satisfying $AB5.$ Assume
that the tensor product commutes with direct sums. \newline
Let $(M,\mu _{M}^{r},\mu _{M}^{l},\rho _{M}^{r},\rho _{M}^{l})$ be in $%
_{H}^{H}\mathcal{M}_{H}^{H}$. Let $T^{c}=T_{H}^{c}(M)$ be the cotensor
coalgebra.\newline
Let $B$ be a graded bialgebra, let $g_{H}:B_{0}\rightarrow H$ be a bialgebra
homomorphism and let $g_{M}:B_{1}\rightarrow M$ be a morphism of $H$%
-bicomodules, where $B_{1}$ is an$\ H$-bicomodule via $g_{H},$ and a
morphism of $B_{0}$-bimodules, where $M\ $is a $B_{0}$-bimodule via $g_{H}$.%
\newline
Then there is a unique coalgebra homomorphism $f:B\rightarrow T_{H}^{c}(M)$
such that
\begin{equation*}
p_{0}^{T^{c}}\circ f=g_{H}\circ p_{0}^{B}\text{\qquad and\qquad }%
p_{1}^{T^{c}}\circ f=g_{M}\circ p_{1}^{B}.
\end{equation*}%
Moreover $f$ is a graded bialgebra homomorphism.
\end{theorem}

\begin{proof}
By Theorem \ref{coro: GRuniv property of cotensor coalgebra}, there is a
unique coalgebra homomorphism $f:B\rightarrow T_{H}^{c}(M)$ such that $%
p_{0}^{C}\circ f=g_{H}\circ p_{0}^{B}$ and $p_{1}^{T^{c}}\circ f=g_{M}\circ
p_{1}^{B}.$ Moreover $f$ is a morphism of graded coalgebras. On the other
hand, set%
\begin{equation*}
f_{H}:=g_{H}\circ p_{0}^{B}\text{\qquad and\qquad }f_{M}:=g_{M}\circ
p_{1}^{B}.
\end{equation*}
By Lemma \ref{lem: forType1}, we have
\begin{equation*}
p_{1}^{B}m_{B}=m_{0,1}\left( p_{0}^{B}\otimes p_{1}^{B}\right)
+m_{1,0}\left( p_{1}^{B}\otimes p_{0}^{B}\right) .
\end{equation*}%
Since $g_{M}:B_{1}\rightarrow M$ is a morphism of $B_{0}$-bimodules, where $%
M\ $is a $B_{0}$-bimodule via $g_{H}$, we get%
\begin{eqnarray*}
f_{M}\circ m_{B} &=&g_{M}\circ p_{1}^{B}\circ m_{B} \\
&=&\mu _{M}^{l}\circ \left( g_{H}\otimes g_{M}\right) \circ \left(
p_{0}^{B}\otimes p_{1}^{B}\right) +\mu _{M}^{r}\circ \left( g_{M}\otimes
g_{H}\right) \circ \left( p_{1}^{B}\otimes p_{0}^{B}\right) \\
&=&\mu _{M}^{l}\circ \left( f_{H}\otimes f_{M}\right) +\mu _{M}^{r}\circ
\left( f_{M}\otimes f_{H}\right) .
\end{eqnarray*}%
Hence, by Proposition \ref{pro: lim for graded coalg} and Theorem \ref{teo:
univ property of cotensor bialgebra}, there is a unique coalgebra
homomorphism $f^{\prime }:B\rightarrow T_{H}^{c}(M)$ such that $%
p_{0}f^{\prime }=f_{H}$ and $p_{1}f^{\prime }=f_{M}$. Moreover $f^{\prime }$
is a bialgebra homomorphism. By uniqueness $f=f^{\prime }.$
\end{proof}

\begin{theorem}
\label{coro: GRuniv property of tensor bialgebra}Let $H$ be a braided
bialgebra in a cocomplete and complete abelian and coabelian braided
monoidal category $\left( \mathcal{M},c\right) $ satisfying AB5. Assume that
the tensor product commutes with direct sums. \newline
Let $(M,\mu _{M}^{r},\mu _{M}^{l},\rho _{M}^{r},\rho _{M}^{l})$ be in $%
_{H}^{H}\mathcal{M}_{H}^{H}$. Let $T=T_{H}(M)$ be the tensor algebra.\newline
Let $B$ be a graded bialgebra, let $g_{H}:H\rightarrow B_{0}$ be a bialgebra
homomorphism and let $g_{M}:M\rightarrow B_{1}$ be a morphism of $H$%
-bimodules, where $B_{1}$ is an $H$-bimodule via $g_{H}$, and a morphism of $%
B_{0}$-bicomodules, where $M$ is a $B_{0}$-bicomodule via $g_{H}.$ \newline
Then there is a unique algebra homomorphism $f:T_{H}(M)\rightarrow B$ such
that
\begin{equation*}
f\circ i_{0}^{T}=i_{0}^{B}\circ g_{H}\text{\qquad and\qquad }f\circ
i_{1}^{T}=i_{1}^{B}\circ g_{M}.
\end{equation*}%
Moreover $f$ is a graded bialgebra homomorphism.
\end{theorem}

\begin{proof}
It is analogous to the proof of Theorem \ref{coro: GRuniv property of
cotensor bialgebra}.
\end{proof}

\begin{claim}
\label{cl: phi_f psi_f}Let $f:X\rightarrow X^{\prime }$ be a graded
homomorphism in $\mathcal{M}$. Hence $f=\oplus _{n\in
\mathbb{N}
}f_{n}$ for suitable morphisms $f_{n}:X_{n}\rightarrow X_{n}^{\prime }$, for
every $n\in
\mathbb{N}
$. Write%
\begin{equation*}
f_{n}=\psi _{n}\circ \varphi _{n}
\end{equation*}%
where $\varphi _{n}:X_{n}\rightarrow \mathrm{Im}\left( f_{n}\right) $ is an
epimorphism and $\psi _{n}:\mathrm{Im}\left( f_{n}\right) \rightarrow
X_{n}^{\prime }$ is a monomorphism. Since $\mathcal{M}$ satisfies $AB5$, it
is in particular an $AB4$ category. Hence coproducts are both left and right
exact. Thus
\begin{equation*}
\varphi _{f}:=\oplus _{n\in
\mathbb{N}
}\varphi _{n}:X\rightarrow \oplus _{n\in
\mathbb{N}
}\mathrm{Im}\left( f_{n}\right) \qquad \text{and}\qquad \psi _{f}:=\oplus
_{n\in
\mathbb{N}
}\psi _{n}:\oplus _{n\in
\mathbb{N}
}\mathrm{Im}\left( f_{n}\right) \rightarrow X^{\prime }
\end{equation*}%
are an epimorphism and a monomorphism respectively. Moreover%
\begin{equation*}
f=\oplus _{n\in
\mathbb{N}
}f_{n}=\left( \oplus _{n\in
\mathbb{N}
}\psi _{n}\right) \circ \left( \oplus _{n\in
\mathbb{N}
}\varphi _{n}\right) =\psi _{f}\circ \varphi _{f}.
\end{equation*}%
In particular we obtain $\mathrm{Im}\left( f\right) =\oplus _{n\in
\mathbb{N}
}\mathrm{Im}\left( f_{n}\right) $ so that $\mathrm{\mathrm{Im}}\left(
f\right) _{n}=\mathrm{\mathrm{Im}}\left( f_{n}\right) $ for every $n\in
\mathbb{N}
$. We would like to outline that the last statements do not in general hold
without assuming $AB5$. Clearly if $f_{n}=\mathrm{Id}_{X_{n}}$ we can choose
$\psi _{n}=\varphi _{n}=\mathrm{Id}_{X_{n}}$.
\end{claim}

\begin{proposition}
\label{pro: Im of GrBrHomo}Let $f:A\rightarrow A^{\prime }$ be a graded
braided bialgebra homomorphism in $\left( \mathcal{M},c\right) $. Then $%
\mathrm{\mathrm{Im}}\left( f\right) $ can be endowed with a unique graded
braided bialgebra structure such that $\varphi _{f}:A\rightarrow \mathrm{Im}%
\left( f\right) $ and $\psi _{f}:\mathrm{Im}\left( f\right) \rightarrow
A^{\prime }$ are graded braided bialgebra homomorphisms.\newline
Furthermore

1) $\varphi _{f}$ is epimorphism and $\psi _{f}$ is a monomorphism in $%
\mathcal{M}$.

2) $\mathrm{\mathrm{Im}}\left( f\right) _{n}=\mathrm{\mathrm{Im}}\left(
f_{n}\right) $ for every $n\in
\mathbb{N}
$.
\end{proposition}

\begin{proof}
Set $B=\mathrm{Im}\left( f\right) $, $B_{n}=\mathrm{Im}\left( f_{n}\right) ,$
$\varphi =\varphi _{f}$ and $\psi =\psi _{f}$.

Since $\left( \mathrm{Im}\left( f\right) ,\varphi \right) =\mathrm{\mathrm{%
\mathrm{\mathrm{Coker}}}}\left( \mathrm{\ker }\left( f\right) \right) ,$ it
is clear that $\mathrm{Im}\left( f\right) $ can be endowed with a unique
algebra structure such that $\varphi $ is an algebra homomorphism.

Since $\left( \mathrm{Im}\left( f\right) ,\psi \right) =\mathrm{\ker }\left(
\mathrm{\mathrm{\mathrm{\mathrm{Coker}}}}\left( f\right) \right) ,$ it is
clear that $\mathrm{Im}\left( f\right) $ can be endowed with a unique
coalgebra structure such that $\psi $ is a coalgebra homomorphism. \newline
Let us check that $\varphi $ is also a coalgebra homomorphism. We have%
\begin{eqnarray*}
\left( \psi \otimes \psi \right) \circ \Delta _{B}\circ \varphi &=&\Delta
_{A^\prime }\circ \psi \circ \varphi =\Delta _{A^\prime }\circ f=\left(
f\otimes f\right) \circ \Delta _{A}=\left( \psi \otimes \psi \right) \circ
\left( \varphi \otimes \varphi \right) \circ \Delta _{A}, \\
\varepsilon _{B}\circ \varphi &=&\varepsilon _{A^\prime }\circ \psi \circ
\varphi =\varepsilon _{A^\prime }\circ f=\varepsilon _{A}.
\end{eqnarray*}%
Since $\psi $ is a monomorphism and $\mathcal{M}$ is a coabelian monoidal
category, we get that $\psi \otimes \psi $ is a monomorphism so that
\begin{equation*}
\Delta _{B}\circ \varphi =\left( \varphi \otimes \varphi \right) \circ
\Delta _{A}\qquad \text{and}\qquad \varepsilon _{B}\circ \varphi
=\varepsilon _{A}
\end{equation*}%
and hence $\varphi $ is also a coalgebra homomorphism.\newline
Similarly it can be proven that $\psi $ is also an algebra homomorphism.
\newline
Let us check that $B$ is a braided bialgebra. Since $\varphi $ is a
bialgebra homomorphism, $A$ is a braided bialgebra and by naturality of the
braiding one computes%
\begin{gather*}
\left( m_{B}\otimes m_{B}\right) \circ \left( B\otimes c_{B,B}\otimes
B\right) \circ \left( \Delta _{B}\otimes \Delta _{B}\right) \circ \left(
\varphi \otimes \varphi \right) =\Delta _{B}\circ m_{B}\circ \left( \varphi
\otimes \varphi \right) \\
\varepsilon _{B}\circ m_{B}\circ \left( \varphi \otimes \varphi \right) =m_{%
\mathbf{1}}\circ \left( \varepsilon _{B}\otimes \varepsilon _{B}\right)
\circ \left( \varphi \otimes \varphi \right) .
\end{gather*}%
Since $\varphi \otimes \varphi $ is an epimorphism we get (\ref{form: def
Braided}) for $B.$\newline
By Proposition \ref{coro: grCoalg} in order to prove that $B=\oplus _{n\in
\mathbb{N}
}B_{n}$ is a graded coalgebra we have to check that there are morphisms $%
\Delta _{a,b}^{B}:B_{a+b}\rightarrow B_{a}\otimes B_{b},$ such that $%
(p_{a}^{B}\otimes p_{b}^{B})\circ \Delta _{B}=\Delta _{a,b}^{B}\circ
p_{a+b}^{B},\text{ for every }a,b\in
\mathbb{N}
. $ Set%
\begin{equation*}
\Delta _{a,b}^{B}:=\left( B_{a+b}\overset{i_{a+b}^{B}}{\rightarrow }B\overset%
{\Delta _{B}}{\rightarrow }B\otimes B\overset{p_{a}^{B}\otimes p_{b}^{B}}{%
\rightarrow }B_{a}\otimes B_{b}\right) .
\end{equation*}%
Since $\varphi $ is a graded coalgebra homomorphism, by (\ref{form: coro
grCoalg1}), one easily checks that%
\begin{equation*}
\Delta _{a,b}^{B}\circ p_{a+b}^{B}\circ \varphi =(p_{a}^{B}\otimes
p_{b}^{B})\circ \Delta _{B}\circ \varphi
\end{equation*}
Since $\varphi $ is an epimorphism, we get $\Delta _{a,b}^{B}\circ
p_{a+b}^{B}=(p_{a}^{B}\otimes p_{b}^{B})\circ \Delta _{B} $ that is $B$ is a
graded coalgebra. \newline
Similarly one proves that $B$ is also a graded algebra and hence a graded
braided bialgebra.

Clearly both $\varphi $ and $\psi $ are graded braided bialgebra
homomorphism.

$1)$ and $2)$ follow by (\ref{cl: phi_f psi_f}).
\end{proof}

\begin{theorem}
\label{teo: F of type1} Let $H$ be a braided bialgebra in a cocomplete and
complete abelian coabelian braided monoidal category $\left( \mathcal{M}%
,c\right) $ satisfying $AB5$. Assume that the tensor product commutes with
direct sums. \newline
Let $(M,\mu _{M}^{r},\mu _{M}^{l},\rho _{M}^{r},\rho _{M}^{l})$ be in $%
_{H}^{H}\mathcal{M}_{H}^{H}$. Let $T=T_{H}(M)$ and $T^{c}=T_{H}^{c}(M).$
\newline
Then there is a unique algebra homomorphism $F=F_{H,M}:T\rightarrow T^{c}$
such that $Fi_{0}^{T}=i_{0}^{T^{c}}$ and $Fi_{1}^{T}=i_{1}^{T^{c}}$:
\begin{equation*}
\diagFofTypeOne%
\end{equation*}%
Moreover $F$ is a graded bialgebra homomorphism such that $F_{0}=\mathrm{Id}%
_{H}$ and $F_{1}=\mathrm{Id}_{M}$. \newline
Write $F=\psi _{F}\circ \varphi _{F}$ as in Proposition \ref{pro: Im of
GrBrHomo}, where $\varphi _{F}:T\rightarrow \mathrm{Im}\left( F\right) $ and
$\psi _{F}:\mathrm{Im}\left( F\right) \rightarrow T^{c}$ are graded braided
bialgebra homomorphisms, $\varphi _{F}$ is epimorphism and $\psi _{F}$ is a
monomorphism in $\mathcal{M}$. We have:

1) $\mathrm{\mathrm{Im}}\left( F\right) $ is a graded braided bialgebra such
that $\mathrm{\mathrm{Im}}\left( F\right) _{0}=H$ and $\mathrm{\mathrm{Im}}%
\left( F\right) _{1}=M$.

2) $\varphi _{F}$ is the unique algebra homomorphism such that $\varphi
_{F}\circ i_{0}^{T}=i_{0}^{\mathrm{Im}\left( F\right) }$ and $\varphi
_{F}\circ i_{1}^{T}=i_{1}^{\mathrm{Im}\left( F\right) }$ as in Theorem \ref%
{teo: phi epi}.

3) $\psi _{F}$ is the unique coalgebra homomorphism such that $%
p_{0}^{T^{c}}\circ \psi _{F}=p_{0}^{\mathrm{Im}\left( F\right) }$ and $%
p_{1}^{T^{c}}\circ \psi _{F}=p_{1}^{\mathrm{Im}\left( F\right) }$ as in
Theorem \ref{teo: psi mono}.
\end{theorem}

\begin{proof}
Set $B=\mathrm{Im}\left( F\right) .$ In view of Theorem \ref{teo: cotensor
bial} $T^{c}=T_{H}^{c}(M)$ is a graded braided bialgebra.

By Corollary \ref{coro: GRuniv property of tensor bialgebra} there is a
unique algebra homomorphism $F:T\rightarrow T^{c}$ such that $F\circ
i_{0}^{T}=i_{0}^{T^{c}}$ and $F\circ i_{1}^{T}=i_{1}^{T^{c}}$. Moreover $F$
is a graded bialgebra homomorphism.\newline
By Proposition \ref{pro: Im of GrBrHomo}, we have that $\mathrm{\mathrm{Im}}%
\left( F\right) $ is a graded braided bialgebra such that $\mathrm{\mathrm{Im%
}}\left( F\right) _{n}=\mathrm{\mathrm{Im}}\left( F_{n}\right) $ for every $%
n\in
\mathbb{N}
$. Moreover we can write $F=\psi \circ \varphi $ where $\varphi =\varphi
_{F}:T\rightarrow \mathrm{Im}\left( F\right) $ and $\psi =\psi _{F}:\mathrm{%
Im}\left( F\right) \rightarrow T^{c}$ are graded braided bialgebra
homomorphisms, $\varphi $ is an epimorphism and $\psi $ is a monomorphism in
$\mathcal{M}$. Since $F$ is graded we have%
\begin{equation*}
F_{0}=p_{0}^{T^{c}}\circ F\circ i_{0}^{T}=p_{0}^{T^{c}}\circ i_{0}^{T^{C}}=%
\mathrm{Id}_{H}\qquad \text{and}\qquad F_{1}=p_{1}^{T^{c}}\circ F\circ
i_{1}^{T}=p_{1}^{T^{c}}\circ i_{1}^{T^{C}}=\mathrm{Id}_{M}.
\end{equation*}%
Then $\mathrm{\mathrm{Im}}\left( F\right) _{0}=\mathrm{Im}\left(
F_{0}\right) =H$ and $\mathrm{\mathrm{Im}}\left( F\right) _{1}=\mathrm{Im}%
\left( F_{1}\right) =M.$ Now, since $F_{0}=\mathrm{Id}_{H}$ and $F_{1}=%
\mathrm{Id}_{M}$ we can choose $\varphi _{0}=\psi _{0}=\mathrm{Id}_{H}$ and $%
\varphi _{1}=\psi _{1}=\mathrm{Id}_{M}$ so that $\varphi \circ
i_{n}^{T}=i_{n}^{B}$ and $p_{n}^{B}=p_{n}^{T^{c}}\circ \psi $ for $n=0,1$.
\end{proof}

\begin{definition}
\label{def: type1} Take the notations and assumptions of Theorem \ref{teo: F
of type1}. Following \cite[page 1533]{Ni}, let
\begin{equation*}
(H[M],i_{H[M]})=\left( \mathrm{Im}(F),\psi _{F}\right) .
\end{equation*}%
Then $H[M]$ is a graded braided bialgebra such that $H[M]_{0}=H$ and $%
H[M]_{1}=M$. \newline
This will be called \emph{the braided bialgebra of type one associated to $H$
and $M$}.
\end{definition}

\begin{theorem}
\label{teo: type One}Let $(B,m_{B},u_{B},\Delta _{B},\varepsilon _{B})$ be a
braided graded bialgebra in a cocomplete and complete abelian coabelian
braided monoidal category $\left( \mathcal{M},c\right) $ satisfying $AB5$.
Assume that the tensor product commutes with direct sums.

\begin{itemize}
\item[1)] There is a unique algebra homomorphism
\begin{equation*}
\varphi :T_{B_{0}}(B_{1})\rightarrow B
\end{equation*}
such that $\varphi \circ i_{0}^{T}=i_{0}^{B}$ and $\varphi \circ
i_{1}^{T}=i_{1}^{B}$. Moreover $\varphi $ is a graded bialgebra homomorphism.

\item[2)] There is a unique coalgebra homomorphism
\begin{equation*}
\psi :B\rightarrow T_{B_{0}}^{c}(B_{1})
\end{equation*}%
such that $p_{0}^{T^{c}}\circ \psi =p_{0}^{B}$ and $p_{1}^{T^{c}}\circ \psi
=p_{1}^{B}$. Moreover $\psi $ is a graded bialgebra homomorphism.

\item[3)] $\psi \circ \varphi =F$ where $F:T_{B_{0}}(B_{1})\rightarrow
T_{B_{0}}^{c}(B_{1})$ is the morphism defined in Theorem \ref{teo: F of
type1}. In particular $F$ is a graded bialgebra homomorphism.
\end{itemize}

Furthermore the following assertions are equivalent

$\left( a\right) $ $B$ is strongly $%
\mathbb{N}
$-graded both as an algebra and a coalgebra.

$\left( b\right) $ $\varphi $ is an epimorphism and $\psi $ is a
monomorphism.

$\left( c\right) $ $\varphi =\varphi _{F}$ and $\psi =\psi _{F}$.

$\left( d\right) $ $B=B_{0}\left[ B_{1}\right] $ is the braided bialgebra of
type one associated to $B_{0}$ and $B_{1}.$
\end{theorem}

\begin{proof}
1)-2) By Lemma \ref{lem: forType1}, we can apply Theorems \ref{coro: GRuniv
property of tensor bialgebra} and \ref{coro: GRuniv property of cotensor
bialgebra}.

3) By Theorem \ref{teo: F of type1}, there is a unique algebra homomorphism $%
F:T_{B_{0}}(B_{1})\rightarrow T_{B_{0}}^{c}(B_{1})$ such that $F\circ
i_{0}^{T}=i_{0}^{T^{c}}$ and $F\circ i_{1}^{T}=i_{1}^{T^{c}}$. Moreover $F$
is a graded bialgebra homomorphism such that $F_{0}=\mathrm{Id}_{H}$ and $%
F_{1}=\mathrm{Id}_{M}$. Since%
\begin{equation*}
\psi \circ \varphi \circ i_{0}^{T}=\psi \circ i_{0}^{B}\circ \varphi
_{0}=i_{0}^{T^{c}}\circ \psi _{0}\circ \varphi _{0}=i_{0}^{T^{c}},\qquad
\psi \circ \varphi \circ i_{1}^{T}=\psi \circ i_{1}^{B}\circ \varphi
_{1}=i_{1}^{T^{c}}\circ \psi _{1}\circ \varphi _{1}=i_{1}^{T^{c}},
\end{equation*}%
we infer that $F=\psi \circ \varphi $. By Proposition \ref{pro: Im of
GrBrHomo}, $\mathrm{\mathrm{Im}}\left( F\right) $ can be endowed with a
unique graded braided bialgebra structure such that $\varphi
_{F}:T\rightarrow \mathrm{Im}\left( F\right) $ and $\psi _{F}:\mathrm{Im}%
\left( F\right) \rightarrow T^{c}$ are graded braided bialgebra
homomorphisms. Furthermore $\mathrm{\mathrm{Im}}\left( F\right) _{n}=\mathrm{%
\mathrm{Im}}\left( F_{n}\right) $ for every $n\in
\mathbb{N}
$, we can write $F=\psi _{F}\circ \varphi _{F}$, $\varphi _{F}$ is an
epimorphism and $\psi _{F}$ is a monomorphism in $\mathcal{M}$. \newline
$\left( a\right) \Leftrightarrow \left( b\right) $ It follows by Theorems %
\ref{teo: phi epi} and \ref{teo: psi mono}.\newline
$\left( b\right) \Rightarrow \left( c\right) $ By hypothesis, since the
category $\mathcal{M}$ is abelian, $B$ identifies as an object with $\mathrm{%
Im}\left( F\right) =B_{0}\left[ B_{1}\right] $. By uniqueness in 2) and 3)
of Theorem \ref{teo: F of type1} and in view of 1) and 2) above, we get $%
\varphi =\varphi _{F}$ and $\psi =\psi _{F}$.\newline
$\left( c\right) \Rightarrow \left( d\right) $ Since the graded braided
bialgebra structure of $B$ and of $\mathrm{Im}\left( F\right) $ make $%
\varphi _{F}:T\rightarrow B$ and $\psi _{F}:B\rightarrow T^{c}$ graded
braided bialgebra homomorphisms, by Proposition \ref{pro: Im of GrBrHomo},
the two structures coincide. \newline
$\left( d\right) \Rightarrow \left( c\right) $ Since $B=\mathrm{Im}\left(
F\right) ,$ it follows by 2) and 3) in Theorem \ref{teo: F of type1}.\newline
$\left( c\right) \Rightarrow \left( b\right) $ It is trivial.
\end{proof}

\begin{remark}
Recall that $\varphi $ is an epimorphism if and only if it fulfills one of
the equivalent conditions in Theorem \ref{teo: phi epi} and $\psi $ is a
monomorphism if and only if it fulfills one of the equivalent conditions in
Theorem \ref{teo: psi mono}. In particular we have the following result.
\end{remark}

\begin{theorem}
\label{teo: Magnum}Let $(B,m_{B},u_{B},\Delta _{B},\varepsilon _{B})$ be a
braided graded bialgebra in a cocomplete and complete abelian coabelian
braided monoidal category $\left( \mathcal{M},c\right) $ satisfying $AB5$.
Assume that the tensor product commutes with direct sums. \newline
Then $B$ is the braided bialgebra of type one $B_{0}\left[ B_{1}\right] $
associated to $B_{0}$ and $B_{1}$ if and only if
\begin{equation*}
(B\left[ 2\right] ,\nu _{2}^{B})=B\left[ 1\right] ^{2}\qquad \text{and}%
\qquad (B\left( 2\right) ,\sigma _{2}^{B})=B_{0}^{\wedge _{B}^{2}}\text{.}
\end{equation*}
\end{theorem}

\section{Radford-Majid bosonization}

In this section $\left( \left( \mathcal{M},\otimes ,\mathbf{1}\right)
,c\right) $ will denote a cocomplete and complete abelian coabelian braided
monoidal category satisfying $AB5$ and such that the tensor product commutes
with direct sums.

\begin{claim}
\label{claim: YDbraiding}Let $(H,m_{H},u_{H},\Delta _{H},\varepsilon
_{H},S_{H})$ be a braided graded Hopf algebra with invertible antipode $%
S_{H} $. By \cite[Theorem 3.4.3]{Bes}, $\left( \left( _{H}^{H}\mathcal{YD}%
\left( \mathcal{M}\right) ,\otimes ,\mathbf{1}\right) ,\Psi \right) $ is a
braided monoidal category, where, for every $\left( U,\mu _{U}^{l},\rho
_{U}^{l}\right) ,\left( V,\mu _{V}^{l},\rho _{V}^{l}\right) \in {_{H}^{H}%
\mathcal{YD}\left( \mathcal{M}\right) }$, the braiding is defined by%
\begin{equation*}
\Psi _{U,V}:=\left( \mu _{V}^{l}\otimes U\right) \circ \left( H\otimes
c_{U,V}\right) \circ \left( \rho _{U}^{l}\otimes V\right)
\end{equation*}%
with inverse given by%
\begin{equation*}
\Psi _{U,V}^{-1}:=\left( U\otimes \mu _{V}^{l}\right) \circ \left( U\otimes
c_{H,V}^{-1}\right) \circ \left( c_{U,V}^{-1}\otimes S_{H}^{-1}\right) \circ
\left( V\otimes c_{U,H}^{-1}\right) \circ \left( V\otimes \rho
_{U}^{l}\right) .
\end{equation*}
\end{claim}

\begin{claim}
\label{cl: BesYD}In view of \cite[Lemma 3.9.4]{Bes}, recall that given an
object $\left( V,\mu _{V}^{l},\rho _{V}^{l}\right) $ in the monoidal
category $\left( _{H}^{H}\mathcal{YD}\left( \mathcal{M}\right) ,\otimes ,%
\mathbf{1}\right) ,$ then
\begin{equation*}
\left( V\otimes H,\mu _{V\otimes H}^{l},\mu _{V\otimes H}^{r},\rho
_{V\otimes H}^{l},\rho _{V\otimes H}^{r}\right)
\end{equation*}%
is an object in $_{H}^{H}\mathcal{M}_{H}^{H}$ where
\begin{eqnarray*}
\mu _{V\otimes H}^{l} &=&\left( \mu _{V}^{l}\otimes m_{H}\right) \circ
\left( H\otimes c_{H,V}\otimes H\right) \circ \left( \Delta _{H}\otimes
V\otimes H\right) ,\qquad \mu _{V\otimes H}^{r}=V\otimes m_{H}, \\
\rho _{V\otimes H}^{l} &=&\left( m_{H}\otimes V\otimes H\right) \circ \left(
H\otimes c_{V,H}\otimes H\right) \circ \left( \rho _{V}^{l}\otimes \Delta
_{H}\right) ,\qquad \rho _{V\otimes H}^{r}=V\otimes \Delta _{H}.
\end{eqnarray*}

Given a braided bialgebra $\left( \left( Q,\mu _{Q}^{l},\rho _{Q}^{l}\right)
,m_{Q},u_{Q},\Delta _{Q},\varepsilon _{Q}\right) $ in the braided monoidal
category $\left( \left( _{H}^{H}\mathcal{YD}\left( \mathcal{M}\right)
,\otimes ,\mathbf{1}\right) ,\Psi \right) $ then the \emph{Radford-Majid
bosonization} of $Q$
\begin{equation*}
Q\rtimes H
\end{equation*}%
is the object $Q\otimes H$ endowed with the following braided bialgebra
structure:%
\begin{eqnarray*}
m_{Q\rtimes H} &=&\left( m_{Q}\otimes m_{H}\right) \circ \left( Q\otimes \mu
_{Q}^{l}\otimes H\otimes H\right) \circ \left( Q\otimes H\otimes
c_{H,Q}\otimes H\right) \circ \left( Q\otimes \Delta _{H}\otimes Q\otimes
H\right) , \\
u_{Q\rtimes H} &=&\left( u_{Q}\otimes u_{H}\right) \circ \Delta _{\mathbf{1}%
}, \\
\Delta _{Q\rtimes H} &=&\left( Q\otimes m_{H}\otimes Q\otimes H\right) \circ
\left( Q\otimes H\otimes c_{Q,H}\otimes H\right) \circ \left( Q\otimes \rho
_{Q}^{l}\otimes H\otimes H\right) \circ \left( \Delta _{Q}\otimes \Delta
_{H}\right) , \\
\varepsilon _{Q\rtimes H} &=&m_{\mathbf{1}}\circ \left( \varepsilon
_{Q}\otimes \varepsilon _{H}\right) .
\end{eqnarray*}%
and the above structures of object in $_{H}^{H}\mathcal{M}_{H}^{H}$.
\end{claim}

\begin{lemma}
\label{lem: SmashGr}Let $(H,m_{H},u_{H},\Delta _{H},\varepsilon _{H},S_{H})$
be a braided graded Hopf algebra with invertible antipode $S_{H}$ in $\left(
\left( \mathcal{M},\otimes ,\mathbf{1}\right) ,c\right) $. \newline
Let $\left( \left( Q=\oplus _{n\in
\mathbb{N}
}Q_{n},\mu ^{l},\rho ^{l}\right) ,m,u,\Delta ,\varepsilon \right) $ be a
graded braided bialgebra in the braided monoidal category $\left( \left(
_{H}^{H}\mathcal{YD}\left( \mathcal{M}\right) ,\otimes ,\mathbf{1}\right)
,\Psi \right) $. Then
\begin{equation*}
Q\rtimes H=\oplus _{n\in
\mathbb{N}
}\left( Q_{n}\otimes H\right)
\end{equation*}%
is a graded braided bialgebra in the monoidal category $\left( \mathcal{M}%
,\otimes ,\mathbf{1}\right) $ and%
\begin{eqnarray*}
m_{a,b}^{Q\rtimes H} &=&\left( m_{a,b}\otimes m_{H}\right) \circ \left(
Q_{a}\otimes \mu _{Q_{b}}^{l}\otimes H\otimes H\right) \circ \left(
Q_{a}\otimes H\otimes c_{H,Q_{b}}\otimes H\right) \circ \left( Q_{a}\otimes
\Delta _{H}\otimes Q_{b}\otimes H\right) , \\
\Delta _{a,b}^{Q\rtimes H} &=&\left( Q_{a}\otimes m_{H}\otimes Q_{b}\otimes
H\right) \circ \left( Q_{a}\otimes H\otimes c_{Q_{b},H}\otimes H\right)
\circ \left( Q_{a}\otimes \rho _{Q_{b}}^{l}\otimes H\otimes H\right) \circ
\left( \Delta _{a,b}\otimes \Delta _{H}\right) .
\end{eqnarray*}%
Moreover, for every $a,b\in
\mathbb{N}
$, we have that

1) $m_{a,b}^{Q\rtimes H}$ is an epimorphism whenever $m_{a,b}$ is an
epimorphism

2) $\Delta _{a,b}^{Q\rtimes H}$ is a monomorphism whenever $\Delta _{a,b}$
is a monomorphism.
\end{lemma}

\begin{proof}
Set $B=Q\rtimes H$ and $B_{n}=Q_{n}\otimes H$ for every $n\in
\mathbb{N}
$. Then we have%
\begin{equation*}
i_{n}^{B}=i_{n}^{Q}\otimes H\qquad \text{and}\qquad
p_{n}^{B}=p_{n}^{Q}\otimes H\text{ }
\end{equation*}%
for every $n\in
\mathbb{N}
$.\newline
By Proposition \ref{coro: grAlg}, in order to prove that $Q\rtimes H$ is a
graded algebra and that $m_{a,b}^{Q\rtimes H}$ is the required morphism, it
is enough to prove that (\ref{form: coro grAlg1}) holds for every $a,b\in
\mathbb{N}
.$ By naturality of the braiding, the fact that $i_{b}^{Q}$ is left $H$%
-linear and $Q$ fulfills (\ref{form: coro grAlg1}), we obtain $m_{B}\circ
(i_{a}^{B}\otimes i_{b}^{B})=i_{a+b}^{B}\circ m_{a,b}^{B}.$ Similarly, by
Proposition \ref{coro: grCoalg}, $Q\rtimes H$ is a graded coalgebra and
hence a graded braided bialgebra.\newline
Now, in view of \ref{claim: YDbraiding}, for every $\left( U,\mu
_{U}^{l},\rho _{U}^{l}\right) ,\left( V,\mu _{V}^{l},\rho _{V}^{l}\right)
\in _{H}^{H}\mathcal{YD}\left( \mathcal{M}\right) $ we have that the morphism%
\begin{equation*}
\Psi _{U,V}:=\left( \mu _{V}^{l}\otimes U\right) \circ \left( H\otimes
c_{U,V}\right) \circ \left( \rho _{U}^{l}\otimes V\right)
\end{equation*}%
is invertible. Recall that, by \cite[Example 4.1.2]{BD}, using the left
adjoint action
\begin{equation*}
ad:=m_{H}\circ \left( m_{H}\otimes H\right) \circ \left( H\otimes
c_{H,H}\right) \circ \left( H\otimes S\otimes H\right) \circ \left( \Delta
_{H}\otimes H\right)
\end{equation*}%
for $H$ and the regular coaction, we have that $H_{ad}:=(H,ad,\Delta _{H})$
is an object in $_{H}^{H}\mathcal{YD}\left( \mathcal{M}\right) $. Dually
using the left coadjoint coaction
\begin{equation*}
coad:=\left( m_{H}\otimes H\right) \circ \left( H\otimes S\otimes H\right)
\circ \left( H\otimes c_{H,H}\right) \circ \left( \Delta _{H}\otimes
H\right) \circ \Delta _{H}
\end{equation*}%
for $H$ and the regular action, we have that $H^{coad}:=(H,m_{H},coad)$ is
an object in $_{H}^{H}\mathcal{YD}\left( \mathcal{M}\right) $. Then we have
\begin{equation*}
m_{a,b}^{Q\rtimes H}=\left( m_{a,b}\otimes m_{H}\right) \circ \left(
Q_{a}\otimes \Psi _{Q_{b},H_{ad}}\otimes H\right) \quad \text{and}\quad
\Delta _{a,b}^{Q\rtimes H}=\left( Q_{a}\otimes \Psi _{H^{coad},Q_{b}}\otimes
H\right) \circ \left( \Delta _{a,b}\otimes \Delta _{H}\right) .
\end{equation*}%
Since $\Psi _{Q_{b},H_{ad}}$ is an isomorphism we get that $%
m_{a,b}^{Q\rtimes H}$ is an epimorphism whenever $m_{a,b}$ is an epimorphism
and $\Delta _{a,b}^{Q\rtimes H}$ is a monomorphism whenever $\Delta _{a,b}$
is a monomorphism.
\end{proof}

\begin{theorem}
\label{teo: TypeOneSmash}Let $(H,m_{H},u_{H},\Delta _{H},\varepsilon
_{H},S_{H})$ be a braided graded Hopf algebra with invertible antipode $%
S_{H} $ in a cocomplete and complete abelian coabelian braided monoidal
category $\left( \left( \mathcal{M},\otimes ,\mathbf{1}\right) ,c\right) $
satisfying $AB5$. Assume that the tensor product commutes with direct sums.
\newline
Let $\left( V,\mu _{V}^{l},\rho _{V}^{l}\right) $ be in $_{H}^{H}\mathcal{YD}%
\left( \mathcal{M}\right) $ and let $\mathbf{1}\left[ V\right] $ be the
braided bialgebra of type one associated to $\mathbf{1}$ and $V\ $in $%
_{H}^{H}\mathcal{YD}\left( \mathcal{M}\right) $. Then $V\otimes H$ is an
object in $_{H}^{H}\mathcal{M}_{H}^{H}$ and
\begin{equation*}
\mathbf{1}\left[ V\right] \rtimes H=H\left[ V\otimes H\right]
\end{equation*}%
is the braided bialgebra of type one associated to $H$ and $V\otimes H$ in $%
\mathcal{M}$.
\end{theorem}

\begin{proof}
$V\otimes H$ is an object in $_{H}^{H}\mathcal{M}_{H}^{H}$ as in \ref{cl:
BesYD}. Let $T_{\mathbf{1}}\left( V\right) $ be the tensor algebra of $V$ in
the monoidal category $\left( _{H}^{H}\mathcal{YD}\left( \mathcal{M}\right)
,\otimes ,\mathbf{1}\right) $. Clearly the isomorphisms
\begin{equation*}
g_{H}:=l_{H}^{-1}:H\longrightarrow \mathbf{1}\otimes H\qquad \text{and}%
\qquad g_{V\otimes H}:=\mathrm{Id}_{V\otimes H}:V\otimes H\longrightarrow
V\otimes H
\end{equation*}%
are a braided bialgebra homomorphism and a morphism in ${_{H}^{H}\mathcal{M}%
_{H}^{H}}$ respectively. By Lemma \ref{lem: SmashGr}, $T_{\mathbf{1}}\left(
V\right) \rtimes H=\oplus _{n\in
\mathbb{N}
}\left( V^{\otimes n}\otimes H\right) $ is a graded braided bialgebra in the
monoidal category $\left( \mathcal{M},\otimes ,\mathbf{1}\right) .$ By
Theorem \ref{coro: GRuniv property of tensor bialgebra}, there is a unique
algebra homomorphism $f:T_{H}(V\otimes H)\rightarrow T_{\mathbf{1}}\left(
V\right) \rtimes H$ such that
\begin{equation*}
f\circ i_{0}^{T_{H}(V\otimes H)}=\left( i_{0}^{T_{\mathbf{1}}\left( V\right)
}\otimes H\right) \circ l_{H}^{-1}\text{\qquad and\qquad }f\circ
i_{1}^{T_{H}(V\otimes H)}=i_{1}^{T_{\mathbf{1}}\left( V\right) }\otimes H.
\end{equation*}%
Moreover $f$ is a graded bialgebra homomorphism. By Theorem \ref{teo: phi
epi}, $f$ is an epimorphism. In fact, by Proposition \ref{pro: tensGrad}, $%
m_{a,b}^{T_{\mathbf{1}}\left( V\right) }$ is an epimorphism for every $%
a,b\in
\mathbb{N}
$ and hence, in view of Lemma \ref{lem: SmashGr}, $m_{a,b}^{T_{\mathbf{1}%
}\left( V\right) \rtimes H}$ is an epimorphism too. \newline
The isomorphisms
\begin{equation*}
g_{H}^{c}:=l_{H}:\mathbf{1}\otimes H\longrightarrow H\qquad \text{and}\qquad
g_{V\otimes H}^{c}:=\mathrm{Id}_{V\otimes H}:V\otimes H\longrightarrow
V\otimes H
\end{equation*}%
are a braided bialgebra homomorphism and a morphism in ${_{H}^{H}\mathcal{M}%
_{H}^{H}}$ respectively. By Lemma \ref{lem: SmashGr}, $T_{\mathbf{1}%
}^{c}\left( V\right) \rtimes H=\oplus _{n\in
\mathbb{N}
}\left( V^{\otimes n}\otimes H\right) $ is a graded braided bialgebra in the
monoidal category $\left( \mathcal{M},\otimes ,\mathbf{1}\right) .$ Then, by
Theorem \ref{coro: GRuniv property of cotensor bialgebra}, there is a unique
coalgebra homomorphism $f^{c}:T_{\mathbf{1}}^{c}\left( V\right) \rtimes
H\rightarrow T_{H}^{c}(V\otimes H)$ such that
\begin{equation*}
p_{0}^{T_{H}^{c}(V\otimes H)}\circ f^{c}=l_{H}\circ \left( p_{0}^{T_{\mathbf{%
1}}^{c}\left( V\right) }\otimes H\right) \text{\qquad and\qquad }%
p_{1}^{T_{H}^{c}(V\otimes H)}\circ f^{c}=\left( p_{1}^{T_{\mathbf{1}%
}^{c}\left( V\right) }\otimes H\right) .
\end{equation*}%
Moreover $f^{c}$ is a graded bialgebra homomorphism. By Theorem \ref{teo:
psi mono} , $f^{c}$ is a monomorphism. In fact, by Proposition \ref{pro:
cotGrad}, $\Delta _{a,b}^{T_{\mathbf{1}}^{c}\left( V\right) }$ is a
monomorphism for every $a,b\in
\mathbb{N}
$ and hence, in view of Lemma \ref{lem: SmashGr}, $\Delta _{a,b}^{T_{\mathbf{%
1}}^{c}\left( V\right) \rtimes H}$ is a monomorphism too. Using the same
notations of Theorem \ref{teo: F of type1}, let us check that%
\begin{equation*}
F_{H,V\otimes H}=f^{c}\circ \left( F_{\mathbf{1},V}\rtimes H\right) \circ f.
\end{equation*}%
By Theorem \ref{teo: F of type1} it is enough to check that
\begin{equation*}
f^{c}\circ \left( F_{\mathbf{1},V}\rtimes H\right) \circ f\circ
i_{0}^{T_{H}(V\otimes H)}=i_{0}^{T_{H}^{c}(V\otimes H)}\qquad \text{and}%
\qquad f^{c}\circ \left( F_{\mathbf{1},V}\rtimes H\right) \circ f\circ
i_{1}^{T_{H}(V\otimes H)}=i_{1}^{T_{H}^{c}(V\otimes H)}.
\end{equation*}%
We have%
\begin{eqnarray*}
f^{c}\circ \left( F_{\mathbf{1},V}\rtimes H\right) \circ f\circ
i_{0}^{T_{H}(V\otimes H)} &=&f^{c}\circ \left( F_{\mathbf{1},V}\otimes
H\right) \circ \left( i_{0}^{T_{\mathbf{1}}\left( V\right) }\otimes H\right)
\circ l_{H}^{-1}=f^{c}\circ \left( i_{0}^{T_{\mathbf{1}}^{c}\left( V\right)
}\otimes H\right) \circ l_{H}^{-1} \\
f^{c}\circ \left( F_{\mathbf{1},V}\rtimes H\right) \circ f\circ
i_{1}^{T_{H}(V\otimes H)} &=&f^{c}\circ \left( F_{\mathbf{1},V}\otimes
H\right) \circ \left( i_{1}^{T_{\mathbf{1}}\left( V\right) }\otimes H\right)
=f^{c}\circ \left( i_{1}^{T_{\mathbf{1}}^{c}\left( V\right) }\otimes
H\right) .
\end{eqnarray*}%
Therefore we have to prove that%
\begin{equation}
f^{c}\circ \left( i_{0}^{T_{\mathbf{1}}^{c}\left( V\right) }\otimes H\right)
\circ l_{H}^{-1}=i_{0}^{T_{H}^{c}(V\otimes H)}\qquad \text{and}\qquad
f^{c}\circ \left( i_{1}^{T_{\mathbf{1}}^{c}\left( V\right) }\otimes H\right)
=i_{1}^{T_{H}^{c}(V\otimes H)}.  \label{form: mah}
\end{equation}%
By construction of $f^{c}$ we get%
\begin{eqnarray*}
p_{0}^{T_{H}^{c}(V\otimes H)}\circ f^{c}\circ \left( i_{0}^{T_{\mathbf{1}%
}^{c}\left( V\right) }\otimes H\right) \circ l_{H}^{-1} &=&l_{H}\circ \left(
p_{0}^{T_{\mathbf{1}}^{c}\left( V\right) }i_{0}^{T_{\mathbf{1}}^{c}\left(
V\right) }\otimes H\right) \circ l_{H}^{-1}=\mathrm{Id}%
_{H}=p_{0}^{T_{H}^{c}(V\otimes H)}\circ i_{0}^{T_{H}^{c}(V\otimes H)}, \\
p_{0}^{T_{H}^{c}(V\otimes H)}\circ f^{c}\circ \left( i_{1}^{T_{\mathbf{1}%
}^{c}\left( V\right) }\otimes H\right) &=&l_{H}\circ \left( p_{0}^{T_{%
\mathbf{1}}^{c}\left( V\right) }i_{1}^{T_{\mathbf{1}}^{c}\left( V\right)
}\otimes H\right) =0=p_{0}^{T_{H}^{c}(V\otimes H)}\circ
i_{1}^{T_{H}^{c}(V\otimes H)}, \\
p_{1}^{T_{H}^{c}(V\otimes H)}\circ f^{c}\circ \left( i_{0}^{T_{\mathbf{1}%
}^{c}\left( V\right) }\otimes H\right) \circ l_{H}^{-1} &=&\left( p_{1}^{T_{%
\mathbf{1}}^{c}\left( V\right) }i_{0}^{T_{\mathbf{1}}^{c}\left( V\right)
}\otimes H\right) \circ l_{H}^{-1}=0=p_{1}^{T_{H}^{c}(V\otimes H)}\circ
i_{0}^{T_{H}^{c}(V\otimes H)}, \\
p_{1}^{T_{H}^{c}(V\otimes H)}\circ f^{c}\circ \left( i_{1}^{T_{\mathbf{1}%
}^{c}\left( V\right) }\otimes H\right) &=&\left( p_{1}^{T_{\mathbf{1}%
}^{c}\left( V\right) }i_{1}^{T_{\mathbf{1}}^{c}\left( V\right) }\otimes
H\right) =\mathrm{Id}_{V}=p_{1}^{T_{H}^{c}(V\otimes H)}\circ
i_{1}^{T_{H}^{c}(V\otimes H)},
\end{eqnarray*}%
so that, by Corollary \ref{coro: uniqueness}, we have proved (\ref{form: mah}%
). Thus $F_{H,V\otimes H}=f^{c}\circ \left( F_{\mathbf{1},V}\rtimes H\right)
\circ f. $ Now, since $f$ is an epimorphism and $f^{c}$ is a monomorphism we
get%
\begin{equation*}
H\left[ V\otimes H\right] =\mathrm{Im}\left( F_{H,V\otimes H}\right) =%
\mathrm{Im}\left( f^{c}\circ \left( F_{\mathbf{1},V}\rtimes H\right) \circ
f\right) =\mathrm{Im}\left( F_{\mathbf{1},V}\rtimes H\right) =\mathrm{Im}%
\left( F_{\mathbf{1},V}\right) \rtimes H=\mathbf{1}\left[ V\right] \rtimes H.
\end{equation*}
\end{proof}

\begin{claim}
Let $(H,m_{H},u_{H},\Delta _{H},\varepsilon _{H},S_{H})$ be a braided graded
Hopf algebra with invertible antipode $S_{H}$ in $\left( \left( \mathcal{M}%
,\otimes ,\mathbf{1}\right) ,c\right) $.\newline
We have an equivalence of braided monoidal categories
\begin{equation*}
F:\left( _{H}^{H}\mathcal{M}_{H}^{H},\otimes _{H},H\right) \rightarrow
\left( _{H}^{H}\mathcal{YD}\left( \mathcal{M}\right) ,\otimes ,\mathbf{1}%
\right) ,
\end{equation*}%
given by
\begin{equation*}
F\left( M\right) =M^{co\left( H\right) }=\mathrm{Eq}\left[ \rho
_{M}^{r},\left( M\otimes u_{H}\right) \circ r_{M}^{-1}\right]
\end{equation*}%
for every $(M,\mu _{M}^{r},\mu _{M}^{l},\rho _{M}^{r},\rho _{M}^{l})\in {%
_{H}^{H}\mathcal{M}_{H}^{H}}$ (see \cite[Theorem 3.9.6]{Bes}). The inverse
of $F$ is the functor $F^{-1}$ defined by%
\begin{equation*}
F^{-1}\left( V\right) =\left( V\otimes H,\mu _{V\otimes H}^{l},\mu
_{V\otimes H}^{r},\rho _{V\otimes H}^{l},\rho _{V\otimes H}^{r}\right)
\end{equation*}%
as in \ref{cl: BesYD} for every $V\in _{H}^{H}\mathcal{YD}\left( \mathcal{M}%
\right) $.
\end{claim}

\begin{theorem}
Let $(H,m_{H},u_{H},\Delta _{H},\varepsilon _{H},S_{H})$ be a braided graded
Hopf algebra with invertible antipode $S_{H}$ in a cocomplete and complete
abelian coabelian braided monoidal category $\left( \left( \mathcal{M}%
,\otimes ,\mathbf{1}\right) ,c\right) $ satisfying $AB5$. Assume that the
tensor product commutes with direct sums. \newline
Let $(M,\mu _{M}^{r},\mu _{M}^{l},\rho _{M}^{r},\rho _{M}^{l})$ be in $%
_{H}^{H}\mathcal{M}_{H}^{H}$ and let $\mathbf{1}\left[ M^{co\left( H\right) }%
\right] $ be the braided bialgebra of type one associated to $\mathbf{1}$
and $M^{co\left( H\right) }\ $in $_{H}^{H}\mathcal{YD}\left( \mathcal{M}%
\right) $. Then
\begin{equation*}
\mathbf{1}\left[ M^{co\left( H\right) }\right] \rtimes H=H[M]
\end{equation*}%
is the braided bialgebra of type one associated to $H$ and $M$ in $\mathcal{M%
}$.
\end{theorem}

\begin{proof}
By Theorem \ref{teo: TypeOneSmash} $\mathbf{1}\left[ M^{coH}\right] \rtimes
H=H\left[ M^{coH}\otimes H\right] \simeq H[M] $ is the braided bialgebra of
type one associated to $H$ and $M$ in $\mathcal{M}$.
\end{proof}


\begin{thebibliography}{}
\bibitem[AS]{AS} N. Andruskiewitsch and H-J. Schneider, \emph{Pointed Hopf
algebras}. New directions in Hopf algebras, 1--68, Math. Sci. Res. Inst.
Publ., \textbf{43}, Cambridge Univ. Press, Cambridge, 2002.

\bibitem[AMS1]{AMS:Cotensor} A. Ardizzoni, C. Menini and D. \c{S}tefan,
\emph{Cotensor Coalgebras in Monoidal Categories}, Comm. Algebra., to
appear. (arXiv:math.QA/0507334)

\bibitem[AMS2]{AMS} A. Ardizzoni, C. Menini and D. \c{S}tefan, \emph{%
Hochschild Cohomology And 'Smoothness' In Monoidal Categories}, J. Pure
Appl. Algebra, \textbf{208} (2007), 297-330.

\bibitem[Bes]{Bes} Yu. N. Bespalov, \emph{Crossed modules and quantum groups
in braided categories}, Appl. Categ. Structures \textbf{5} (1997), 155--204.

\bibitem[BD]{BD} Y. Bespalov, B. Drabant, \emph{Hopf (bi-)modules and
crossed modules in braided monoidal categories}, J. Pure Appl. Algebra
\textbf{123} (1998), no. 1-3, 105--129.

\bibitem[Ka]{Kassel} Kassel, \emph{Quantum Groups, }Graduate Text in
Mathematics \textbf{155},\textbf{\ }Springer, 1995.

\bibitem[Mj1]{Majid} S. Majid, \emph{Foundations of quantum group theory},
Cambridge University Press, 1995.

\bibitem[NT]{NT} C. N\u{a}st\u{a}sescu, B. Torrecillas, \emph{Graded
coalgebras}, Tsukuba J. Math. \textbf{17} (1993), 461-479.

\bibitem[Ni]{Ni} Nichols, W. D. Bialgebras of type one. Comm. Algebra 6
(1978), no. 15, 1521--1552.

\bibitem[Po]{Po} N. Popescu, \emph{Abelian Categories with Application to
Rings and Modules}, Academic Press, London \& New York, (1973).

\bibitem[Ro]{Ro} M. Rosso, \emph{Quantum groups and quantum shuffles},
Invent. Math. \textbf{133} (1998), no. 2, 399--416.

\bibitem[Sch]{Schauenburg1} P. Schauenburg, \emph{A characterization of the
Borel-like subalgebras of quantum enveloping algebras} Comm. Algebra \textbf{%
24} (1996), no. 9, 2811--2823.
\end{thebibliography}
\end{document}